\documentclass[11pt]{amsart}
\usepackage{amstext,amssymb,amsmath,amsbsy}
\usepackage[left=2.6cm,right=2.6cm,top=2.7cm,bottom=2.7cm]{geometry}

\usepackage{tikz}
\usepackage{hyperref}
\usepackage{amscd}
\usepackage{amsfonts,dsfont}
\usepackage{indentfirst}
\usepackage{verbatim}
\usepackage{amsmath}
\usepackage{amsthm}
\usepackage{enumerate}
\usepackage{graphicx}
\usepackage{subfig}
\usepackage{cleveref}
\usepackage{xcolor}

\usepackage[OT1]{fontenc}
\usepackage[latin1]{inputenc}
\usepackage[english]{babel}
\usepackage{amssymb}
\newtheorem{theorem}{Theorem}[section]
\newtheorem{lemma}{Lemma}[section]

\newtheorem{proposition}{Proposition}[section]

\newtheorem{remark}{Remark}[section]

\newtheorem{definition}{Definition}[section]
\newcommand{\cF}{{\mathcal F}}
\newcommand{\cB}{{\mathcal B}}

\newcommand{\cT}{{\mathcal T}}

\newcommand{\tw}{\widetilde w}

\newcommand{\hw}{\hat w}

\newcommand{\tr}{^\mathsf{T}}

\numberwithin{equation}{section}

\newcommand{\dsp}{\displaystyle}

\newcommand{\cH}{{\mathcal H}}

\newcommand{\supp}{\operatorname{supp}}

\newcommand{\eps}{\varepsilon}

\newcommand{\mN}{\mathbb{N}}

\newcommand{\mR}{\mathbb{R}}

\newcommand{\mC}{\mathbb{C}}

\newcommand{\mc}{\mathrm{c}}

\newcommand{\hI}{ I}

\newcommand{\bC}{{\bf C}}

\newcommand{\bw}{{\bf w}}

\newcommand{\cX}{\mathcal X}

\newcommand{\cK}{{\mathcal K}}

\newcommand{\ta}{\widetilde \alpha}

\newcommand{\tb}{\widetilde \beta}

\newcommand{\ha}{\hat \alpha}

\newcommand{\hb}{\hat \beta}

\newcommand{\cL}{\mathcal L}

\newcommand{\cY}{{\mathcal Y}}

\newcommand{\sign}{\mathrm{sign }}
\date{03/03/2021}

\title[Controllability of time-dependent linear hyperbolic systems]
{On the optimal controllability time for linear hyperbolic systems with time-dependent coefficients}%\author{Jean-Michel Coron \thanks{Sorbonne Universit\'{e}, Universit\'{e} Paris-Diderot SPC, CNRS, INRIA, Laboratoire Jacques-Louis Lions, \'{e}quipe Cage, Paris, France, coron@ann.jussieu.fr.} \and Hoai-Minh Nguyen \thanks{Ecole Polytechnique F\'ed\'erale de Lausanne, EPFL,  SB, CAMA, Station 8,  CH-1015 Lausanne, Switzerland,  hoai-minh.nguyen@epfl.ch. }}

\author[J.-M. Coron]{Jean-Michel Coron}

\author[H.-M. Nguyen]{Hoai-Minh Nguyen}

\address[J.-M. Coron]{Sorbonne Universit\'{e} \newline\indent
	Universit\'{e} Paris-Diderot SPC, CNRS, INRIA, \newline\indent
	Laboratoire Jacques-Louis Lions, \'{e}quipe Cage, Paris, France}
\email{coron@ann.jussieu.fr}

\address[H.-M. Nguyen]{Department of Mathematics \newline\indent
	EPFL SB CAMA \newline\indent
	Station 8 CH-1015 Lausanne, Switzerland}
\email{hoai-minh.nguyen@epfl.ch}

\begin{document}
\maketitle

\begin{abstract} The optimal time for the controllability of linear hyperbolic systems in one dimensional space with one-side controls has been obtained recently  for time-independent coefficients in our previous works.
In this paper, we consider linear hyperbolic systems with time-varying zero-order terms. We show the possibility that  the optimal time  for the null-controllability becomes significantly larger than the one of the time-invariant setting even when the zero-order term is indefinitely differentiable.  When the analyticity with respect to time is  imposed for the zero-order term, we also establish that the optimal time is the same as in the time-independent setting.
\end{abstract}



\tableofcontents

\section{Introduction and statement of the main results}

Hyperbolic systems in one dimensional space are frequently used
in the modeling of many systems such as traffic flow \cite{SFA08}, heat exchangers \cite{XS02},  fluids in open channels \cite{GL03, dHPCAB03, GLS04, DP08}, and phase transition \cite{Goatin06}. Many other interesting examples can be found in \cite{BC16} and the references therein.  The optimal time for the controllability of hyperbolic systems in one dimensional space with one-side controls has been derived recently  for time-independent coefficients \cite{CoronNg19, CoronNg19-2}.
In this paper, we consider hyperbolic systems with time-varying zero-order terms. It is known that these systems are controllable in some positive time. In this paper,  we show the possibility that  the optimal time  for the null-controllability becomes significantly larger than the one of the time-invariant setting even when the zero-order term is indefinitely differentiable.  When the analyticity with respect to time is  imposed for the zero-order term, we also establish that the optimal time is the same as in the time-independent setting.   The first result is quite surprising since the zero-order term does not interfere with the characteristic flows of the system. The later result complement to the first one can be then viewed as an extension of a well-known controllability property of linear differential equations: if a linear control system is controllable in some positive time and  is {\it analytic}, then it is controllable in any time greater than the optimal time, which is 0.

Let us first briefly discuss known results  for the time-independent coefficients to underline the phenomena. Consider the system
\begin{equation}\label{Sys-1}
\partial_t u (t, x) =  \Sigma(x) \partial_x u (t, x) + C(x) u(t, x) \mbox{ for } (t, x)  \in \mR_+ \times (0, 1).
\end{equation}
Here $u = (u_1, \cdots, u_n)\tr : \mR_+ \times (0, 1) \to \mR^n$ ($n \ge 2$),   $\Sigma$ and $C$ are   $(n \times n)$ real, matrix-valued functions defined in $[0,1]$. We assume that, for every $x \in [0, 1]$,  the matrix $\Sigma(x)$ is diagonalizable with $m \ge 1$  distinct positive eigenvalues and $k = n - m \ge 1$  distinct negative eigenvalues. Using Riemann coordinates, one might assume that $\Sigma(x)$ is of the form
\begin{equation}\label{form-A}
\Sigma(x) = \mbox{diag} \big(- \lambda_1(x), \cdots, - \lambda_{k}(x),  \lambda_{k+1}(x), \cdots,  \lambda_{n}(x) \big),
\end{equation}
where
\begin{equation}\label{relation-lambda}
-\lambda_{1}(x) < \cdots <  - \lambda_{k} (x)< 0 < \lambda_{k+1}(x) < \cdots < \lambda_{k+m}(x).
\end{equation}
In what follows, we assume that
\begin{equation}\label{cond-lambda}
\mbox{$\lambda_i$  is of class $C^2$ on $[0, 1]$  for $1 \le i \le n \,  (= k + m)$,}
\end{equation}
and denote
$$
\mbox{$u_- = (u_1, \cdots, u_k)\tr $ and $u_+ = (u_{k+1}, \cdots, u_{k+m})\tr$.}
$$

We are interested in the following type of boundary conditions and boundary controls. The boundary conditions at $x = 0$ are given by
\begin{equation}\label{bdry-0}
u_- (t, 0) = B  u_+ (t, 0) \mbox{ for } t \ge 0,
\end{equation}
for some $(k \times m)$ real {\it constant} matrix $B$, and   at $x = 1$
\begin{equation}\label{bdry-1}
u_+(t, 1) \mbox{ is controlled } \mbox{ for } t \ge 0.
\end{equation}

 Let us recall that the control system \eqref{Sys-1}, \eqref{bdry-0}, and \eqref{bdry-1} is null-controllable (resp. exactly controllable) at time $T>0$ if, for every initial datum $u_0: (0,1)\to \mathbb{R}^n$ in $[L^2(0,1)]^n$ (resp. for every initial datum $u_0: (0,1 )\to \mathbb{R}^n$ in $[L^2(0,1)]^n$  and for every (final) state $u_T:  (0,1 )\to \mathbb{R}^n$  in $[L^2(0,1)]^n$), there is a control $U:(0,T)\to \mathbb{R}^m$ in $[L^2(0,T)]^m$ such that the solution
of \eqref{Sys-1}, \eqref{bdry-0}, and \eqref{bdry-1} (with $u_+ = U$) satisfying $u(t=0,x)=u_0(x)$ vanishes (resp. reaches $u_T$) at the time $T$: $u(t=T, \cdot)=0$ (resp. $u(t = T, \cdot)= u_T$).

Throughout this paper, we consider broad solutions in $L^2$ with respect to $t$ and $x$ for an initial datum in $[L^2(0, 1)]^n$ and a control in $[L^2(0, T)]^m$ (see, for example,
\cite[Section 3]{1981-Li-Wen-Chi-Shen-CAM}). In particular,  the solutions belong to  $C([0, T]; [L^2(0, 1)]^n)$ and $C([0, 1]; [L^2(0, T)]^n)$.
The well-posedness for broad solutions for system \eqref{Sys-1}, \eqref{bdry-0}, and \eqref{bdry-1} even when $\Sigma$ and $C$ depending also on $t$ is standard.

Set
\begin{equation}\label{def-tau}
\tau_i :=  \int_{0}^1 \frac{1}{|\lambda_i(\xi)|}  \, d \xi  \mbox{ for } 1 \le i \le n.
\end{equation}

The exact controllability, the null-controllability, and the boundary stabilization problem of  hyperbolic system in one dimensional space  have been widely investigated in the literature for almost half a century,  see,  e.g.,  \cite{BC16} and the references therein.
Concerning the exact controllability and  the null-controllability related to \eqref{Sys-1}, \eqref{bdry-0} and \eqref{bdry-1},
the pioneer works date back to the ones of Rauch and Taylor \cite{RT74} and Russell \cite{Russell78}.   In particular, it was shown, see \cite[Theorem 3.2]{Russell78},  that system \eqref{Sys-1}, \eqref{bdry-0}, and \eqref{bdry-1} is null-controllable for  time $
\tau_k + \tau_{k+1}$,
and is exactly controllable at the same time if $k=m$ and $B$ is invertible. The extension of  this result for quasilinear systems was initiated by
Greenberg and Li \cite{GL84} and Slemrod \cite{Slemrod83}.

A recent efficient way in the study of the  stabilisation and the  controllability of system \eqref{Sys-1}, \eqref{bdry-0}, and \eqref{bdry-1} is via a backstepping approach.  The backstepping approach for the control of partial differential equations was pioneered by  Miroslav Krstic and his coauthors (see \cite{Krstic08} for a concise introduction).  The backstepping method is now frequently used for various control problems, modeling by partial differential equations in one dimensional space. For example, it has been used to stabilize the  wave equations  \cite{KGBS08, SK09, SCK10}, the parabolic equations in \cite{SK04,SK05},  nonlinear parabolic equations  \cite{Vazquez08},  and to  obtain the null-controllability of the heat equation \cite{CoronNg17}. The standard backstepping approach relies on the Volterra transform of the second kind. It is worth noting that, in some situations, more general transformations have to be  considered as for Korteweg-de Vries equations  \cite{CoronC13}, Kuramoto--Sivashinsky equations \cite{Coron15}, Schr\"{o}dinger's equation \cite{CoronGM18}, and hyperbolic equations with internal controls \cite{2019-Zhang-preprint}.

The use of backstepping approach for the hyperbolic system in one dimensional space  was first proposed by Coron et al. \cite{CVKB13} for $2\times 2$ system $(m=k=1)$.  Later,  this approach has been extended and now can be applied for general pairs $(m,k)$, see \cite{MVK13, HMVK16, AM16, CHO17, CoronNg19, CoronNg19-2,HO19}.

Set
\begin{equation}\label{def-Top}
T_{opt} := \left\{ \begin{array}{cl}  \dsp \max \big\{ \tau_1 + \tau_{m+1}, \dots, \tau_k + \tau_{m+k}, \tau_{k+1} \big\} & \mbox{ if } m \ge k, \\[6pt]
\dsp \max \big\{ \tau_{k+1-m} + \tau_{k+1},  \tau_{k+2-m} + \tau_{k+2},  \dots, \tau_{k} + \tau_{k+m} \big\} &  \mbox{ if } m < k.
\end{array} \right.
\end{equation}
Involving the backstepping technique, we established \cite{CoronNg19, CoronNg19-2}
that  the null-controllability holds at $T_{opt}$ for generic $B$ and $C$,  and the null-controllability holds for any $T> T_{opt}$ under the condition  $B \in \cB$. Here
\begin{equation}\label{def-B}
{\mathcal B}: = \Big\{B \in \mR^{k \times m}; \mbox{ such that  \eqref{cond-B-1} holds for  $1 \le i \le  \min\{k, m-1\}$} \Big\},
\end{equation}
where
\begin{multline}\label{cond-B-1}
\mbox{ the $i \times i$  matrix formed from the last $i$ columns and the last $i$ rows of $B$  is invertible.}
\end{multline}
Roughly speaking, the condition $B \in \cB$ allows us to implement $l$ controls corresponding to the fastest  positive speeds to control  $l$ components corresponding to the lowest negative speeds \footnote{The $i$  direction ($1 \le i \le n$) is called positive (resp. negative) if $\Sigma_{ii}$ is positive (resp. negative).}. It is clear that $B \in \cB$ for almost every $k \times m$ matrix  $B$.  It is worthy noting that the condition $T > T_{opt}$ is necessary, see \cite[Assertion 2) of Theorem 1.1]{CoronNg19}.  The optimality of $T_{opt}$ was established under the additional condition \eqref{cond-B-1} being valid with $i = m$  when $k \ge m$, see \cite[Proposition 1.6]{CoronNg19}.
Our results improved the time to reach the null-controllability obtained previously.
Similar conclusions hold for the exact controllability under the natural conditions $m \ge k$ and  \eqref{cond-B-1} for $1 \le i \le k$ (see \cite{CoronNg19, CoronNg19-2,HO19}).
When the system is homogeneous, i.e., $C \equiv 0$,  we established  that the null-controllability can be achieved via a time-independent feedback  even for the  quasilinear setting  \cite{CoronNg20}.  We also constructed Lyapunov functions which yield the null-controllability for such a system at the optimal time $T_{opt}$ \cite{CoronNg20-L}.

\medskip

In this paper, we are interested in hyperbolic systems with time-dependent coefficients in one dimensional space. More precisely, instead of
\eqref{Sys-1}, \eqref{bdry-0}, and \eqref{bdry-1}, we deal with
\begin{equation}\label{Sys-1-t}
\partial_t u (t, x) =  \Sigma(x) \partial_x u (t, x) + C(t, x) u(t, x) \mbox{ for } (t, x)  \in \mR_+ \times (0, 1),
\end{equation}
and \eqref{bdry-0}, and \eqref{bdry-1}.

The first result of the paper reveals that  the optimal time  for the null-controllability of system \eqref{Sys-1-t}, \eqref{bdry-0}, and \eqref{bdry-1} might be significantly larger than the one for the time-independent setting even when $\Sigma$ is constant and $C$ is indefinitely differentiable.  More precisely, we have

\begin{theorem}\label{thm1} Let $k \ge 1$,  $m \ge 2$,  and  $\Sigma$ be constant such that
\eqref{relation-lambda} holds.  Assume that
\begin{equation}\label{assumption-B}
B_{k, 1} \neq 0,  \quad B_{k, \ell} \neq 0, \quad B_{k, j} = 0 \mbox{ for } 2 \le j \le m \mbox{ with } j \neq \ell,
\end{equation}
for some $2 \le \ell \le m$. There exists  $C \in C^\infty([0, + \infty) \times [0, 1])$ such that for all
$\eps > 0$, system \eqref{Sys-1-t}, \eqref{bdry-0}, and \eqref{bdry-1} is {\bf not} null-controllable at time
\begin{equation}\label{def-T}
T = \tau_{k} + \tau_{k+1} - \eps.
\end{equation}
\end{theorem}

\begin{remark} \rm The definition of the null-controllability for system \eqref{Sys-1-t}, \eqref{bdry-0}, and \eqref{bdry-1} is similar to the one corresponding to \eqref{Sys-1}, \eqref{bdry-0}, and \eqref{bdry-1}.
\end{remark}

\begin{remark} \rm There are infinitely many matrices $B \in \cB$ satisfying \eqref{assumption-B}.
\end{remark}

In a recent work, Coron et al. \cite{CHOS20} establish the null-controllability of \eqref{Sys-1-t}, \eqref{bdry-0}, and \eqref{bdry-1}  for time $\tau_k  + \tau_{k+1}$ for all $k \times m$ matrices $B$. They also obtain
stabilizing  feedbacks and derive similar results when $\Sigma$ depends on $t$.  Combining \Cref{thm1} and their results, one obtains the optimality for the time $\tau_{k} + \tau_{k+1}$ when $m \ge 2$ and $k \ge 1$, and for a large class of $B$.

The proof of \Cref{thm1} is based on constructing counter-examples for the associated observability inequality. The construction is inspired by the one given in the proof of \cite[Assertion 2) of Theorem 1.1]{CoronNg19} but much more involved.

\medskip
When the analyticity of $C$ with respect to time is imposed, the situation changes dramatically. To state our results in this direction, we first introduce some notations.  For a non-empty interval  $(a, b)$ of $\mR$ and a Banach space $\cX$, we denote
$$
\cH\big( (a, b); \cX \big) = \Big\{\Phi: (a, b) \to \cX; \; \Phi \mbox{ is analytic} \Big\}.
$$
When the space $\cX$ is clear, we simply call a $\Phi  \in \cH\big( (a, b); \cX \big) $ that $\Phi$ is analytic in $(a, b)$. For $m \ge k$, set
\begin{equation}\label{def-Be}
{\mathcal B}_e: = \Big\{B \in \mR^{k \times m}; \mbox{ such that  \eqref{cond-B-1} holds for  $1 \le i \le k $} \Big\}.
\end{equation}
Denote
\begin{equation}\label{def-T1}
T_1 =  \tau_{k} + \tau_{k+1}.
\end{equation}

Our main results for the analytic setting are the following two theorems. The first one on  the null-controllability is:

\begin{theorem}\label{thm-A} Let $k \ge m \ge 1$, and let $B \in \cB$ be such that \eqref{cond-B-1} holds for $i =m$. Assume that $C \in \cH\big(I; [L^\infty(0, 1)]^{n \times n}\big)$ for some open interval $I$ containing  $[0, T_1]$.  System  \eqref{Sys-1-t}, \eqref{bdry-0}, and \eqref{bdry-1} starting from time $0$ is null-controllable at any time $T> T_{opt}$.
\end{theorem}

The second one on  the exact-controllability is:

\begin{theorem}\label{thm-E} Let $m \ge k \ge 1$, and let $B \in \cB_e$.  Assume, for some open interval $I$ containing  $[T_{opt} - T_1, T_{opt}]$, that $C \in \cH\big(I; [L^\infty(0, 1)]^{n \times n} \big)$
System  \eqref{Sys-1-t}, \eqref{bdry-0}, and \eqref{bdry-1} starting from time $0$ is exact-controllable at any time $T> T_{opt}$.
\end{theorem}

Except for the case where $m = 1$  for which  $T_1 = T_{opt}$,  \Cref{thm-A,thm-E} are new to our knowledge. \Cref{thm1,thm-A,thm-E} reveal the crucial role of the analytic assumption of the coefficients  on the optimal controllability time. It is  well-known  that a linear control system modeled by  differential equations is controllable in some time $T$ and is {\it analytic},  then it is controllable in any time greater than the optimal time, which is 0, see, e.g., \cite[Chapter 1]{Coron07} or \cite[Chapter 3]{Sontag98}.  \Cref{thm-A,thm-E}, which are complement to \Cref{thm1},   can be thus viewed as an extension of this well-known result for linear hyperbolic systems in one dimensional space.

A related context to  \Cref{thm-E} is the one of the wave equation. For the wave equation with time varying,  first and  zero-order terms being analytic in time, it is known that  the
controllability holds under  a sharp geometric control condition,  introduced in \cite{BLR92} (see also \cite{RT74}). This can be obtained by combining the results in \cite{BLR92}, on  the propagation of singularities for the wave equation,  and the unique continuation principle for the wave equations with coefficients analytic in time using Carleman's estimates due to
Tataru-H\"ormander-Robbiano-Zuily \cite{Tataru95,Hor97,RZ98} (see also \cite{LL19} for a discussion). Related results concerning the Schr\"odinger equation are due to 
Nalini Anantharaman, Matthieu L\'{e}autaud, and Fabricio Maci\`a  \cite{ALM16}.

We now say a few words on the proof of \Cref{thm-A,thm-E}.   \Cref{thm-E} is derived from \Cref{thm-A}  using  our arguments in the proof of  \cite[Theorem 3]{CoronNg19-2}.
The proof of \Cref{thm-A} is inspired from the analysis  in \cite{CoronNg19-2}, in which we established similar result for the time-independent setting. The crucial part of the analysis   is then to locate  the essential, analytic nature of the system, the smoothness is not sufficient as shown previously in \Cref{thm1}. This is done by exploring both the orignal system and its dual one.  The proof also  involves  the theory of perturbations of analytic compact operators, see, e.g., \cite{Kato}. As a consequence of our analysis, we also obtain the unique continuation principle for hyperbolic systems for the optimal time in the analytic setting (see \Cref{pro-UCP}), which has its own interest. The strategy of the proof is described in more details at the beginning of \Cref{sect-thmA}.

\medskip
The paper is organized as follows.  \Cref{thm1,thm-A,thm-E} are given in \Cref{sect-thm1,sect-thmA,sect-thmE}, respectively. In the appendix, we establish properties of  hyperbolic systems used in \Cref{sect-thmA}. The situation is non-standard in the sense that the domains considered are not rectangle and the boundary conditions are involved.  The analysis is  delicate and has its own interest.

\section{Analysis in the smooth setting - Proof of \Cref{thm1}}\label{sect-thm1}

The starting point of the proof of \Cref{thm1} is the  equivalence between the null-controllability of  system  \eqref{Sys-1-t}, \eqref{bdry-0}, and \eqref{bdry-1} and its corresponding observability inequality. To this end,
we first introduce some notations and recall this property.

Fix $T>0$ and  define
$$
\begin{array}{cccc}
\cF_{T}: &  [L^2(0, T)]^m &  \to & [L^2(0, 1)]^n \\[6pt]
& \cF_T (U) &  \mapsto & u(T, \cdot),
\end{array}
$$
where $u$ is the unique solution of system \eqref{Sys-1-t}, \eqref{bdry-0}, and \eqref{bdry-1} with $u_+(\cdot, 1) = U$  and with $u(0, \cdot) = 0$.  Denote
$$
\Sigma_- = \mbox{diag} (- \lambda_1, \cdots, - \lambda_k) \quad \mbox{ and } \quad \Sigma_+ = \mbox{diag} (\lambda_{k+1}, \cdots, \lambda_{k+m}).
$$

As usual,  we have

\begin{lemma}\label{lem-FT} Let $T>0$. We have, for $\varphi \in [L^2(0, 1)]^n$,
\begin{equation*}
\cF_{T}^*(\varphi) = \Sigma_+ (1) v_+(\cdot, 1) \mbox{ in } (0, T),
\end{equation*}
where $v$ is the unique solution of the system
\begin{equation}\label{eq-v}
\partial_t v (t, x) =  \Sigma(x) \partial_x v (t, x) + \big(\Sigma'(x) - C\tr(t, x) \big)  v(t, x)
\mbox{ for } (t, x)  \in (0, T) \times (0, 1),
\end{equation}
with, for $0 < t < T$,
\begin{equation}\label{bdry-v1}
v_-(t, 1)  = 0,
\end{equation}
\begin{equation}\label{bdry-v0}
\Sigma_+ (0) v_+(t, 0) = - B\tr \Sigma_- (0) v_- (t, 0),
\end{equation}
and
\begin{equation}\label{initial-v}
v(t = T, \cdot) = \varphi \mbox{ in } (0, 1).
\end{equation}
\end{lemma}

The proof of \Cref{lem-FT} is standard and omitted, see, e.g., \cite[the proof of Lemma 1]{CoronNg19-2} for a closely related context.

From \Cref{lem-FT}, one derives the following characterization of the null-controllability of system
\eqref{Sys-1-t}, \eqref{bdry-0}, and \eqref{bdry-1} in  time $T$, whose proof is standard and omitted, see, e.g., \cite[Section 2.3]{Coron07}.

\begin{lemma}\label{lem2} Let $T>0$.  System \eqref{Sys-1-t}, \eqref{bdry-0}, and \eqref{bdry-1} starting at time 0 is null-controllable in  time $T$ if and only if there exists a positive constant $C_T$ such that
\begin{equation*}
\| v_+ (\cdot, 1)\|_{L^2(0, T)} \ge C_T \|  v(0, \cdot)\|_{L^2(0, 1)},
\end{equation*}
for all solutions $v$  of system \eqref{eq-v}, \eqref{bdry-v1}, and \eqref{bdry-v0}.
\end{lemma}

%\begin{remark} \rm Given $v(T, \cdot) \in L^2(0, 1)$, there exists a unique broad solution $v \in C^0\big( [0, T]; L^2(0, 1)  \big) \cap C^0\big( [0, 1]; L^2(0, T)  \big)$ of the system \eqref{eq-v}, \eqref{bdry-v1}, and \eqref{bdry-v0}; see \cite{CoronNg19, CHOS20}.
%\end{remark}

We are ready to give

\begin{proof}[Proof of \Cref{thm1}] In what follows, we will assume that
$$
T \ge \max\big\{ T_{opt}, \tau_k + \tau_{k+\ell}\big\},
$$
where $T_{opt}$ is defined by \eqref{def-Top}; hence $\eps$ is assumed to be sufficiently small (note that $\tau_{k} + \tau_{k+1} > \max\big\{ T_{opt}, \tau_k + \tau_{k+\ell}\big\}$  since $2 \le \ell \le m$).
We will consider the coefficient $C(t, x)$ satisfying the following structure:
\begin{equation}\label{cond-C}
C_{i, j} (t, x) =  \left\{ \begin{array}{cl}
- \alpha (t, x) & \mbox{ if } (i, j) = (k, k + \ell), \\[6pt]
- \beta (t, x) & \mbox{ if } (i, j) = (k + 1, k+ \ell), \\[6pt]
0 &  \mbox{ otherwise},
\end{array}\right.
\end{equation}
where $\alpha$ and $\beta$ are two smooth functions defined later.

Since $\Sigma$ is constant and $C$ satisfies \eqref{cond-C},  system \eqref{eq-v} is equivalent to, for $(t, x) \in (0, T) \times (0, 1)$,
\begin{equation}\label{eq-vj}
\partial_t v_{j} (t, x) = \Sigma_{j, j} \partial_x v_{j} (t, x) \mbox{ if } 1 \le j \le n \mbox{ with }
j \neq k + \ell,
\end{equation}
and
\begin{equation}\label{eq-v-k+2}
\partial_t v_{k+\ell} (t, x) = \lambda_{k+\ell} \partial_x v_{k+\ell} (t, x) + \alpha (t, x) v_{k} (t, x) + \beta (t, x) v_{k+1} (t, x)
\end{equation}
($\Sigma_{j, j} = -\lambda_j$ if $1 \le j \le k$ and $\Sigma_{j, j} = \lambda_j$ otherwise).

Under appropriate choices of  $\alpha$ and $\beta$ determined later, we will construct a smooth solution $v$ of system \eqref{eq-vj} and \eqref{eq-v-k+2} for which, for $t \in (0, T)$,
\begin{equation}\label{bdry-vv1}
v_-(t, 1)  = 0,
\end{equation}
\begin{equation}\label{bdry-vv0}
\Sigma_+ (0) v_+(t, 0) = - B\tr \Sigma_- (0) v_- (t, 0),
\end{equation}
and  $v$ satisfies the following {\it additional} conditions:
\begin{equation}\label{main-vv}
v_+(\cdot, 1) =0 \quad \mbox{ and } \quad v(0, \cdot) \not \equiv 0.
\end{equation}
By \Cref{lem2}, the conclusion of \Cref{thm1} follows from this construction.

\medskip
We now construct $\alpha$ and $\beta$. To this end, we first derive their constrains.  From \eqref{eq-v-k+2}, we have,  for $\tau_{k+\ell} \le t + \tau_{k+\ell} \le T$ and $0 \le s \le 1$,
\begin{multline*}
\frac{d}{ds} \Big( v_{k+\ell} (t + \tau_{k+\ell} s, 1 - s)  \Big) = \tau_{k+\ell} \alpha(t + \tau_{k+\ell} s, 1- s) v_k (t + \tau_{k+\ell} s, 1- s)  \\[6pt]
+ \tau_{k+\ell} \beta(t + \tau_{k+\ell}s, 1-s) v_{k+1} (t + \tau_{k+\ell} s, 1- s).
\end{multline*}
This implies, for $\tau_{k+\ell} \le t + \tau_{k+\ell} \le T$,
\begin{multline}\label{cond-v4-*}
v_{k+\ell}(t + \tau_{k+\ell}, 0) = \int_0^1 \tau_{k+\ell} \alpha(t + \tau_{k+\ell} s, 1- s) v_k (t + \tau_{k+\ell} s, 1- s)  \,  ds \\[6pt]
+ \int_0^1 \tau_{k+\ell} \beta(t + \tau_{k+\ell}s, 1-s) v_{k+1} (t + \tau_{k+\ell} s, 1- s) \, ds + v_{k+\ell} (t, 1).
\end{multline}
It follows that, {\it if} $v_{k+\ell} (t, 1) = 0$ for $\tau_{k+\ell} \le t + \tau_{k+\ell} \le T$, then
\begin{multline}\label{cond-v4}
v_{k+\ell}(t + \tau_{k+\ell}, 0) = \int_0^1 \tau_{k+\ell} \alpha(t + \tau_{k+\ell} s, 1- s) v_k (t + \tau_{k+\ell} s, 1- s)  \, ds \\[6pt]
+ \int_0^1 \tau_{k+\ell} \beta(t + \tau_{k+\ell}s, 1-s) v_{k+1} (t + \tau_{k+\ell} s, 1- s) \, ds
\end{multline}
for $\tau_{k+\ell} \le t + \tau_{k+\ell} \le T$.

We will assume that, for $t \in (0, T)$,
\begin{equation}\label{choice-a-1}
 \tau_{k+\ell} \alpha(t + \tau_{k+\ell} s, 1- s) =   \ta (t + \tau_{k+\ell})  \mbox{ for } s \in [0, 1],
\end{equation}
and
\begin{equation}\label{choice-b-1}
 \tau_{k+\ell} \beta(t + \tau_{k+\ell} s, 1- s)  = \tb (t + \tau_{k+\ell}) \mbox{ for } s \in [0, 1],
\end{equation}
for some functions $\ta$ and $\tb$ constructed later;  this implies that the LHS of \eqref{choice-a-1} and \eqref{choice-b-1} are constant with respect to  $s \in [0, 1]$. Given $\ta$ and $\tb$ defined in $\mR$, one can verify that \eqref{choice-a-1} and \eqref{choice-b-1} hold if
\begin{equation}\label{alpha-ta}
\alpha(t, x) =  \tau_{k+\ell}^{-1} \ta (t + \tau_{k+\ell} x) \quad \mbox{ and } \quad \beta(t, x) =  \tau_{k+\ell}^{-1}  \tb (t + \tau_{k+\ell} x).
\end{equation}

Under conditions \eqref{choice-a-1} and \eqref{choice-b-1}, by  replacing first $s$ by $1-s$ and then $t + \tau_{k+\ell}$ by $t$, identity \eqref{cond-v4} can be then written as, for $t \in (\tau_{k+\ell}, T)$,
\begin{equation}\label{eq-v4}
v_{k+\ell}(t , 0) = \ta(t)  \int_0^1 v_{k} (-\tau_{k+\ell} s  + t, s) \,   ds
+ \tb(t )  \int_0^1 v_{k+1} (-\tau_{k+\ell} s  + t , s) \, ds.
\end{equation}

\noindent We write \eqref{bdry-vv0} as
\begin{equation}\label{bdry-vv-M}
v_+ (t, 0) =  - \Sigma_+^{-1} B\tr \Sigma_- v_-(t, 0).
\end{equation}

\noindent In what follows, we consider the solution $v$ satisfying
\begin{equation}\label{final-1}
v_{1}(T, \cdot) = \cdots = v_{k-1} (T, \cdot) = v_{k+1} (T, \cdot) = \cdots = v_{k+ m} (T, \cdot) = 0,
\end{equation}
and
\begin{equation}\label{final-2}
v_k(T, x) = 0 \mbox{ for } 0 \le x \le \frac{T - \tau_{k+1}}{\tau_k} < 1 \mbox{ since } T < \tau_k + \tau_{k+1}.
\end{equation}
From the system of $v$ \eqref{eq-vj}, \eqref{eq-v-k+2}, \eqref{bdry-vv1}, and \eqref{bdry-vv0}, the solution $v$ is then uniquely determined by $v_{k}(T, x)$ for $\frac{T - \tau_{k+1}}{\tau_k} < x \le 1$.

Since, for $t \in (0, T)$,
$$
v_1(t, 1) = \cdots = v_{k-1} (t, 1) = 0 \quad  (\mbox{by } \eqref{bdry-vv1})
$$
and, for $(t, x) \in (0, T) \times (0, 1)$,
$$
\partial_t v_j(t, x) = -  \lambda_j \partial_x v_j(t, x) \mbox{ for } 1 \le j \le k-1 \quad (\mbox{by } \eqref{eq-v-k+2}),
$$
it follows from \eqref{final-1} that, for $t \in (0, T)$,
$$
v_1(t, 0) = \cdots = v_{k-1} (t, 0) = 0.
$$
We then derive from \eqref{assumption-B} and \eqref{bdry-vv-M} that, for $t \in (0, T)$,
\begin{equation}\label{eq-vvv-0}
v_{k+1}(t, 0) = \gamma_{k+1} v_{k}(t, 0) \quad  \mbox{ and } \quad  v_{k+\ell}(t, 0) = \gamma_{k+\ell} v_{k}(t, 0),
\end{equation}
where
\begin{equation}\label{def-gamma}
\gamma_{k+1} :=  \lambda_{k+1}^{-1}  \lambda_k B_{k, 1}  \mathop{\neq}^\eqref{assumption-B} 0 \quad \mbox{ and } \quad \gamma_{k+\ell} :=  \lambda_{k+\ell}^{-1}  \lambda_k B_{k, \ell} \mathop{\neq}^\eqref{assumption-B} 0.
\end{equation}
Since
$$
\partial_{t} v_{k}(t,x) \mathop{=}^{\eqref{eq-vj}} - \lambda_k \partial_x v_k(t, x),
$$
$v_{k}(t, 1) = 0$ for $t \in (0,  T)$ by \eqref{bdry-vv1}, and $T \ge \tau_{k} + \tau_{k+l}$,
one has, for $t \in (\tau_{k+\ell},  \tau_{k+1})$,
\begin{equation*}
\int_0^1 v_{k} (-\tau_{k+\ell} s  + t , s) \,  ds = \int_0^{\gamma_{k}(t)} v_{k} (-\tau_{k+\ell} s  + t , s) \,  ds
\end{equation*}
(see \Cref{fig1} for the definition of $\gamma_k(t)$).  This implies, by \eqref{eq-vj} applied with $i = k$, for $t \in (\tau_{k+\ell},  \tau_{k+1})$,
\begin{equation}\label{replacement-p1}
\int_0^1 v_{k} (-\tau_{k+\ell} s  + t , s)  ds = \theta_k  \int_{\tau_{k+\ell}}^{t } v_{k} (s, 0) \, ds,
\end{equation}
where
\begin{equation}\label{theta-k}
\theta_k =   \frac{\gamma_k(t)}{t - \tau_{k+\ell}}  \quad \left( =  \frac{1}{\sqrt{1 + \tau_{k+\ell}^2}} \frac{\sqrt{1 + \tau_{k+\ell}^2}}{\tau_{k} + \tau_{k+\ell}} = \frac{1}{\tau_k + \tau_{k+\ell}}\right)\quad  \mbox{: independent of $t$}.
\end{equation}

Similarly, since $$
\partial_{t} v_{k+1}(t,x) \mathop{=}^\eqref{eq-vj}  \lambda_{k+1} \partial_x v_{k+1}(t, x),
$$
$$
v_{k+1}(t, 1) = 0 \mbox{ for } t \in (0, T)
$$
 thanks to  $ \dsp v_{k}(t, 0) \mathop{=}^{\eqref{eq-vj}, \eqref{final-2}} 0$ for  $t \in (\tau_{k+1}, T)$ and \eqref{eq-vvv-0},  and $\dsp v_{k+1}(T, \cdot) \mathop{=}^\eqref{final-1} 0$,  we obtain, for $ t \in (\tau_{k+\ell}, \tau_{k+1})$,
\begin{equation*}
\int_0^1 v_{k+1} (-\tau_{k+\ell} s + t, s) \, ds =\int_0^{\gamma_{k+1}(t)} v_{k+1} (-\tau_{k+\ell} s  + t, s) \, ds
\end{equation*}
(see \Cref{fig1} for the definition of $\gamma_{k+1}(t)$).
\begin{figure}
\centering
\begin{tikzpicture}[scale=2.8]

\newcommand\x{0.5}

\draw[] (1+\x,0) -- (1 +\x ,2.0);
\draw[] (0,2) -- (1+ \x,2);

\draw (0, 2) node[left]{$T$};

\draw[] (1 +\x,0) -- (0, 0.2);
\draw (0, 0.2) node[left]{$\tau_{k+\ell}$};

\draw (0, 1) node[left]{$\tau_{k+1}$};

\draw (0, 0.5) node[left]{$t$};
\draw (1+\x, 0.3) node[right]{$t - \tau_{k+\ell}$};

   %%%%%%%%%%%%%%%%%
\draw [<->,thick] (0,2.2) node (yaxis) [above] {$t$}
        |- (1.7,0) node (xaxis) [right] {$x$};
    % Draw two intersecting lines
    \draw[] (0,0.2) coordinate (a_1) -- (1+ \x,1.5) coordinate (a_2);
  \draw[] (1+ \x,1.5) node[right]{$\tau_{k+\ell} + \tau_k$};
    \draw[dashed] (0,0.5) coordinate (b_1) -- (1+\x,0.3) coordinate (b_2);
    % Calculate the intersection of the lines a_1 -- a_2 and b_1 -- b_2
    % and store the coordinate in c.
    \coordinate (c) at (intersection of a_1--a_2 and b_1--b_2);
    % Draw lines indicating intersection with y and x axis. Here we use
    % the perpendicular coordinate system

    \draw[dashed] (c) -- (xaxis -| c) node[below] {$\gamma_k(t)$};

%%%%%%%%%%%%%

 \draw[] (0,1) coordinate (c_1) -- (1+ \x,0) coordinate (c_2);
   % \draw[dashed] (0,0.5) coordinate (b_1) -- (1+\x,0.3) coordinate (b_2);
    % Calculate the intersection of the lines a_1 -- a_2 and b_1 -- b_2
    % and store the coordinate in c.
    \coordinate (d) at (intersection of c_1--c_2 and b_1--b_2);
    % Draw lines indicating intersection with y and x axis. Here we use
    % the perpendicular coordinate system

    \draw[dashed] (d) -- (xaxis -| d) node[below] {$\gamma_{k+1}(t)$};

\draw (0, 0) node[left, below]{$0$};
\draw (1 + \x, 0) node[below]{$1$};

\end{tikzpicture}
\caption{On the definition of $\gamma_{k}$ and $\gamma_{k+1}$ for $t \in (\tau_{k+\ell}, \tau_{k+1})$:  $\gamma_k(t)$ is the abscise of the intersection of the line passing $(0, t)$ and $(1, t - \tau_{k+\ell})$,  and the line passing $(0, \tau_{k+\ell})$ and $(1, \tau_{k+\ell} + \tau_k)$;  $\gamma_{k+1}(t)$ is the abscise of the intersection of the line passing $(0, t)$ and $(1, t - \tau_{k+\ell})$,  and the line passing $(0, \tau_{k+1})$ and $(1, 0)$.}\label{fig1}
\end{figure}
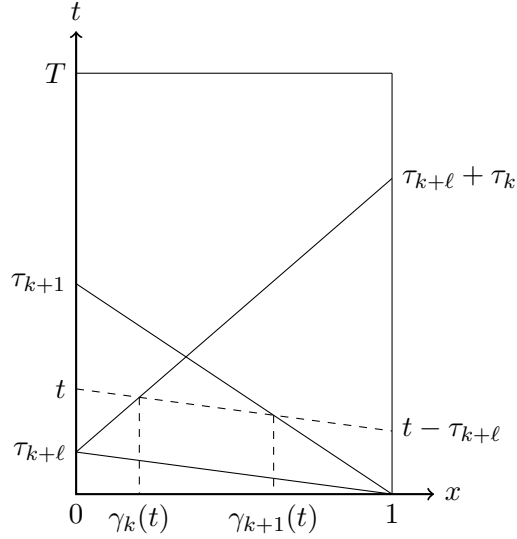

 This implies, by \eqref{eq-vj} applied with $i = k+1$, for $t \in (\tau_{k+\ell},  \tau_{k+1})$,
\begin{equation}\label{replacement-p2}
\int_0^1 v_{k+1} (-\tau_{k+\ell} s  + t, s) \, ds = \theta_{k+1} \int_{t}^{ \tau_{k+1}} v_{k+1}(s, 0) \, ds,
\end{equation}
where
\begin{equation}\label{theta-k+1}
\theta_{k+1} =  \frac{\gamma_{k+1}(t)}{ \tau_{k+1} - t}  \quad \left( = \frac{1}{\sqrt{1 + \tau_{k+\ell}^2}} \frac{\sqrt{1 + \tau_{k+\ell}^2}}{\tau_{k+1} - \tau_{k+\ell}} = \frac{1}{\tau_{k+1} - \tau_{k+\ell}} \right)  \quad  \mbox{: independent of $t$}.
\end{equation}

Using \eqref{replacement-p1} and \eqref{replacement-p2}, we derive from \eqref{eq-v4} that
\begin{equation*}
v_{k+\ell}(t, 0) =  \ta(t) \theta_k   \int_{\tau_{k+\ell}}^{t } v_{k} (s, 0) \,  ds +  \tb(t) \theta_{k+1}  \int_{t}^{ \tau_{k+1}} v_{k+1} (s, 0) \, ds \quad   \mbox{ for }  t \in (\tau_{k+\ell},  \tau_{k+1}).
\end{equation*}
This implies, by \eqref{eq-vvv-0},
\begin{equation}\label{eq-v4-p2}
v_{k+\ell}(t, 0) = \ha(t)  \int_{\tau_{k+\ell}}^{t } v_{k+\ell} (s, 0)  \, ds +  \hb(t)  \int_{t}^{ \tau_{k+1}} v_{k+\ell} (s, 0) \, ds \quad  \mbox{ for }  t \in (\tau_{k+\ell},  \tau_{k+1}),
\end{equation}
where
\begin{equation}\label{ta-ha}
\ha= \gamma_{k+\ell}^{-1} \theta_k \ta \quad \mbox{ and } \quad \hb = \gamma_{k+1} \gamma_{k+\ell}^{-1} \theta_{k+1} \tb.
\end{equation}

Since
$$
\tau_k + \tau_{k+1}  - \eps \mathop{=}^{\eqref{def-T}}  T,
$$
it follows that, at least if $\eps>0$ is small enough so that $T>\tau_{k+\ell}$,
$$
I: = (\tau_{k+\ell}, \tau_{k+1}) \cap (T - \tau_k, T) \not =  \emptyset.
$$

Fix $\varphi \in C^\infty_{\mc}(\mR)$ such that
\begin{equation}\label{cond-varphi}
\supp \varphi \subset I  \quad \mbox{ and } \quad \int_I \varphi = 1.
\end{equation}
Set
\begin{equation}\label{choice-1}
v_{k+\ell} (t, 0) =  \varphi (t)  \mbox{ for } t \in (\tau_{k+\ell}, \tau_{k+1}) \quad \mbox{ and } \quad  \ha (t) = \hb (t)= \varphi (t) \mbox{ for } t \in \mR.
\end{equation}
One can check that \eqref{eq-v4-p2} holds for this choice. From \eqref{eq-vvv-0}, we have
\begin{equation}\label{choice-2}
v_{k}(t, 0) = \gamma_{k+\ell}^{-1} \varphi (t) \quad  \mbox{ and } \quad  v_{k+1}(t, 0) =  \gamma_{k+1}  \gamma_{k+\ell}^{-1}\varphi(t) \quad \mbox{ for } t \in (\tau_{k+\ell}, \tau_{k+1}).
\end{equation}

We have just presented arguments for a choice of $\alpha$ and $\beta$,  and a choice of $v(T, \cdot)$ so that  \eqref{bdry-vv1}, \eqref{bdry-vv0}, and \eqref{main-vv} hold. We now proceed in the opposite direction to  rigorously establish this.

Consider $\alpha$ and $\beta$ defined by, for $(t, x) \in (0, T) \times (0, 1)$,
\begin{equation}\label{aaa}
\alpha(t, x) =  \lambda_{k+\ell} \gamma_{k+\ell}  \theta_k^{-1} \varphi (t + \tau_{k+\ell} x)
\end{equation}
and
\begin{equation}\label{bbb}
 \beta(t, x) =  \lambda_{k+\ell} \gamma_{k+1}^{-1}  \gamma_{k+\ell}  \theta_{k+1}^{-1} \varphi (t + \tau_{k+\ell} x),
\end{equation}
as suggested by \eqref{alpha-ta},  \eqref{ta-ha}, and \eqref{choice-1}, where $\varphi$ is determined as above.

Let $v(T, \cdot) \in C^\infty_{\mc}(0, 1)$ be such that \eqref{final-1} holds and
$v_{k}(T, \cdot)$ is chosen such that
\begin{equation*}
v_k(t, 0) = \gamma_{k+\ell}^{-1} \varphi (t) \mbox{ for } t \in (T - \tau_k, T),
\end{equation*}
as suggested by \eqref{final-2} and \eqref{choice-2}.  This implies, by \eqref{eq-vj} applied with $i=k$ and the fact $v_{k}(t, 1) = 0$ for $t \in (0, T)$ (see \eqref{bdry-v1}),
\begin{equation}\label{vk-0}
v_k(t, 0) = \gamma_{k+\ell}^{-1} \varphi (t) \mbox{ for } t \in (0, T)
\end{equation}
since $\supp \varphi \subset I \subset (T - \tau_k, T)$.  One can check that \eqref{eq-vvv-0} holds
by the same arguments used to derive it as before.  One  can also check that \eqref{eq-v4-p2} holds by \eqref{cond-varphi}.
Using \eqref{eq-vvv-0}, one then obtains
\eqref{eq-v4} for $t \in (\tau_{k+\ell}, \tau_{k+1})$, which implies \eqref{cond-v4} for $t \in (0, \tau_{k+1} - \tau_{k + \ell})$. From \eqref{cond-v4} being valid for $t \in (0, \tau_{k+1} - \tau_{k + \ell})$, and \eqref{eq-v-k+2} (see also \eqref{cond-v4-*}), we derive that
\begin{equation}\label{v-k+2-p1}
v_{k+\ell}(t, 1) = 0 \mbox{ for } t \in (0, \tau_{k+1} - \tau_{k+\ell}).
\end{equation}

Since $\ha = \hb = 0$ for $t \in (\tau_{k+1}, T)$, which implies  $\ta = \tb = 0$ for $t \in (\tau_{k+1}, T)$, it follows from \eqref{cond-v4-*} (see also \eqref{eq-v4}) that
\begin{equation}
v_{k+\ell}(t + \tau_{k+\ell}, 0) = v_{k+\ell} (t, 1) \mbox{ for } t \in (\tau_{k+1} - \tau_{k+\ell}, T - \tau_{k+ \ell}).
\end{equation}
This implies, by \eqref{eq-vvv-0} and \eqref{vk-0},
\begin{equation}\label{v-k+2-p2}
v_{k+\ell} (t, 1) = 0 \mbox{ for } t \in (\tau_{k+1} - \tau_{k+\ell}, T - \tau_{k+\ell}).
\end{equation}
Similarly, since $v_{k + \ell}(T, \cdot) = 0$, we derive from \eqref{aaa} and \eqref{bbb} \footnote{Since $\varphi(t) = 0$ for $t > \tau_{k+1}$, \eqref{aaa} and \eqref{bbb} imply that $\alpha = \beta = 0$ in the region of $(t,x)$ which is  below the characteristic flow of $v_{k+l}$ passing $(0, T)$ in the $xt$-plane.} that
\begin{equation}\label{v-k+2-p3}
v_{k+\ell} (t, 1) = 0 \mbox{ for } t \in (T - \tau_{k+\ell}, T).
\end{equation}
Combining \eqref{v-k+2-p1},  \eqref{v-k+2-p2}, and  \eqref{v-k+2-p3} yields
\begin{equation}\label{v-k+2-1}
v_{k+\ell} (t, 1) = 0 \mbox{ for } t \in (0, T).
\end{equation}

From the choice of $v(T, \cdot)$ in \eqref{final-1}, the property of $v$ given in \eqref{eq-vj}, and the fact $v_{-}(t, 1) = 0$ for $t \in (0, T)$, we have, for $1 \le j \le k-1$,
$$
v_{j}(t, 0) = 0 \mbox{ for } t \in [0, T].
$$
Since $B_{k, j} = 0$ for $2 \le j \le m$ with $j \neq \ell$ by \eqref{assumption-B}, it follows from \eqref{bdry-vv-M} that, for $k+2 \le j \le k+ m$ with $j \neq k + \ell$,
\begin{equation}\label{coucou-vj}
v_{j}(t, 0) = 0 \mbox{ for } t \in [0, T].
\end{equation}
We  derive from \eqref{eq-vj}, the choice of $v(T, \cdot)$  in \eqref{final-1}, and \eqref{coucou-vj} that, for  $k+2 \le j \le k+ m$ with $j \neq k + \ell$,
\begin{equation}\label{vj-1}
v_{j}(t, 1) = 0 \mbox{ for }  t \in [0, T].
\end{equation}
From \eqref{vk-0} (see also \eqref{final-2}), we obtain that $v_{k}(t, 0) = 0$ for $t \in (\tau_{k+1}, T)$. This implies, by \eqref{bdry-vv-M} and \eqref{coucou-vj} (see also \eqref{eq-vvv-0}),
\begin{equation*}
v_{k+1}(t, 0) = 0 \mbox{ for }  t \in (\tau_{k+1}, T).
\end{equation*}
We derive that, by using  \eqref{eq-vj} and   \eqref{final-1},
\begin{equation}\label{v-k+1-1}
v_{k+1}(t, 1) = 0 \mbox{ for }  t \in [0, T].
\end{equation}

We have, by  \eqref{eq-vvv-0} and \eqref{vk-0},
$$
v_{k+1}(t, 0) =  \gamma_{k+1}   \gamma_{k+\ell}^{-1}   \varphi(t) \mbox{ for } t \in (0, T).
$$
This implies, by  \eqref{eq-vj},
$$
v_{k+1}(0, x) =\gamma_{k+1}   \gamma_{k+\ell}^{-1} \varphi(\tau_{k+1} x) \mbox{ for } x \in [0, 1].
$$
We thus arrive, since $\supp \varphi \subset (0, \tau_{k+1})$,
\begin{equation}\label{cl2}
v(0, \cdot) \not \equiv 0.
\end{equation}

From \eqref{v-k+2-1}, \eqref{vj-1}, \eqref{v-k+1-1} and \eqref{cl2}, we reach
\begin{equation*}
v_+(\cdot, 1) =0 \quad \mbox{ and } \quad v(0, \cdot) \not \equiv 0.
\end{equation*}
The proof is complete.
\end{proof}

\section{Null-controllability in the  analytic setting - Proof of  \Cref{thm-A}} \label{sect-thmA}

This section is devoted to the proof of \Cref{thm-A}. The proof  is divided into three steps described below:

\medskip
$\bullet$ {Step 1:} for each $\tau$, we characterize the space $H(\tau)$ ($\subset [L^2(0, 1)]^n$), which is of finite dimension,  for which one can steer any element in $H(\tau)^\perp$ \footnote{Here and in what follows, for a closed subspace $E$ of $[L^2(0, 1)]^n$, we denote $\mbox{Prof}_E$ the projection to $E$, and $E^\perp$ its orthogonal complement, both with respect to the standard $L^2(0, 1)$-scalar product.} at time $\tau$ to 0 in time $T_{opt}$. (In particular, from this definition of $H(\tau)$, one cannot steer any element in $H(\tau) \setminus \{0 \}$ at time $\tau$ to 0 in time  $T_{opt}$.) Moreover,  we show that $H(\cdot)$ is analytic in a neighborhood $I_1$ of $[0, T_1 - T_{opt}]$ except for a discrete subset, which is removable  \footnote{The analyticity of $H(\tau)$ is understood via the analyticity of the mapping $\mbox{Prof}_{H(\tau)}$. This convention is used throughout the paper.}.

$\bullet$ {Step 2:}  For each $\tau \in I_1$,  we characterize the subspace $J(\tau)$ of $H(\tau)$ for which one can steer every  element $\varphi$ in $J(\tau)$ from time $\tau$ to 0 in  time $T_{opt, +}$, i.e., in time $T_{opt} + \delta$ for all $\delta > 0$.
Let $M(\tau)$ be the orthogonal complement of $J(\tau)$ in $H(\tau)$.  We also show that there exists a constant $\eps_0$ such that,  roughly speaking, the following property holds:  if $\tau \in I_1$ and $\varphi \in M(\tau) \setminus \{0\}$, then one {\it cannot} steer $\varphi$ from time $\tau$ to 0 in time $T_{opt} + \eps_0$.

$\bullet$ {Step 3:} We give the proof of \Cref{thm-A} using Steps 1 and 2.

\medskip
Let us make some comments on these three steps  before proceeding them.  Concerning Step 1,  the fact that $H(\tau)$ is of finite dimension already appeared in our previous analysis \cite{CoronNg19-2}. Some necessary conditions on $H(\tau)$ are derived in \cite{CoronNg19-2} and the starting point of the analysis there is the backstepping technique. In this paper, the (complete) characterization of $H(\tau)$ is given and it plays a crucial role in our proof of  \Cref{thm-A}.  This characterization  can be obtained by first applying the backstepping technique (and then by using similar ideas given here). However, this way requires a quite strong assumption on the analyticity of $C$  in the step of using backstepping technique (see \Cref{rem-backstepping}). To avoid it, we implement a new approach applied directly to the original system. The analysis is though strongly inspired/guided by our understanding in the form obtained via the backstepping.  A part of  technical points in this step is to establish the well-posedness of hyperbolic equations with unusual boundary conditions (the boundary condition of a component can be given both on the left at $x=0$ for some interval of time and on the right at  $x=1$ for some other interval of time),  and in a domain which is not necessary to be a rectangle in $xt$ plane. The analysis is interesting but delicate,  and presented in  the appendix.  After characterizing $H(\cdot) $, the analyticity of $H(\cdot)$ is established by suitably applying  the theory of perturbations of analytic compact operators, see, e.g., \cite{Kato}. These results are given in \Cref{pro-A1} in \Cref{sect-A1}. Concerning Steps 1 and  2, the characterizations of  all states for which one can steer from time $\tau$ to 0 in time $T_{opt}$ or in time $T_{opt, +}$ can be done for $C \in \big[L^\infty ( I \times (0, 1)) \big]^{n \times n}$. The analyticity of $C$ is not required for this purpose.  It is in the proof of the existence of $\eps_0$, given in Step 2, that the analyticity of $C$  plays a crucial role. The analysis of Step 3 is also based on  a technical lemma (\Cref{lem-AAA}). The approach proposed in this paper is quite robust and might be applied to other contexts.

\medskip

The rest of this section containing four subsections is organized as follows. In the first section, we introduce notations and present preliminary results related to observability inequalities, which are the starting point of our analysis. Steps 1, 2, and 3 are then given in  the second, the third, and the fourth subsection, respectively.

\subsection{Preliminaries}\label{sect-A1}

Fix $\tau \in I$ and $T>0$ such that $[\tau, \tau + T] \subset I$. Define
$$
\begin{array}{cccc}
\cF_{\tau, T}: &  [L^2(\tau, \tau +  T)]^m &  \to & [L^2(0, 1)]^n \\[6pt]
& U &  \mapsto & u(\tau +T, \cdot),
\end{array}
$$
where $u$ is the unique solution of the system
 \begin{equation}\label{eq-wtA}
\partial_t u (t, x) =  \Sigma(x) \partial_x u (t, x) + C(t, x) u(t, x)   \mbox{ for } (t, x)  \in  (\tau, \tau + T)  \times (0, 1),
\end{equation}
\begin{equation}\label{bdry-wtA}
u_- (t, 0)  = B u_+(t, 0) \mbox{ for } t \in (\tau, \tau + T),
\end{equation}
\begin{equation}\label{control-wtA}
u_+ (t, 1)= U(t) \mbox{ for } t \in (\tau, \tau + T),
\end{equation}
\begin{equation}\label{intial-wtA}
u (t = \tau, \cdot) = 0 \mbox{ in } (0, 1).
\end{equation}

Set, for $(t, x) \in I \times (0, 1)$,
\begin{equation}\label{def-hC}
\bC(t, x) = \Sigma'(x) - C\tr(t, x).
\end{equation}
The following result provides the formula for the adjoint $\cF_{\tau, T}^*$ of $\cF_{\tau, T}$.

\begin{lemma}\label{lem-FT*} We have, for $\varphi \in [L^2(0, 1)]^n$,
\begin{equation*}
\cF_{\tau,  T}^*(\varphi) = \Sigma_+ (1) v_+(\cdot, 1) \mbox{ in } (\tau, \tau +  T),
\end{equation*}
where $v$ is the unique solution of the system
\begin{equation}\label{eq-vA}
\partial_t v (t, x) =  \Sigma(x) \partial_x v (t, x) +  \bC(t, x)   v(t, x)
\mbox{ for } (t, x)  \in (\tau, \tau +  T) \times (0, 1),
\end{equation}
with, for $0 < t < T$,
\begin{equation}\label{bdry-v1A}
v_-(t, 1)  = 0,
\end{equation}
\begin{equation}\label{bdry-v0A}
\Sigma_+ (0) v_+(t, 0) = - B\tr \Sigma_- (0) v_- (t, 0),
\end{equation}
and
\begin{equation}\label{initial-vA}
v(t = \tau +  T, \cdot) = \varphi \mbox{ in } (0, 1).
\end{equation}
\end{lemma}

The proof of \Cref{lem-FT*} is quite standard and similar to the one of \cite[Lemma 1]{CoronNg19-2}. The details are  omitted.

\medskip

Using the same method, we also obtain the following two results, see,  e.g.,  the proof of \cite[Lemma 2]{CoronNg19-2} for the analysis.

\begin{lemma}\label{lem-ST1}  Assume that $u$ is a solution of \eqref{eq-wtA}-\eqref{control-wtA}  such that  $u_+(\cdot, 1) = 0$ in $(\tau, \tau + T)$. Then, for $\varphi \in [L^2(0, 1)]^n$, we have \footnote{The notation  $\langle \cdot, \cdot \rangle$ stands for  the Euclidean scalar product in $\mR^\ell$ for $\ell \ge 1$.}
\begin{equation*}
\int_0^1 \langle u(\tau + T, x),  v(\tau + T, x) \rangle  \, dx  = \int_0^1 \langle u (\tau, x), v(\tau, x ) \rangle  \, dx,
\end{equation*}
where $v$ is a solution of \eqref{eq-vA}-\eqref{initial-vA}.
\end{lemma}

\begin{lemma}\label{lem-ST2} Assume that $u$ is a solution of \eqref{eq-wtA}-\eqref{control-wtA}. Then
\begin{equation*}
\int_0^1 \langle u(\tau + T, x),  v(\tau + T, x) \rangle  \, dx  = \int_0^1 \langle u (\tau, x), v(\tau, x ) \rangle  \, dx,
\end{equation*}
where $v$ is a solution of \eqref{eq-vA}-\eqref{bdry-v0A} satisfying $v_+(\cdot, 1) = 0$.
\end{lemma}

Applying the Hilbert uniqueness method, see e.g. \cite[Chapter 2]{Coron07} and \cite{Lions88-VolI}, we have

\begin{lemma}\label{lem-observability} Let $E$ be a closed subspace of $[L^2(0, 1)]^n$. System \eqref{eq-wtA}-\eqref{control-wtA} is null controllable at the time $\tau + T$ for initial datum at time $\tau$ in $E$ if and only if, for some positive constant $C_{\tau, T}$,
\begin{equation}\label{observability}
\int_{\tau}^{\tau + T} |v_+(t, 1)|^2 \, dt \ge C_{\tau, T} \int_0^1 |\mbox{Proj}_{E} v(\tau, x)|^2 \, dx \quad  \forall \varphi \in [L^2(0, 1)]^n,
\end{equation}
where $v$ is the solution of \eqref{eq-vA}-\eqref{initial-vA}.
\end{lemma}

\subsection{Characterization of states at time $\tau$ steered  to 0 in time  $T_{opt}$.} \label{sect-A2}

In what follows in this section, we assume that $I =  (\alpha, \beta)$ is  an open bounded interval containing $[0, T_1]$ and set $I_1 = (\alpha, \beta - T_{opt})$.

\medskip
We first characterize states which can be steered at time $\tau$ to 0 in  time $T_{opt}$.  The following proposition is the key result of this section and is the starting point of our analysis in the analytic setting.

\begin{proposition}\label{pro-A1}
 Let $k \ge m \ge 1$ and let $B \in \cB$ be such that \eqref{cond-B-1} holds for $i =m$.  Assume that $C \in \big[L^\infty \big( I \times (0, 1)\big) \big]^{n \times n}$. There exist a
compact operator $\cK(\tau): [L^2(0, 1)]^n \to  [L^2(0, 1)]^n$ and a continuous linear operator $\cL(\tau): [L^2(0, 1)]^n \to  [L^2(0, T_{opt} - \tau_{k-m+1})]^m$ defined for $\tau \in I_1$ such that they are uniformly bounded in $I_1$ and,  with
\begin{equation}\label{def-H}
H(\tau) : = \Big\{\varphi \in [L^2(0, 1)]^n; \varphi + \cK(\tau) \varphi = 0 \mbox{ and } \cL(\tau) \varphi = 0 \Big\},
\end{equation}
the following two facts, concerning  system  \eqref{eq-wtA}-\eqref{control-wtA}, hold
\begin{itemize}
\item[i)] one can steer $\varphi \in H(\tau)^\perp$ at time $\tau$ to 0 at  time $\tau + T_{opt}$.

\item[ii)] one cannot steer any element $\varphi$ in  $H(\tau) \setminus \{0 \}$ at time $\tau$ to 0 at time $\tau + T_{opt}$.
\end{itemize}
Assume in addition that $C \in \cH(I, [L^\infty(0, 1)]^{n \times n})$. Then $\cK$ and $\cL$ are analytic in $I_1$.
\end{proposition}

We also obtain an explicit characterization of the space $H(\tau)$ in \Cref{pro-A1} via the dual system. In fact, the characterization of $H(\tau)$ in \eqref{def-H} is proved using such a characterization.   For the later use in the proof of \Cref{thm-A}, we state it in a slightly more general form:

\begin{lemma}\label{lem-AAA}  Let $k \ge m \ge 1$ and let $B \in \cB$ be such that \eqref{cond-B-1} holds for $i =m$.  Let $\varphi \in [L^2(0, 1)]^n$, $\tau \in I$,  and $T \ge T_{opt}$ be such that $\tau + T \in I$. There exists a subspace $H(\tau, T)$ of $H(\tau)$ such that
the following two facts, concerning system  \eqref{eq-wtA}-\eqref{control-wtA}, hold
\begin{itemize}
\item[i)] one can steer $\varphi \in H(\tau, T)^\perp$ at time $\tau$ to 0 at  time $\tau + T$.

\item[ii)] one cannot steer any element $\varphi$ in  $H(\tau, T) \setminus \{0 \}$ at time $\tau$ to 0 at time $\tau + T$.
\end{itemize}
Moreover, $\varphi \in H(\tau, T)$ if and only if there exists a solution $v$ of the system
\begin{equation}\label{AAA-1}
\partial_t v (t, x) =  \Sigma(x) \partial_x v (t, x) +  \bC(t, x)   v(t, x)
\mbox{ for } (t, x)  \in (\tau, \tau +  T) \times (0, 1),
\end{equation}
with, for $\tau < t < \tau + T$,
\begin{equation}\label{AAA-2}
v_-(t, 1)  = 0,
\end{equation}
\begin{equation}\label{AAA-3}
\Sigma_+ (0) v_+(t, 0) = - B\tr \Sigma_- (0) v_- (t, 0),
\end{equation}
\begin{equation}\label{AAA-4}
v_+(t, 1) = 0,
\end{equation}
and
\begin{equation}\label{AAA-5}
v(\tau, \cdot)  = \varphi.
\end{equation}
\end{lemma}

\begin{remark}  \rm Assume that the assumptions in \Cref{thm-E} hold. Let $\tau \in I$ and $T> T_{opt}$. Assume that  $\tau + T_{1} \in I$. We later prove that $H(\tau, T) = \{0\}$ (see \Cref{pro-UCP}) which is the unique continuation principle corresponding to  \eqref{AAA-1}-\eqref{AAA-4}.
\end{remark}

%\medskip

As a consequence of \Cref{pro-A1}
and the theory of analytic compact operators, see,  e.g.,  \cite{Kato}, we can prove

\begin{lemma}\label{lem-A2} Assume that  $C \in \cH(I, [L^\infty(0, 1)]^{n \times n})$. Then $H(\tau)$ is analytic in $I_1$ except for a discrete set, which is removable \footnote{The analyticity of $H(\tau)$ means the analyticity of $\mbox{Proj}_{H(\cdot)}$. }.
\end{lemma}

The proofs of \Cref{pro-A1}, \Cref{lem-AAA,lem-A2}  are given in the next three subsections, respectively.

Before entering the details of the proof, we introduce some notations on the characteristic flows which are used several times later. Extend $\lambda_i$ in $\mR$  with $1 \le i \le k+m$ by $\lambda_i(0)$ for $x< 0$ and $\lambda_i(1)$ for $x >1$.  For $(s, \xi) \in [0, + \infty) \times [0, 1]$, define $x_i(t, s, \xi)$ for $t \in \mR$ by
\begin{equation}\label{def-xi-1}
\frac{d}{d t} x_i(t, s, \xi) = \lambda_i \big(x_i(t, s, \xi) \big) \mbox{ and }  x_i(s, s, \xi) = \xi \mbox{ if } 1 \le i \le k,
\end{equation}
and
\begin{equation}\label{def-xi-2}
\frac{d}{d t} x_i(t, s, \xi) = - \lambda_i \big(x_i(t, s, \xi) \big) \mbox{ and }  x_i(s, s, \xi) = \xi \mbox{ if } k+1 \le i \le k+ m.
\end{equation}

\subsubsection{Proof of \Cref{pro-A1}}

\begin{figure}
\centering
\begin{tikzpicture}[scale=2.5]

\newcommand\x{3}

\draw[->] (0+\x, -0.2) -- (0+\x, 2.2);
\draw (0+\x, 2.2) node[left]{$t$};

\draw[->] (0+\x, 0) -- (1.55+\x, 0);
\draw (1.55+\x, 0) node[above, right]{$x$};

\draw[] (0+\x, 2) -- (1.4+\x, 2);
\draw (0+\x, 2) node[above, left]{$T_{opt}$};

\draw[-] (0+\x, 0) -- (1.4+\x, 0);
\draw (-0.1+\x,0) node[left, below]{$0$};

\draw (1.4+\x, 0) -- (1.4+\x, 2);
\draw (1.4+\x, 0) node[below]{$1$};

\draw[-, blue] (0+\x, 1.5) -- (1.4+\x, 2);
\draw[blue] (0.4+\x, 1.7) node[above]{$w_{k-m + 1}$};
\draw[blue] (0+\x, 1.5) node[left]{$T_{opt}- \tau_{k-m+ 1}$};

\draw [thick,-,fill=gray!5] (0+\x, 0) -- (0+\x, 1.5) -- (1.4+\x, 2) -- (1.4+\x, 0)-- cycle;

\draw (0.7 + \x, 0.8) node{$\Omega$};

\draw (0.7+\x, 0) node[below]{$w_+ = g$};

%\draw (0.7 + \x, -0.3) node[below]{$1 \le l \le k-m + 1$};

\draw[-, violet]  (1.4 + \x, 1.0) node[right]{$w = f$};
%\draw[-, violet]  (1.4 + \x, 0.8) node[right]{$w_{+} = f_{+}$};

%%%%%%%%%

\draw[->] (0, -0.2) -- (0, 2.2);
\draw (0, 2.2) node[left]{$t$};

\draw[->] (0, 0) -- (1.55, 0);
\draw (1.55, 0) node[above, right]{$x$};

\draw[] (0, 2) -- (1.4, 2);
\draw (0, 2) node[above, left]{$T_{opt}$};

\draw[-] (0, 0) -- (1.4, 0);
\draw (-0.1,0) node[left, below]{$0$};

\draw (0.7, -0.3) node[below]{$k-m+2 \le j \le k$};

\draw (1.4, 0) -- (1.4, 2);
\draw (1.4, 0) node[below]{$1$};

\draw [thick,-,fill=gray!5] (0, 0) -- (0, 1.5) -- (1.4, 2) -- (1.4, 0)-- cycle;
\draw (0.7, 0.8) node[]{$\Omega$};

%\draw[-,orange] (0, 0.7) -- (1.4, 0);
%\draw (0, 0.7) node[left]{$\textcolor{orange}{\tau_{j+m}}$};
%\draw[orange] (0.7, 0.4) node[above]{$w_{j+m}$};
\draw[-, orange] (0, 0.9) -- (1.4, 2);
\draw[-, orange] (1.1, 1.7)  node[below]{$w_j$};
%\draw[-, green]  (1.4, 1.8) node[right]{$\tau_{j+m} + \tau_j$};
%\draw[dashed, green]  (1.4, 2) -- (0, 0.9);

%node[right]{$\tau_{j+m} + \tau_j$};

%\draw[-, violet]  (1.4, 0.9) node[right]{$w_j = f_j$};
%\draw[-, violet, thick] (1.4, 0) --  (1.4, 1.8);

\draw[-, blue] (0, 1.5) -- (1.4, 2);
\draw[blue] (0.5, 1.74) node[above]{$w_{k-m + 1}$};
\draw[blue] (0, 1.5) node[left]{$T_{opt}- \tau_{k-m+ 1}$};
\draw[-, thick, red] (0, 0.9) -- (0, 1.5);
\draw[red] (-0.01, 1.2) node[left]{$w_{j}$};
\draw[red] (-0.01, 0.9) node[left]{$T_{opt} - \tau_{j}$};

\end{tikzpicture}
\caption{Geometry of the setting considered in the proof of \Cref{pro-A1} when $\Sigma$ is constant.   The boundary conditions  imposed  at $x=0$ for $w_j$ with $k-m+ 2 \le j \le k$ are given on the left, and  the boundary conditions  imposed at $x=1$ and $t=0$  are given on the right.
}\label{fig2}
\end{figure}
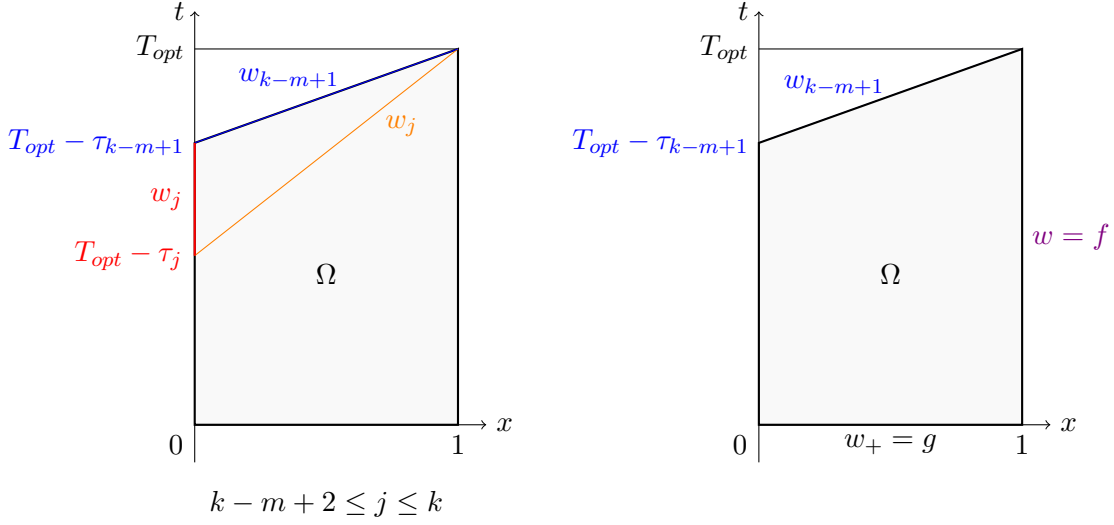

Fix $\tau \in I_1$.  Let $v$ be a solution of the system
\begin{equation}\label{eq-vAt}
\partial_t v (t, x) =  \Sigma(x) \partial_x v (t, x) +  \bC(t + \tau, x)   v(t, x)
\mbox{ for } (t, x)  \in (0, T_{opt}) \times (0, 1),
\end{equation}
with, for $0 < t < T_{opt}$,
\begin{equation}\label{bdry-vAt}
v_-(t, 1)  = 0,
\end{equation}
\begin{equation}\label{bdry-v0At}
\Sigma_+ (0) v_+(t, 0) = - B\tr \Sigma_- (0) v_- (t, 0),
\end{equation}
such that
\begin{equation}\label{lem-compact-cond-vA}
v_+(t, 1) = 0 \mbox{ for }  t \in (0, T_{opt}).
\end{equation}
Recall that $\bC$ is defined in \eqref{def-hC}.

The proof is now divided into two steps:

\medskip
\begin{itemize}
\item Step 1: We give a  characterization  of  $v(0, \cdot)$ where $v$ is a solution of \eqref{eq-vAt}-\eqref{bdry-v0At} satisfying \eqref{lem-compact-cond-vA}.

\item Step 2: We establish assertions i) and ii).
\end{itemize}

Step 1 is the key part of the proof. The operators $\cK(\tau)$ and $\cL(\tau)$ will be introduced in Step 1. The proof of Step 2 is quite standard after Step 1 and  the results in \Cref{sect-A1}.

\medskip

We now proceed with Steps 1 and 2.

\medskip
\noindent $\bullet$ {Step 1}:   For $1 \le i \le k \le j \le k+m$, we denote, for a vector $v \in \mR^{k+m}$,
$$
v_{-, \ge i} = (v_i, \cdots, v_k)
$$
and
$$
v_{< i, \ge j} = (v_1, \dots, v_{i-1}, v_{j}, \cdots, v_{k+m}).
$$

Using condition \eqref{cond-B-1} with $i=1$,
one can write the last equation of \eqref{bdry-v0At} in an equivalent form: 
\begin{equation}\label{pro-A1-vk}
v_{-, \ge k}(t, 0) = Q_k v_{<k, \ge k+m} (t, 0),
\end{equation}
for some $1 \times k$  matrix $Q_k$.

Using condition \eqref{cond-B-1} with $i=2$, one can write the last two equations of \eqref{bdry-v0At} in an equivalent form:
\begin{equation}\label{pro-A1-vk-1}
v_{-, \ge k-1}(t, 0) = Q_{k-1} v_{<k-1, \ge k+m-1} (t, 0),
\end{equation}
for some $2 \times k$ matrix $Q_{k-1}$.

\dots

Using condition \eqref{cond-B-1} with $i=m-1$, one can write the last $(m-1)$ equations of \eqref{bdry-v0At} in an equivalent form:
\begin{equation}\label{pro-A1-vk-m+2}
v_{-, \ge k-m+2}(t, 0) = Q_{k-m+2} v_{<k-m+2, \ge k+2} (t, 0),
\end{equation}
for some $(m-1) \times k$ matrix $Q_{k-m+2}$.

Using condition \eqref{cond-B-1} with $i=m$, one can write the last $m$ equations of \eqref{bdry-v0At} in an equivalent form:
\begin{equation}\label{pro-A1-vk-m+1}
v_{-, \ge k-m+1}(t, 0) = Q_{k-m+1} v_{<k-m+1, \ge k+1} (t, 0),
\end{equation}
for some $m \times k$ matrix $Q_{k-m+1}$.

Let $\Omega$ be the region of points $(t, x) \in (0, + \infty) \times (0, 1)$ such that in the $xt$-plane they are below the characteristic flow of $v_{k-m+1}$ passing the point $(1, T_{opt})$, see \Cref{fig2}.

Given $f \in [L^2(0, T_{opt})]^n$ and  $g \in [L^2(0, 1)]^m$,  we  consider the system
\begin{equation}\label{w1}
w_t (t, x)  = \Sigma (x) \partial_x w (t, x)+ \bC(t + \tau, x) w (t, x) \mbox{ for } (t, x) \in \Omega,
\end{equation}
\begin{equation}\label{w2}
w (\cdot, 1)= f  \mbox{ in } (0, T_{opt}),
\end{equation}
\begin{equation}\label{w3}
w_+(0, \cdot) = g \mbox{ in } (0, 1),
\end{equation}
\begin{equation}\label{w4}
w_{-, \ge k}(t, 0) = Q_k w_{<k, \ge k+m} (t, 0) \mbox{ for } t \in (T_{opt} - \tau_{k}, T_{opt} - \tau_{k-1}),
\end{equation}
\begin{equation}\label{w5}
w_{-, \ge k-1}(t, 0) = Q_{k-1} w_{<k-1, \ge k+m-1} (t, 0) \mbox{ for } t \in  (T_{opt} - \tau_{k-1}, T_{opt} - \tau_{k-2}),
\end{equation}
\dots
\begin{equation}\label{w6}
w_{-, \ge k-m+ 2}(t, 0) = Q_{k-m + 2} w_{<k-m+2, \ge k+2} (t, 0) \mbox{ for } t \in (T_{opt} - \tau_{k-m+2},  T_{opt} - \tau_{k - m + 1}),
\end{equation}
(see \Cref{fig2}).

For $\tau \in I_1$, define
\begin{equation}
\begin{array}{cccc}
\cT(\tau):  & [L^2(0, 1)]^n \times  [L^2(0, 1)]^m  & \to  & [L^2 (\Omega)]^{n} \\[6pt]
 & (f, g) & \mapsto & w,
\end{array}
\end{equation}
where $w$ is the (broad) solution of  \eqref{w1}-\eqref{w6} (see \Cref{def-WP}  for the definition of broad solutions  and \Cref{thm-WP} for their existence and uniqueness, both in the appendix).

\medskip
We claim that
\begin{equation}\label{claim-Step1}
\mbox{$v$ is a solution of \eqref{eq-vAt}-\eqref{bdry-v0At} satisfying \eqref{lem-compact-cond-vA} if and only if $v(0, \cdot) \in H(\tau)$ defined in \eqref{def-H},}
\end{equation}
where $\cK(\tau)$ and $\cL(\tau)$ are determined below.

\medskip
We now introduce $\cK$ and $\cL$.  Set  $w = \cT (\tau) \big(0, v_+(0, \cdot)\big)$. By noting that $v = w$,

\medskip

\noindent $\bullet$ $DK)$ the operator $\cK(\tau)$ is determined/summarized (the details are given below) by

\begin{enumerate}
\item[] $DK_m)$ the $m$ equations of  system \eqref{bdry-v0At}  imposed for $w$ in $(0, \tau_{k+m})$,

\item[] $DK_{m-1})$ the first $m-1$ equations  of system \eqref{bdry-v0At}  imposed for $w$  in $(\tau_{k+m}, \tau_{k+m-1})$,

\dots

\item[] $DK_1)$ the first equation of system \eqref{bdry-v0At}   imposed for $w$ in $(\tau_{k+2}, \tau_{k+1})$,

$\big($these above conditions are on $v_+(0, \cdot)$$\big)$,  and

\item[] $DK_-)$ $v_-(0, \cdot) = w_-(0, \cdot)$.

\end{enumerate}

\medskip
\noindent $\bullet$ $DL)$ $\cL(\tau)$ is defined by \eqref{bdry-v0At} in $(0, T_{opt} - \tau_{k-m+1})$.

\begin{figure}
\centering
\begin{tikzpicture}[scale=2.5]

\newcommand\z{0.0}

\draw[->] (0.0 + \z,0) -- (1.7+\z,0);
\draw[->] (0.0 +\z,0) -- (0.0 +\z,2.0);
\draw[] (1.5+\z,0) -- (1.5+\z,2.0);

\draw[green!50!black] (0+\z, 0) -- (1.5 + \z, 1);
\draw[green!50!black] (1.5 + \z, 1) node[right]{$\tau_\ell$};

\draw[green!50!black, dashed] (0.75 + \z, 0) -- (1.5 + \z, 0.5);
\draw (0.75+\z, 0) node[below]{$x$};
\draw[green!50!black] (\z + 1.5, 0.5) node[right]{$\tau(\ell, x)$};
\draw (1.5+\z, 0) node[below]{$1$};

\draw (\z, 2.0) node[left]{$t$};

\draw (\z, 0) node[left, below]{$0$};

%\draw[orange, dashed] (0.375+\zz, 0) -- (\zz, 0.4);
%\draw (\zz + 0.4, 0) node[below, orange]{$a_{i, j}(x)$};

\draw[] (\z + 0.85, -0.6) node[left]{$b)$};

%***********************
%Figure 2
%***********************

\newcommand\zz{-2.5}

\draw[->] (0+\zz,0) -- (1.7+\zz,0);
\draw[->] (0+\zz,0) -- (0+\zz,2.0);
\draw[] (1.5+\zz,0) -- (1.5+\zz,2.0);

\draw[blue] (1.5+\zz, 0) -- (\zz, 0.75);
\draw[blue] (\zz, 0.75) node[left]{$\tau_i$};

\draw[orange] (1.5+\zz, 0) -- (\zz, 1.5);
\draw[orange] (\zz + 0.5, 1.2) node[above]{\small $x_j(\cdot, 0, 1)$};
\draw[orange] (\zz, 1.5) node[left]{$\tau_j$};

\draw[orange, dashed] (0.75+\zz, 0) -- (\zz, 0.75);
\draw (\zz + 0.75, -0.2) node[below, orange]{$x_j (0, \tau_i, 0)$};

%\draw[blue] (1.5+\zz, 0) -- (\zz, 1.2);
%\draw[] (\zz + 1.1, 0.7) node[above]{\small $x_i(\cdot, 1, 0)$};

\draw (1.5+\zz, 0) node[below]{$1$};

\draw (\zz, 2.0) node[left]{$t$};

\draw (\zz, 0) node[left, below]{$0$};

\draw[blue, dashed] (0.3*15/7.5+\zz, 0) -- (\zz, 0.3);
\draw (0.3*15/7.5+\zz, 0) node[below]{$x$};

%\draw[orange, dashed] (0.375+\zz, 0) -- (\zz, 0.4);
%\draw (\zz + 0.4, 0) node[below, orange]{$a_{i, j}(x)$};

\draw[] (\zz, 0.3) node[left]{$\tau(i, x)$};
\draw[blue] (\zz + 0.55, 0.56) node[above]{\small $x_i(\cdot,  0, 1)$};

\draw[] (\zz + 0.85, -0.6) node[left]{$a)$};

\end{tikzpicture}
\caption{$\Sigma$ is constant; a) the definition of $x_i(\cdot, 0, 1)$, $x_j(\cdot, 0, 1)$, $x_j(0, \tau_i, 0)$, and $\tau_j(x)$ for $k+1 \le j < i \le k+m$; b) the definition of $x_\ell(\cdot, 0, 0)$ and $\tau(\ell, x)$ for $1 \le \ell \le k$.}\label{fig1*}
\end{figure}
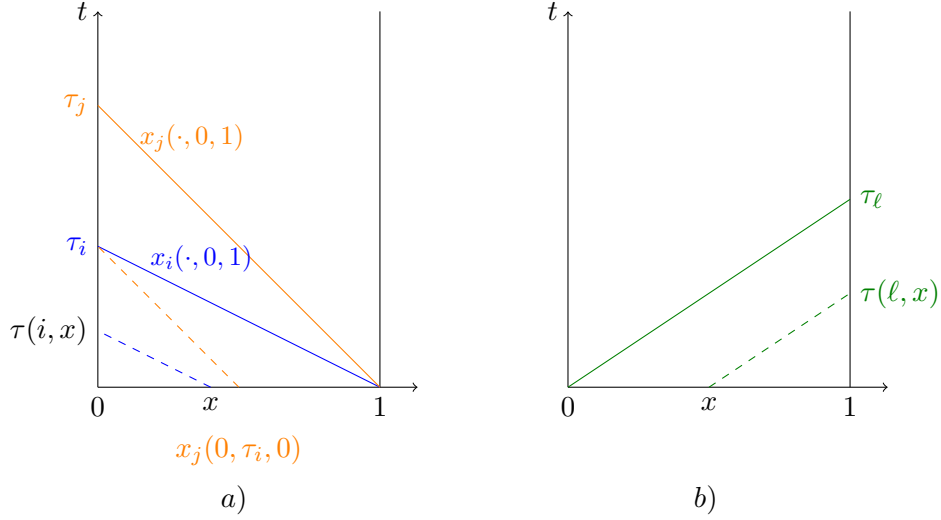

\medskip
Let us explain how to define the operators $\cK(\tau)$ and $\cL(\tau)$ from these conditions. To this end,  we first introduce some notations. For $x \in [0, 1]$ and $1 \le j \le k+m$,  let $\tau(j, x) \in [0, + \infty)$ be such that
$$
x_{j} \big(\tau(j, x), 0, x \big) = 0 \mbox{ for }  k+1 \le j \le k + m,
$$
and
$$
x_{j} \big(\tau(j, x), 0, x \big) = 1 \mbox{ for }  1 \le j \le k
$$
(see \Cref{fig1*}). Recall that $x_j(t, s, \xi)$ is defined in \eqref{def-xi-1} and \eqref{def-xi-2}.

We now consider $\cK(\tau)$ and first deal with the conditions $DK_m)$, \dots, $DK_1)$. The condition $DK_m)$ can be understood as follows. We have, for $1 \le j \le m$,
\begin{equation*}
\frac{d}{dt} w_{k+j} \big(t, x_{k+j} (t, 0, x) \big) = \Big( \bC \big(t + \tau, x_{k+j} (t, 0, x) \big) w \big(t, x_{k+j} (t, 0, x)\big) \Big)_{k+j}.
\end{equation*}
Integrating from $0$ to $\tau(k+j, x)$ yields, for $1 \le j \le m$ and for $x \in \big(0, x_{k+j} (0, \tau_{k+m}, 0\big)$,
$$
w_{k+j} (0, x) =   w_{k+j} \big(\tau(k+j, x), 0 \big) -  \int_0^{\tau(k+j, x)} \Big( \bC \big(t + \tau, x_{k+j} (t, 0, x) \big) w \big(t, x_{k+j} (t, 0, x) \big) \Big)_{k+j} \, dt.
$$
Using the $m$ equations of  system \eqref{bdry-v0At}, one has, for $1 \le j \le m$ and for $x \in \big(0, x_{k+j} (0, \tau_{k+m}, 0) \big)$,
\begin{multline}\label{DKm-A}
w_{k+j} (0, x) =  - \Big( \Sigma_+(0)^{-1} B\tr \Sigma_-(0) w_-(\tau(k+j, x), 0) \Big)_{j} \\[6pt] -  \int_0^{\tau(k+j, x)} \Big( \bC \big(t + \tau, x_{k+j} (t, 0, x) \big) w \big(t, x_{k+j} (t, 0, x) \big) \Big)_{k+j}\, dt.
\end{multline}

Similarly, the condition $DK_{m-1})$ can be written as,
for $1 \le j \le m-1$ and for \\$x \in \big(x_{k+j} (0, \tau_{k+m}, 0), x_{k+j} (0, \tau_{k+m-1}, 0) \big)$,
\begin{multline}\label{DKm-1-A}
w_{k+j} (0, x) =  - \Big( \Sigma_+(0)^{-1} B\tr \Sigma_-(0) w_-(\tau(k+j, x), 0) \Big)_{j} \\[6pt] -  \int_0^{\tau(k+j, x)} \Big( \bC \big(t + \tau, x_{k+j} (t, 0, x) \big) w \big(t, x_{k+j} (t, 0, x) \big) \Big)_{k+j} \, dt,
\end{multline}
\dots,  and the condition $DK_{1})$ can be written as, for $x \in \big(x_{k+1} (0, \tau_{k+2}, 0), x_{k+1} (0, \tau_{k+1}, 0) \big) =  \big(x_{k+1} (0, \tau_{k+2}, 0), 1 \big)$,
\begin{multline}\label{DKm-2-A}
w_{k+1} (0, x) =  - \Big( \Sigma_+(0)^{-1} B\tr \Sigma_-(0) w_-(\tau(k+1, x), 0) \Big)_{1} \\[6pt] -  \int_0^{\tau(k+1, x)} \Big( \bC \big(t + \tau, x_{k+1} (t, 0, x) \big) w \big(t, x_{k+1} (t, 0, x) \big) \Big)_{k+1} \, dt.
\end{multline}

We now deal with the condition $DK_-)$. We have, for $1 \le j \le k$,
\begin{equation*}
\frac{d}{dt} w_{j} \big(t, x_{j} (t, 0, x) \big) = \Big( \bC \big(t + \tau, x_{j} (t, 0, x) \big) w_{j} \big(t, x_{j} (t, 0, x)\big) \Big)_{j}.
\end{equation*}
Integrating from $0$ to $\tau(j, x)$ yields, for $1 \le j \le k$ and for $x \in (0, 1)$,
\begin{equation*}
w_{j} (0, x) =   w_{j} \big(\tau(j, x), 1 \big) -  \int_0^{\tau(j, x)} \Big( \bC \big(t + \tau, x_{j} (t, 0, x) \big) w \big(t, x_{j} (t, 0, x) \big) \Big)_j \, dt.
\end{equation*}
Since $f= 0$, it follows that, for $1 \le j \le k$ and $x \in (0, 1)$,
\begin{equation}\label{DK-A}
w_{j} (0, x) =   -  \int_0^{\tau(j, x)}  \Big( \bC \big(t + \tau, x_{j} (t, 0, x) \big) w \big(t, x_{j} (t, 0, x) \big) \Big)_j \, dt.
\end{equation}

The operator $\cK(\tau)$ is then defined via \eqref{DKm-A}-\eqref{DK-A}, with $v_+(0, \cdot) = \varphi_+$ and $w =  \cT\big(\tau, 0, v_+(0, \cdot)\big)$ as follows:

- for $1 \le j \le m$ and $x \in (0, 1)$, 
\begin{multline}
\Big(\cK(\tau) (\varphi) (x) \Big)_{k+j}=   \Big( \Sigma_+(0)^{-1} B\tr \Sigma_-(0) w_-(\tau(k+j, x), 0) \Big)_{j} \\[6pt] +  \int_0^{\tau(k+j, x)} \Big( \bC \big(t + \tau, x_{k+j} (t, 0, x) \big) w \big(t, x_{k+j} (t, 0, x) \big) \Big)_{k+j} \, dt. 
\end{multline}
%for $x \in (0, 1)$.  
%$x \in \big(x_{k+j} (0, \tau_{k+\ell +1}, 0), x_{k+j} (0, \tau_{k+\ell}, 0)  \big) $  with $j \le \ell \le m $ (the convention $\tau_{k+m + 1} = 0$ is used).

- for $1 \le j \le k$ and $x \in (0, 1)$,
\begin{equation}
\Big(\cK(\tau) (\varphi) (x) \Big)_{j} =    \int_0^{\tau(j, x)} \Big( \bC \big(t + \tau, x_{j} (t, 0, x) \big) w \big(t, x_{j} (t, 0, x) \big) \Big)_j \, dt.
\end{equation}

\medskip
Using \Cref{pro-WP} in the appendix, one can derive that  $\cK(\tau)$ is uniformly bounded in $I_1$ and is analytic in $I_1$ if $C \in \cH(I, [L^\infty(0, 1)]^{n \times n})$.

\medskip
The definition and the properties of $\cL(\tau)$ follow from $DL)$, with $v_+(0, \cdot) = \varphi_+$ and $w =  \cT\big(\tau, 0, v_+(0, \cdot)\big)$ as follows:
$$
\cL(\tau)(\varphi) = \Sigma_+ (0) w_+(t, 0) +  B\tr \Sigma_- (0) w_- (t, 0) \mbox{ in } (0, T_{opt} - \tau_{k-m+1}).
$$

It is clear that $H(\tau) \subset \Big\{\varphi \in [L^2(0, 1)]^n; \varphi + \cK(\tau) \varphi = 0 \mbox{ and } \cL(\varphi) = 0 \Big\}$. It remains to prove that
\begin{equation}\label{pro-AA1-p1}
\Big\{\varphi \in [L^2(0, 1)]^n; \varphi + \cK(\tau) \varphi = 0 \mbox{ and } \cL (\tau) \varphi = 0\Big\}   \subset H(\tau).
\end{equation}

To this end, we introduce another operator $\hat \cT$ related to $\cT$. Consider the system, for $(f, g) \in [L^2(0, 1)]^n \times [L^2(0, 1)]^m$,
\begin{equation}\label{hw1}
\partial_t \hw (t, x)  = \Sigma (x) \partial_x \hw (t, x)+ \bC(t + \tau, x) \hw (t, x) \mbox{ for } (t, x) \in (0, T_{opt}) \times (0, 1),
\end{equation}
\begin{equation}\label{hw2}
\hw (\cdot, 1)= f  \mbox{ in } (0, T_{opt}),
\end{equation}
\begin{equation}\label{hw3}
\hw_+(0, x) = g(x) \mbox{ in }  (0, 1),
\end{equation}
\begin{equation}\label{hw4}
\hw_i(T_{opt}, \cdot) = 0 \mbox{ in }  (0, 1), \mbox{ for } 1 \le i \le k-m,
\end{equation}
\begin{equation}\label{hw5}
\hw_{-, \ge k}(t, 0) = Q_k \hw_{<k, \ge k+m} (t, 0)  \mbox{ for } t \in (T_{opt} - \tau_{k}, T_{opt} - \tau_{k-1}),
\end{equation}
\begin{equation}\label{hw6}
\hw_{-, \ge k-1}(t, 0) = Q_{k-1} \hw_{<k-1, \ge k+m-1} (t, 0)  \mbox{ for } t \in (T_{opt} - \tau_{k-1}, T_{opt} - \tau_{k-2}),
\end{equation}
\dots
\begin{equation}\label{hw7}
\hw_{-, \ge k-m+ 2}(t, 0) = Q_{k-m + 2} \hw_{<k-m+2, \ge k+ 2} (t, 0)  \mbox{ for } t \in (T_{opt} - \tau_{k-m+2}, T_{opt} - \tau_{k-m+1}),
\end{equation}
\begin{equation}\label{hw8}
\hw_{-, \ge k-m+ 1}(t, 0) = Q_{k-m+ 1} \hw_{<k-m+1, \ge k+1} (t, 0) \mbox{ for } t \in (T_{opt} - \tau_{k-m+1}, T_{opt})
\end{equation}
(it is at this stage that the condition \eqref{cond-B-1} with $i =m$ is required!).

For $\tau \in I_1$, define
\begin{equation}
\begin{array}{cccc}
\hat \cT(\tau)  : &   [L^2(0, 1)]^n \times  [L^2(0, 1)]^m  & \to  & [L^2 \big( (0, T_{opt}) \times (0, 1) \big)]^{n} \\[6pt]
 & (f, g) & \mapsto & \hw,
\end{array}
\end{equation}
where $\hw$ is the unique broad solution of  \eqref{hw1}-\eqref{hw8} (see \Cref{thm-hWP} in the appendix for the existence and uniqueness of broad solutions; the definition of broad solutions is similar to \Cref{def-WP}).

It is clear that
\begin{equation}\label{restriction}
\mbox{$\cT(\tau) (0, g)$ is the restriction of $\hat \cT(\tau) (0, g)$ in $\Omega$ for $g \in [L^2(0, 1)]^m$. }
\end{equation}

Fix
\begin{equation}\label{pro-V}
\varphi_0 \in  \Big\{\varphi \in [L^2(0, 1)]^n; \varphi + \cK(\tau) \varphi = 0 \mbox{ and } \cL(\tau) \varphi = 0 \Big\}.
\end{equation}
Denote
$$
w = \cT(\tau) (0, \varphi_{0, +}) \quad \mbox{ and } \quad \hw = \hat \cT(\tau) (0, \varphi_{0, +}).
$$
Then, by \eqref{restriction},
\begin{equation}\label{w-hw}
\hw = w \mbox{ in } \Omega.
\end{equation}
Since $\varphi + \cK(\tau) (\varphi) = 0$ (see also the condition $DK_-)$), we have
$$
w(0, \cdot) = \varphi_0 \mbox{ in } (0, 1).
$$
Since $\cL(\tau) (\varphi) = 0$, we obtain
\begin{equation}\label{hw9}
\Sigma_+ (0) \hw_+(t, 0) = - B\tr \Sigma_- (0) \hw_- (t, 0) \mbox{ for } t \in (0, T_{opt} - \tau_{k-m+1}).
\end{equation}
On the other hand, by the definition of $\hat \cT$ (in particular, condition \eqref{hw8}), one has,
\begin{equation}\label{hw10}
\Sigma_+ (0) \hw_+(t, 0) = - B\tr \Sigma_- (0) \hw_- (t, 0) \mbox{ for } t \in (T_{opt} - \tau_{k-m+1}, T_{opt}).
\end{equation}
Combining \eqref{hw9} and \eqref{hw10} yields
\begin{equation}\label{hw-p2}
\Sigma_+ (0) \hw_+(t, 0) = - B\tr \Sigma_- (0) \hw_- (t, 0) \mbox{ for } t \in (0, T_{opt}).
\end{equation}
Thus $\hat w$ is a solution of \eqref{eq-vAt}-\eqref{bdry-v0At} satisfying \eqref{lem-compact-cond-vA} with $\hat w(0, \cdot) = w (0, \cdot) = \varphi_0$.

\medskip
\noindent $\bullet$ {Step 2:} We derive i) and ii). We begin with assertion ii). Let $\varphi \in H(\tau) \setminus \{0\}$ be arbitrary. By Step 1, there exists a solution  $v$ of
\eqref{eq-vAt}-\eqref{bdry-v0At} such that
\begin{equation}
v_+(\cdot, 1) = 0 \mbox{ in } (0, T_{opt}) \quad \mbox{ and } \quad v(0, \cdot) = \varphi \mbox{ in } (0, 1).
\end{equation}
Set
$$
v^{(\tau)}(t, x) = v(t - \tau, x) \mbox{ for } (t, x) \in (\tau, \tau + T_{opt}).
$$
Let $w$ be a solution of \eqref{eq-wtA} - \eqref{control-wtA} \footnote{Condition\eqref{control-wtA}  means that $w_+(t, 1) \in [L^2(\tau, \tau+ T)]^m$.} with $T = T_{opt}$, in which $ u$ is replaced by $w$,  with $w(\tau, \cdot) = v^{(\tau)}(\tau, \cdot) = \varphi$.
By \Cref{lem-ST2}, we have
\begin{equation*}
\int_0^1 \langle w(T_{opt} + \tau, x),  v^{(\tau)}(T_{opt} + \tau, x) \rangle  \, dx  = \int_0^1 \langle w (\tau, x), v^{(\tau)}(\tau, x ) \rangle  \, dx = \int_0^1 |\varphi |^2 \neq 0.
\end{equation*}
Therefore, one cannot steer $\varphi$ from time $\tau$ to $0$ at time $\tau + T_{opt}$.

We next establish assertion i) by a contradiction argument. Assume that this is not true. By \Cref{lem-observability} with $E = H(\tau)^\perp$, there exists a sequence of solutions $(v_N)$ of
\eqref{eq-vAt}-\eqref{bdry-v0At} such that
\begin{equation}
\lim_{N \to + \infty}\| v_{N, +}(\cdot, 1)\|_{L^2(0, T_{opt})} = 0   \quad \mbox{ and } \quad \| \mbox{Proj}_{H(\tau)^\perp} v_N(0, \cdot)\|_{L^2(0, 1)} = 1.
\end{equation}

Set
$$
\varphi_N = \mbox{Prof}_{H(\tau)} v_N(0, \cdot) \in H(\tau) \subset  [L^2(0, 1)]^n.
$$
Define
$$
w_N = \cT\big(\tau) (0, \varphi_{N, +} \big) \mbox{ in } \Omega \quad \mbox{ and } \quad  \hw_N = \hat \cT\big(\tau) (0, \varphi_{N, +} \big) \mbox{ in } (0, T_{opt}) \times (0, 1).
$$
Since $\varphi_N \in H(\tau)$, it follows from the definition of $\cK(\tau)$ that
$$
w_N(0, \cdot) = \varphi_N \mbox{ in } (0, 1).
$$
We derive from \eqref{restriction} that
$$
\hw_N(0, \cdot) = \varphi_N \mbox{ in } (0, 1).
$$
Replacing $v_N$ by $v_N - \hw_N$ if necessary,  without loss of generality,  one can assume in addition that $v_N(0, \cdot) \in H(\tau)^\perp$,  which  yields in particular that $\| v_N(0, \cdot)\|_{L^2(0, 1)} = \| \mbox{Proj}_{H(\tau)^\perp} v_N(0, \cdot) \|_{L^2(0, 1)} = 1$. This will be assumed from now on.

Consider $f_N \in  [L^2(0, T_{opt})]^n$ defined by
$$
f_N = v_N (\cdot, 1).
$$
Since $v_{N, -} (\cdot, 1 ) =  0$ in $(0, T_{opt})$ and $
\lim_{N \to + \infty}\| v_{N, +}(\cdot, 1)\|_{L^2(0, T_{opt})} = 0
$, it follows that
\begin{equation}\label{fN}
\lim_{N \to + \infty} f_N = 0 \mbox{ in $ [L^2(0, T_{opt})]^n$. }
\end{equation}
Set, in $\Omega$,
$$
u_{N}  =  \cT(\tau) (f_N, v_{N, +}(0, \cdot)).
$$
Then
\begin{equation}\label{uv-N}
u_N = v_N \mbox{ in } \Omega.
\end{equation}
Since $v_N(0, \cdot) + \cK(\tau) v_N(0, \cdot) = 0$, $\| v_N(0, \cdot)\|_{L^2(0, 1)} =1$,
and $\cK(\tau)$ is compact, it follows that, for some subsequence,  $v_{N_k}(0, \cdot) \to \varphi$ in $[L^2(0, 1)]^n$ and hence $\varphi \in H(\tau)^\perp$ by \eqref{fN} and the continuity of $\cT(\tau)$ (see \Cref{pro-WP} in the appendix). Set $u = \cT(\tau) (0, \varphi_+)$. Since, by \eqref{uv-N},
$$
\Sigma_+(0) u_{N, +} (t, 0)  = - B\tr \Sigma_-(0) u_{N, +} (t, 0) \mbox{ for } t \in (0, T_{opt}),
$$
and
$$
v_{N} (0, \cdot) = \cT(\tau) (f_N, v_{N, +} (0, \cdot)) \mbox{ in } (0, 1),
$$
we derive from \eqref{fN}, \eqref{uv-N} and  the continuity of $\cT(\tau)$ (see \Cref{pro-WP} in the appendix) that
$$
\Sigma_+(0) u_{+} (t, 0)  = - B\tr \Sigma_-(0) u_{+} (t, 0) \mbox{ for } t \in (0, T_{opt}),
$$
and
$$
\varphi(x)  =   \cT(\tau) (0, \varphi_+)(0,x) \mbox{ for } x \in  (0, 1).
$$
This implies that  $\varphi \in H(\tau)$.  It follows that  $\varphi = 0$ since $\varphi  \in H(\tau)^\perp$. We deduce  that
$$
0 =\| \mbox{Proj}_{H(\tau)^\perp} \varphi \| = \lim_{k \to + \infty}\| \mbox{Proj}_{H(\tau)^\perp} v_{N_k}(0, \cdot)\| = 1.
$$
We have a contradiction. Assertion i) is proved.

\medskip
The proof is complete. \qed

\subsubsection{Proof of \Cref{lem-AAA}}  The proof of \Cref{lem-AAA} follows from the one of \Cref{pro-A1} by replacing $T_{opt}$ by $T$. One just notes here that $H(\tau, T)$ is the set of  $v(0, \cdot)$ where $v$ is a solution of \eqref{eq-vAt}-\eqref{bdry-v0At} satisfying \eqref{lem-compact-cond-vA} with $T_{opt}$ replaced by $T$. The details of the proof are omitted.

\subsubsection{Proof of \Cref{lem-A2}}  \label{proof-lem-A2}
Set, for $\tau \in I_1$,
\begin{equation}\label{def-E}
E(\tau) = \Big\{\varphi \in [L^2(0, 1)]^n; \varphi + \cK(\tau) \varphi = 0 \Big\},
\end{equation}
and
\begin{equation}\label{def-P}
\mbox{let $P(\tau)$ be the generalized eigenspace  of $\cK(\tau)$ with respect to the eigenvalue $-1$}.
\end{equation}
From \eqref{def-H}, we have, for $\tau \in I_1$,
\begin{equation}
H(\tau) = E(\tau) \cap  \Big\{\varphi \in [L^2(0, 1)]^n; \cL(\tau) \varphi  = 0 \Big\}.
\end{equation}

Applying the theory of the perturbation of analytic compact operators, see e.g. \cite{Kato}, one derives that
\begin{equation}\label{lem-A2-p1}
P(\tau) \mbox{ is analytic in $I_1$ except for a discrete subset, which is removable.}
\end{equation}
Indeed, since $K(\tau)$ is compact, it follows that the eigenvalue $-1$ of $K(\tau)$ is isolated for each $\tau \in I_1$. From \cite[Section 3 of Chapter 7]{Kato} (see also \cite[Section 3 of Chapter 2]{Kato}),  for each $\tau$ there is a $\gamma > 0$  ($\gamma$ depends on $\tau$) such that the sum of the eigenprojections for all the eigenvalues of $\cK(\tau)$ lying inside $\{z \in \mC: |z + 1|< \gamma \}$ is analytic. We now can apply the theory of the perturbation of analytic operators in a finite dimensional space,  see  \cite[Chapter 2]{Kato}, involving the theory of algebraic functions, see e.g. \cite[Section 2 of Chapter 8]{Ahlfors},    to derive \eqref{lem-A2-p1}.

We have
$$
E(\tau) \mathop{=}^{\eqref{def-E}} \Big\{\varphi \in [L^2(0, 1)]^n; \varphi + \cK(\tau) \varphi = 0 \Big\} \mathop{=}^{\eqref{def-P}} \Big\{\varphi \in P(\tau); \varphi + \cK(\tau) \varphi = 0 \Big\}.
$$
It follows that
\begin{equation}
H(\tau)  =  \Big\{\varphi \in P(\tau); \varphi + \cK(\tau) \varphi = 0 \mbox{ and } \cL(\tau) \varphi = 0  \Big\}.
\end{equation}
We now can use the theory of the perturbation of the null-space  of analytic matrices. Applying \cite[Theorem S6.1 on page 388-389]{GLR-Matrix} and using \eqref{lem-A2-p1},
we derive that  \footnote{One way to  apply the theory of the perturbation of the null-space  of analytic matrices can be done as follows. One can first locally choose an analytic orthogonal  basis $\Big\{ \varphi_1(\tau), \cdots, \varphi_{\ell}(\tau) \Big\} $ of $P(\tau)$. We then represent the operator $Id + \cK(\tau)$ (where $Id$ denotes the identity map) in this basis after noting that it is an application from $P(\tau)$ into $P(\tau)$. We also represent $\cL(\tau)$ using the set $\Big\{ \cL (\tau) (\varphi_1(\tau)), \dots, \cL(\tau) (\varphi_{\ell} (\tau)) \Big\}$ for which the dimension of the span space  is constant outside a discrete set, which is removable. }
\begin{equation}\label{lem-A2-p3}
H(\tau) \mbox{ is analytic in $I_1$ except for a discrete subset, which is removable.}
\end{equation}
The proof is complete.    \qed

\subsection{Characterization of states at time $\tau$ steered to 0 in  time  $T_{opt, +}$} \label{sect-B2}

Fix  $\gamma_0 > 0$ such that $[0, T_1] \subset (\alpha + \gamma_0, \beta -\gamma_0)$. Set
\begin{equation}\label{def-I2}
I_2 = (\alpha + \gamma_0, \beta -\gamma_0 - T_{opt}).
\end{equation}

Given $0 < \eps < \gamma_0$ and $\tau \in I_2$, consider the system, for  $V \in [L^2(0, \eps)]^m$,
\begin{equation}\label{Sys-UA}
\left\{\begin{array}{c}
\partial_t v (t, x) =  \Sigma(x) \partial_x v (t, x) + C(t + \tau, x) v(t, x) \mbox{ for } (t, x)  \in (0,  \eps)  \times (0, 1), \\[6pt]
v_-(t, 0) = B v_+(t, 0) \mbox{ for } t \in (0,  \eps), \\[6pt]
v_+(t, 1) = V(t) \mbox{ for } t \in (0,  \eps), \\[6pt]
v(0, \cdot) = 0 \mbox{ in } [0, 1].
\end{array}\right.
\end{equation}
Define \footnote{The sub-index $c$ means that controls are used.}
$$
\begin{array}{cccc}
\cT^c_{\tau, \eps}: & [L^2(0, \eps)]^m & \to &  [L^2(0, 1)]^n \\[6pt]
& V & \mapsto & v(\eps, \cdot),
\end{array}
$$
where $v$ is the solution of \eqref{Sys-UA}. Consider two subsets $Y_{\tau, \eps}$ and $A_{\tau, \eps}$ \footnote{The letter $A$ means the attainability.} of $[L^2(0,1)]^n$ defined by
\begin{equation}\label{def-YA}
Y_{\tau, \eps} = \cT^c_{\tau, \eps} \Big\{ [L^2(0, \eps)]^m \Big\} \quad \mbox{ and } \quad A_{\tau, \eps}  =  \mbox{Proj}_{H(\tau + \eps)}  \Big\{ Y_{\tau, \eps} \Big\}.
\end{equation}

Given $0< \eps < \gamma_0$ and $\tau \in I_2$, we also define \footnote{The sub-index $I$ means that initial data are considered.}
$$
\begin{array}{cccc}
\cT^I_{\tau, \eps}: & [L^2(0, 1)]^n & \to &  [L^2(0, 1)]^n \\[6pt]
& \varphi  & \mapsto & w(\eps, \cdot),
\end{array}
$$
where $w$ is the solution of
\begin{equation}\label{Sys-UA0}
\left\{\begin{array}{c}
\partial_t w (t, x) =  \Sigma(x) \partial_x w (t, x) + C(t + \tau, x) w(t, x) \mbox{ for } (t, x)  \in (0, \eps)  \times (0, 1), \\[6pt]
w_-(t, 0) = B w_+(t, 0) \mbox{ for } t \in (0, \eps), \\[6pt]
w_+(t, 1) = 0 \mbox{ for } t \in (0, \eps), \\[6pt]
w(0,  \cdot) = \varphi \mbox{ in } [0, 1].
\end{array}\right.
\end{equation}

Set, for $0< \eps < \gamma_0$ and for $\tau \in I_2$,
\begin{equation}\label{def-J}
J(\tau, \eps) : = \Big\{\varphi \in H(\tau);  \mbox{Proj}_{H(\tau + \eps)} \cT^I_{\tau, \eps} (\varphi)
 \in A_{\tau, \eps} \Big\}.
\end{equation}

The motivation for the definition of $\cT^c_{\tau, \eps}$ and $\cT^I_{\tau, \eps}$ is:

\begin{lemma} \label{lem-J} Let $0< \eps < \gamma_0$ and  $\tau \in I_2$. Then
$J(\tau, \eps)$ is the space of (functions) states in $H(\tau)$ such that one can steer them from time $\tau$ to 0 at time $\tau + T_{opt} + \eps$.  As a consequence, for $\tau \in I_2$,
\begin{equation}\label{lem-J-st1}
J(\tau, \eps') \subset J(\tau, \eps) \mbox{ for } 0 < \eps' < \eps < \gamma_0,
\end{equation}
and the limit $J(\tau)$ of $J(\tau, \eps)$ as $\eps \to 0_+$ exists.
\end{lemma}

\begin{remark}  \rm The monotone property of $J(\tau, \eps)$ with respect to $\eps$ given in \eqref{lem-J-st1} will play a role in our analysis.
\end{remark}

\begin{remark} \rm The analyticity of $C$ in $I$ is not required in \Cref{lem-J}.
\end{remark}

\begin{proof}[Proof of \Cref{lem-J}]
Given $\varphi \in J(\tau, \eps)$, by the definition of $J(\tau, \eps)$,  there exists $\hat V \in [L^2(0, \eps)]^m$ such that
$$
\mbox{Proj}_{H(\tau + \eps)} \hw(\eps, \cdot)  = 0,
$$
where $\hw$ defined in $(0, \eps) \times (0, 1)$ is the solution of the system
\begin{equation}
\left\{\begin{array}{c}
\partial_t \hw (t, x) =  \Sigma(x) \partial_x \hw (t, x) + C(t + \tau, x) \hw(t, x) \mbox{ for } (t, x)  \in (0, \eps)  \times (0, 1), \\[6pt]
\hw_-(t, 0) = B \hw_+(t, 0) \mbox{ for } t \in (0, \eps), \\[6pt]
\hw_+(t, 1) = \hat V \mbox{ for } t \in (0, \eps), \\[6pt]
\hw(0,  \cdot) = \varphi \mbox{ in } [0, 1].
\end{array}\right.
\end{equation}
It follows that, by the characterization of $H(\tau + \eps)$,  there exists  $\widetilde V \in [L^2(\eps,T_{opt} + \eps)]^m$ such that
$$
\tw(T_{opt} + \eps, \cdot) = 0 \mbox{ in } (0, 1),
$$
where $\tw$ defined in $(\eps, T_{opt} + \eps) \times (0, 1)$ is the solution  of the system
\begin{equation}
\left\{\begin{array}{c}
\partial_t \tw (t, x) =  \Sigma(x) \partial_x \tw (t, x) + C(t + \tau, x) \tw(t, x) \mbox{ for } (t, x)  \in (\eps, T_{opt} +  \eps )  \times (0, 1), \\[6pt]
\tw_-(t, 0) = B \tw_+(t, 0) \mbox{ for } t \in (\eps, T_{opt} +  \eps ) , \\[6pt]
\tw_+(t, 1) = \widetilde V \mbox{ for } t \in (\eps, T_{opt} +  \eps ) , \\[6pt]
\tw(\eps,  \cdot) = \hw (\eps, \cdot) \mbox{ in } [0, 1].
\end{array}\right.
\end{equation}

Let $w$ be defined in $(0, T_{opt} + \eps) \times (0, 1)$ by $\hw$ in $(0, \eps) \times (0, 1)$  and by  $\tw$ in $(\eps, T_{opt} + \eps) \times (0, 1)$. Set
$$
\bw (t, x)  = w (t - \tau, x) \mbox{ in } (\tau, \tau +  T_{opt} + \eps) \times (0, 1).
$$
Then $\bw$
 is a solution starting from $\varphi$ at time $\tau$ and arriving at 0 at time $\tau + T_{opt} + \eps$, i.e.,
 \begin{equation}
\left\{\begin{array}{c}
\partial_t \bw (t, x) =  \Sigma(x) \partial_x \bw (t, x) + C(t, x) \bw(t, x) \mbox{ for } (t, x)  \in (\tau, \tau+ T_{opt} +  \eps )  \times (0, 1), \\[6pt]
\bw_-(t, 0) = B \bw_+(t, 0) \mbox{ for } t \in (\tau, \tau+ T_{opt} +  \eps ) , \\[6pt]
\bw(\tau,  \cdot) = \varphi  \mbox{ and } \bw(\tau + T_{opt} + \eps, \cdot) = 0 \mbox{ in } [0, 1].
\end{array}\right.
\end{equation}
We have thus proved that one can steer  $\varphi \in J(\tau, \eps)$ at time $\tau$ to 0 at time $\tau + T_{opt} + \eps$.

Conversely, let $\varphi \in H(\tau)$ be such that one can steer $\varphi$ at  time $\tau$ to $0$ at time $\tau + T_{opt} + \eps$ using a  control  $W \in [L^2(\tau, \tau + T_{opt} + \eps)]^m$. Let $\bw$ be the corresponding solution, and set $w(t, x) = \bw(t + \tau,x)$ in  $(0, T_{opt} + \eps) \times (0, 1)$.  Since $\bw(\tau + \eps, \cdot)$ is steered from time $\tau + \eps$ to 0 at time $\tau + T_{opt} + \eps$, it follows from the characterization of $H(\tau + \eps)$ that
$$
\mbox{Proj}_{H(\tau + \eps)} \bw(\tau + \eps, \cdot)  = 0.
$$
In other words,
$$
\mbox{Proj}_{H(\tau + \eps)} w(\eps, \cdot)  = 0.
$$
This yields  that $\varphi \in J(\tau, \eps)$.

\medskip
We thus proved that $J(\tau, \eps)$ is the space of (functions) states in $H(\tau)$ such that one can steer them from time $\tau$ to 0 at time $\tau + T_{opt} + \eps$. The other conclusions of \Cref{lem-J} are direct consequences of this fact and the details of the proof are omitted.
\end{proof}

Concerning $A_{\tau, \eps}$, we  have

\begin{lemma}\label{lem-Atau} Let $0 < \eps < \gamma_0$. Assume that $C$ is analytic in $I$. We have
$$
A_{\tau, \eps} \mbox{ is analytic in $I_2$ except for a discret set, which is removable}.
$$
\end{lemma}

Recall that $A_{\tau, \eps}$ is defined in \eqref{def-YA}.

\begin{proof} Denote
$$
l = \mathop{\max_{\tau \in I_2}}_{H \mathrm{\; is \; continuous \; at \; } \tau + \eps } \dim A_{\tau, \eps} < + \infty.
$$
Fix $\tau_0 \in I_2$ such that $\dim A_{\tau_0, \eps}
 = l$ and fix $\xi_1, \cdots, \xi_l \in [L^2(0, \eps)]^m$ such that
$$
\Big\{ \mbox{Proj}_{H(\tau_0 +\eps)} \cT^c_{\tau_0, \eps} (\xi_j); 1 \le j \le l \Big\} \mbox{ is an orthogonal basis of $A_{\tau_0, \eps}$}.
$$
Since,  for fixed $\eps$, $\cT^{c}_{\cdot, \eps}$ is analytic in $I_2$ and $H(\cdot + \eps)$ is analytic in $I_2$ except for a discrete subset which is removable, it follows that
\begin{equation}\label{lem-Atau-p1}
\dim \mbox{span} \Big\{ \mbox{Proj}_{H(\tau +\eps)} \cT^c_{\tau, \eps} (\xi_j); 1 \le j \le l \Big\} = l \mbox{ in $I_2$ except for a discrete subset}.
\end{equation}
This in turn implies, by the property of $l$,
\begin{equation}\label{lem-Atau-p2}
A(\tau, \eps) = \mbox{span} \Big\{ \mbox{Proj}_{H(\tau +\eps)} \cT^c_{\tau, \eps} (\xi_j); 1 \le j \le l \Big\} \mbox{ in $I_2$ except for a discrete subset}.
\end{equation}
Combining \eqref{lem-Atau-p1} and \eqref{lem-Atau-p2} yields the conclusion.
\end{proof}

Let
\begin{equation}
\mbox{$M(\tau)$ be the orthogonal complement of $J(\tau)$ in $H(\tau)$}.
\end{equation}
It is clear that for each $\tau \in I_1$, there exists some $\eps_\tau > 0$  such that
one cannot steer any $\varphi \in M(\tau) \setminus \{0 \}$ at time $\tau$ to $0$ at time $\tau + T_{opt} + \eps_\tau$. The constant  $\eps_\tau$ can be chosen independently of $\varphi \in M(\tau) \setminus \{0\}$, for example, one can take $\eps_\tau$ so that $J(\tau, \eps) = J(\tau)$ for $0 \le \eps \le \eps_\tau/2$. The analyticity of $C$ is not required for this purpose. Nevertheless, when the analyticity of $C$ in $I$ is imposed, one can obtain a {\it uniform} lower bound for $\eps_\tau$ for $\tau \in I_2$ in a sense which will be precise now. The uniform lower bound of $\eps_\tau$ will play a crucial role in our proof of \Cref{thm-A}.  To establish this property, for $0< \eps < \gamma_0$ and $\tau \in I_2$, we first write $J(\tau, \eps)$ under the form
\begin{equation}
J(\tau, \eps)  = \Big\{\varphi \in H(\tau);  \mbox{Proj}_{A_{\tau, \eps}} \mbox{Proj}_{H(\tau + \eps)}  \cT^I_{\tau, \eps} (\varphi)  -  \mbox{Proj}_{H(\tau + \eps)}  \cT^I_{\tau, \eps} (\varphi)  = 0
\Big\}.
\end{equation}
Since the operator
\begin{multline*}
 \mbox{Proj}_{A_{\cdot, \eps}} \mbox{Proj}_{H(\cdot + \eps)}  \cT^I_{\cdot, \eps}  -  \mbox{Proj}_{H(\cdot + \eps)}  \cT^I_{\cdot, \eps}  \mbox{ is analytic in $I_2$,} \\[6pt]
 \mbox{except for a discret subset, which is removable},
\end{multline*}
one has, as in the proof of \Cref{lem-A2},
$$
J(\cdot, \eps) \mbox{ is analytic in $I_2$ except for a discret set, which is removable}.
$$
We derive that for each $n \in \mN$ with $1/n < \gamma_0$, there exists a discrete subset $D_n$ of $I_2$ such that
$$
J(\tau, 1/n)   \mbox{ is analytic in $I_2$ except for a discrete set  $D_n$, which is removable } \footnote{Replacing $\gamma_0$ by $\gamma_0/2$ if necessary, one can even assume that $D_n$ is finite.}.
$$
As a consequence, one has
\begin{equation}\label{key-fact1}
\dim J(\cdot, 1/n) \mbox{ is constant in $I_2 \setminus D_n$}.
\end{equation}

Set
\begin{equation}\label{def-D}
D = \bigcup_{n \in \mN; 1/n < \gamma_0} D_n
\end{equation}
and fix $\tau_0 \in I_2 \setminus D$. There exists $0< \eps_0 < \gamma_0$ such that
$$
J(\eps, \tau_0) = J(\tau_0) \mbox{ for } 0 < \eps < \eps_0.
$$
It follows from \Cref{lem-J} and \eqref{key-fact1} that, for $0 < \eps < \eps_0$ and $\tau, \tau' \in I_2 \setminus D$, one has
\begin{equation}\label{key-fact2}
J(\tau, \eps) = J(\tau) \quad \mbox{ and } \quad \dim J(\tau) = \dim J(\tau').
\end{equation}

We thus proved

\begin{lemma}\label{lem-JM} There exists a discrete set $D$ \footnote{The set mentioned here is the union of the set $D$ given in \eqref{def-D} and the set of $\tau \in I_2$ such that $\dim H(\tau)$ is constant, which is discrete. For notational ease, we still use the same notation $D$.} and $0< \eps_0 < \gamma_0$ such that
$$
\dim M(\tau) = \dim M(\tau') \mbox{ for } \tau, \tau' \in  I_2 \setminus D,
$$
and one cannot steer any $v \in M(\tau) \setminus \{0\}$ from time $\tau$ to 0 at time $\tau + T_{opt} + \eps_0 \mbox{ for } \tau \in I_2 \setminus D$.
\end{lemma}

We now summarize the results which have been derived in this section:

\begin{proposition} \label{pro-B}  There exist an  orthogonal decomposition of $H(\tau)$ via $H(\tau) = J(\tau) \otimes M(\tau)$ for $\tau \in I_1$, a discrete subset $D$ of $I_2$, and a constant $\eps_0 > 0$ such that the following four properties hold:
\begin{itemize}

\item[i)] For $\varphi \in J(\tau)$, one can steer $v$ at time $\tau$ to $0$ at  time $\tau + T_{opt} + \delta$ for all $\delta > 0$.

\item[ii)] For $\varphi \in M(\tau) \setminus \{0\}$, there exists $\eps_\tau > 0$ such that one cannot  steer $\varphi$ at time $\tau$ to $0$ at time $\tau + T_{opt} + \delta$ for $0<  \delta < \eps_\tau$.

\item[iii)]
$$
\dim M(\tau) = \dim M(\tau') \mbox{ for } \tau, \tau' \in I_2 \setminus D,
$$

\item[iv)] For $\tau \in I_2 \setminus D$, and  $\varphi \in M(\tau) \setminus \{0\}$, one cannot steer $\varphi$ at time $\tau$ to 0 at time $\tau + T_{opt} + \eps_0$.

\end{itemize}
\end{proposition}

\Cref{pro-B} also gives the characterization of states which can be steered at time $\tau$ to 0 at time $\tau + T_{opt} + \delta$ for all $\delta > 0$. Indeed, one has, for $\tau \in I_1$,

\begin{itemize}

\item For $v \in H(\tau) \cup J(\tau)$, one can steer $v$ at  time $\tau$ to $0$ at  time $\tau + T_{opt} + \delta$ for all $\delta > 0$.

\item For $v \in M(\tau) \setminus \{0\}$,  there exists $\eps_\tau > 0$ such that one cannot steer $v$ at time $\tau$ to 0 at  time $\tau + T_{opt} + \delta$ for $0 < \delta < \eps_\tau$.
\end{itemize}

%\begin{remark} \rm
%The uniformity of $\eps_0$ given in iv) of \Cref{pro-B} plays an important role in the proof of \Cref{thm-A}
%\end{remark}

\subsection{Null-controllability in time $T_{opt, +}$ - Proof of \Cref{thm-A}}

We first assume that $0 \not \in D$.  We will prove that $M(0) = \{0\}$ by contradiction,  and the conclusion follows from \Cref{pro-B}. Assume that there exists  $\varphi \in M(0) \setminus \{0 \}$. Since $M(0) \subset H(0, T_{opt} + \eps_0)$ by assertion iv) of \Cref{pro-B}, it follows from  \Cref{lem-AAA} that there exists a solution $v^{(0)}$ of the system
\begin{equation}\label{thmA-vp1}
\partial_t v^{(0)} (t, x) =  \Sigma(x) \partial_x v^{(0)} (t, x) + \bC (t, x)   v^{(0)}(t, x)
\mbox{ for } (t, x)  \in (0, T_{opt} + \eps_0) \times (0, 1),
\end{equation}
with, for $t \in (0, T_{opt} + \eps_0)$,
\begin{equation}\label{thmA-vp2}
v^{(0)}(t, 1)  = 0,
\end{equation}
\begin{equation}\label{thmA-vp3}
\Sigma_+ (0) v^{(0)}_+(t, 0) = - B\tr \Sigma_- (0) v^{(0)}_- (t, 0),
\end{equation}
\begin{equation}\label{thmA-vp4}
v^{(0)}(t = 0, \cdot) = \varphi \mbox{ in } (0, 1).
\end{equation}

Fix $t_1 \in (\eps_0/3, \eps_0/2) \setminus D$ (recall that $D$ is discrete). By \Cref{lem-AAA}, one has
$$
v^{(0)}(t_1, \cdot) \in H(t_1, T_{opt} + \eps_0 - t_1).
$$
This in turn implies that, since  $H(t_1, T_{opt} + \eps_0 - t_1) = M(t_1) = H(t_1, T_{opt} + \eps_0)$ by assertion iv) of \Cref{pro-B},
$$
v^{(0)}(t_1, \cdot) \in H(t_1, T_{opt} + \eps_0).
$$

By \Cref{lem-AAA} again, there exists a solution $v^{(1)}$ of the system
\begin{equation}\label{sys-v1}
\partial_t v^{(1)} (t, x) =  \Sigma(x) \partial_x v^{(1)} (t, x) + \bC(t, x)  v^{(1)}(t, x)
\mbox{ for } (t, x)  \in (t_1, t_1 + T_{opt} + \eps_0) \times (0, 1),
\end{equation}
with, for $t \in (t_1, t_1 + T_{opt} + \eps_0)$,
\begin{equation}\label{bdry1-v1}
v^{(1)}(t, 1)  = 0,
\end{equation}
\begin{equation}\label{bdry0-v1}
\Sigma_+ (0) v^{(1)}_+(t, 0) = - B\tr \Sigma_- (0) v^{(1)}_- (t, 0),
\end{equation}
\begin{equation}
v^{(1)}(t = t_1, \cdot) = v^{(0)}(t_1, \cdot)  \mbox{ in } (0, 1).
\end{equation}

Consider the solution $v$ of system \eqref{eq-vA}-\eqref{bdry-v0A}  for the time interval $(0, t_1 + T_{opt} + \eps_0)$ with  $v(t_1 + T_{opt} + \eps_0, \cdot) =  v^{(1)}(t_1 + T_{opt} + \eps_0, \cdot)$ (backward system). One can check that
$$
v (t, \cdot)= v^{(1)} (t, \cdot) \mbox{ for } t \in (t_1, t_1 + T_{opt} + \eps_0)
$$
and, since $v^{(1)} (t_1, \cdot) = v^{(0)} (t_1, \cdot)$,
$$
v (t, \cdot)= v^{(0)} (t, \cdot) \mbox{ for } t \in (0, t_1).
$$

For notational ease, we will denote this $v$ by $v^{(1)}$. We thus proved that there exists a solution $v^{(1)}$ of  \eqref{eq-vA}-\eqref{bdry-v0A}  such that
$$
v^{(1)}(\cdot, 1) = 0 \mbox{ in } (0, t_1 + T_{opt} + \eps_0), 
$$
and
$$
v^{(1)}(0, \cdot) = \varphi \mbox{ in } (0, 1).
$$

Continuing this process, there exist  $0 = t_0 < t_1 < \dots < t_{N-1} \le T_1 - T_{opt} <  t_N < \beta - T_{opt}$ and a family  of  $v^{(\ell)}$ with $1 \le \ell \le N $ such that $t_\ell \in \hI \setminus D$,
\begin{equation}
\partial_t v^{(\ell)} (t, x) =  \Sigma(x) \partial_x v^{(\ell)} (t, x) + \bC(t, x)  v^{(\ell)}(t, x)
\mbox{ for } (t, x)  \in (0, t_\ell + T_{opt} + \eps_0) \times (0, 1),
\end{equation}
with, for $t \in (0, t_\ell + T_{opt} + \eps_0)$,
\begin{equation}
v^{(\ell)}(t, 1)  = 0,
\end{equation}
\begin{equation}
\Sigma_+ (0) v^{(\ell)}_+(t, 0) = - B\tr \Sigma_- (0) v^{(\ell)}_- (t, 0),
\end{equation}
\begin{equation}
v^{(\ell)}(t = 0, \cdot) = \varphi (\cdot)\mbox{ in } (0, 1),
\end{equation}
and
$$
\eps_0 / 3 \le t_\ell - t_{\ell-1} \le \eps_0/2.
$$
This implies, by \Cref{lem-AAA},  that one cannot steer $\varphi$ from time $0$ to 0 at  time $T_1$. We have a contradiction since the system is null-controllable at the time $T_1$. The conclusion follows in the case $0 \in I_2 \setminus D$.

\medskip
The proof in the general case can be derived from the previous case
by noting that, using the same arguments, one has
$$
\mbox{ $M(\tau_0)  = \{0 \}$ for  $\tau_0 \in I_2 \setminus D$ and $\tau_0$ is close to $0$. }
$$
The details are omitted.

\medskip
The proof is complete. \qed

\medskip
The proof of \Cref{thm-A} also yields the following unique continuation principle:

\begin{proposition}\label{pro-UCP}  Let $k \ge m \ge 1$ and let $B \in \cB$ be such that \eqref{cond-B-1} holds for $i =m$.  Assume that $C_1 \in \cH(I, [L^\infty(0, 1)]^{n \times n})$. Let $\tau \in I$ and $T > T_{opt}$. Assume that  $\tau + T_1 \in I$. Let $v$ be a solution of system
\begin{equation}\label{UCP-1}
\partial_t v (t, x) =  \Sigma(x) \partial_x v (t, x) +  C_1(t, x)   v(t, x)
\mbox{ for } (t, x)  \in (\tau, \tau +  T) \times (0, 1),
\end{equation}
with, for $\tau < t < \tau + T$,
\begin{equation}\label{UCP-2}
v_-(t, 1)  = 0,
\end{equation}
\begin{equation}\label{UCP-3}
\Sigma_+ (0) v_+(t, 0) = - B\tr \Sigma_- (0) v_- (t, 0),
\end{equation}
\begin{equation}\label{UCP-4}
v_+(t, 1) = 0.
\end{equation}
Then $v = 0$.
\end{proposition}

Recall that $T_1 = \tau_k + \tau_{k+1}$, see \eqref{def-T1}. 

\begin{proof} The conclusion of \eqref{pro-UCP} follows from the proof of \Cref{thm1} applied to  $C(t, x)$ defined by  $\Sigma'(x) - C(t, x)\tr  = C_1(t, x)$.
\end{proof}

The unique continuation result stated in \Cref{pro-UCP} can be seen as a variant of the unique continuation principle for the wave equations whose first and zero-order terms  are analytic in time due to Tataru-H\"ormander-Robbiano-Zuily. Our strategy was mentioned at the beginning of \Cref{sect-thmA}.  We do not know if such a unique continuation principle can be proved using Carleman's estimate as in the wave setting. It is worth noting that if this is possible then
the analyticity of $C_1$ in time  must be taken into account by \Cref{thm1}. More importantly the conditions $B \in \cB$ and \eqref{cond-B-1} holding for $i =m$ have to be essentially used in the proof process since it is known that the unique continuation does not hold without this assumption even in the case $C_1 \equiv 0$. The advantage of Carleman's estimate might be that  the analyticity of $C_1$ is only required for a neighborhood of $[0, T_{opt}]$ instead of $[0, T_1]$.

\begin{remark}\label{rem-backstepping} \rm It is natural to compare the direct approach here with the one involving the backstepping technique. In the time-invariant setting, both approaches yield the same result since \eqref{cond-B-1} with $i=m$ is not imposed to establish the compactness of $\cK(\tau)$ (see Step 1 of the proof of \Cref{pro-A1}). Nevertheless, (equivalent) control-forms obtained from the backstepping approach are easier to handle/understand. The analysis in this paper is strongly inspired/guided by such control-forms. In the time-varying setting, one might derive the same conclusion under the assumption that $C$ is analytic in $\mR$ and   its holomorphic extension in  $\{z \in \mC; |\Im(z)| < \gamma \}$ is bounded for some $\gamma > 0$. This quite strong assumption on the analyticity of $C$ comes from the construction of the kernel in the step of using backstepping and might not be necessary.
 \end{remark}

\section{Exact controllability in the analytic setting -  Proof of  \Cref{thm-E}} \label{sect-thmE}

\Cref{thm-E} can be derived from \Cref{thm-A}, as in the proof of  \cite[Theorem 3]{CoronNg19-2}. For the convenience of the reader, we reproduce the proof.

\medskip
We first consider the case $m  = k$. Let $T > T_{opt}$ be such that $T \in I$. Set
$$
\tw (t, x) = w(T -t, x) \mbox{ for } t \in (0, T), \, x \in (0, 1).
$$
Then
$$
\tw_-(t, 0) = \widetilde B^{-1} \tw_{+}(t, 0),
$$
with $\tw_-(t, \cdot)= (w_{2k}, \dots, w_{k+1})\tr (T-t, \cdot)$, and  $\tw_+ (t, \cdot) = (w_k, \dots, w_1)\tr (T -t, \cdot)$, and $\widetilde B_{ij} = B_{pq}$ with $p = k-i$ and $q = k-j$. Note that the $i \times i$  matrix formed from the first $i$ columns and  rows of $\widetilde B$  is invertible.
Using the Gaussian elimination method, one can find  $(k \times k)$ matrices  $T_1, \dots, T_N$ such that
$$
T_N \dots T_1 \widetilde B = U,
$$
where $U$ is a $(k \times k)$ upper triangular matrix, and $T_i$ ($1 \le i \le N$) is the matrix given by  the operation which replaces a row $p$ by itself plus a multiple of a row $q$ for some $1\le q<p \le N$. It follows that
$$
\widetilde B^{-1} = U^{-1} T_N \dots T_1.
$$
One can check that $U^{-1}$ is an invertible,  upper triangular matrix,  and $T_N \dots T_1$ is an invertible,  lower triangular matrix. It follows that the
$i \times i$  matrix formed from the last $i$ columns and  rows of $\widetilde B^{-1}$ is the product of the matrix formed from the last $i$ columns and  rows of $U^{-1}$
and the matrix formed from the last $i$ columns and  rows of $T_N \dots T_1$. Therefore,  $\widetilde B^{-1} \in \cB$. One can also check  that the exact controllability of the system for $w(\cdot, \cdot)$ at the time $T$ from time $0$ is equivalent to the null-controllability  of the system for $\tw(\cdot, \cdot)$ at the same time from time $0$. The conclusion of \Cref{thm-E} now follows from \Cref{thm-A} by noting that $C(\cdot - T, \cdot)$ is analytic in a neighborhood of $[0, T_1]$.

The case $m>k$ can be obtained from the case $m = k$ as follows. Consider $\hw (\cdot, \cdot)$ the solution of the system
$$
\partial_t \hw (t, x) = \hat \Sigma(x) \partial_x \hw (t, x) + \hat C(t, x) \hw (t, x),
$$
$$
\hw_-(t, 0) = \hat B \hw_+ (t, 0),  \quad \mbox{ and } \quad \hw_+(t, 1) \mbox{ are controls}.
$$
Here
$$
\hat \Sigma = \mbox{diag} (- \hat \lambda_1, \dots, -\hat \lambda_m, \hat \lambda_{m+1}, \dots \hat \lambda_{2m}),
$$
with $\hat \lambda_j =  - (1 + m -k - j ) \eps^{-1}$  for $1 \le j \le m-k$ with positive small $\eps$, $\hat \lambda_{j} = \lambda_{j - (m-k)}$ if $ m-k + 1 \le j \le m$, and $\hat \lambda_{j+m} = \lambda_{j + k}$ for $1 \le j \le m$,
$$
\hat C(t, x) = \left(\begin{array}{cc} 0_{m-k, m-k} & 0_{m-k, n}\\[6pt]
0_{n, m-k} & C (t, x)
\end{array}\right),
$$
and
$$
\hat B = \left(\begin{array}{cc} I_{m-k} & 0_{m-k, m}\\[6pt]
0_{m-k, m} & B
\end{array}\right),
$$
where $I_{\ell}$ denotes the identity matrix of size $\ell \times \ell$ for $\ell \ge 1$. Here $0_{i, j}$ denotes the zero matrix of size $i \times j$ for $i, j, \ell \ge 1$. Then the exact controllability of $w$ at the time $T$ from time $0$ can be derived from the exact controllability of $\hw$ at the same time from time $0$. One then can deduce the conclusion of \Cref{thm-E} from the case $m=k$ using \Cref{thm-A} by noting that the optimal time for the system of $\hw$ converges to the optimal time for the system of  $w$ as $\eps \to 0_+$. \qed

\appendix

\section{Hyperbolic systems in non-rectangle domains}

In this section, we give the meaning  of broad solutions used to define  $\cT(\tau)$ and $\hat \cT(\tau)$ and study their well-posedness. We also establish the boundedness and the analyticity of $\cT(\tau)$ under appropriate assumptions. The key point of the analysis is to find suitable weighted norms in order to apply the fixed point arguments. This matter is subtle (see \Cref{rem-Weight}). In this section, we assume that $k \ge m \ge 1$ although the arguments are quite robust and also work for the case $m > k \ge 1$ under appropriate modifications.

\medskip
Let $F \in [L^\infty(\Omega)]^{n\times n}$, $(f, g) \in [L^2(0, T_{opt})]^n \times [L^2(0, 1)]^m$, and $\gamma \in [L^2(\Omega)]^n$.
We first deal with the following system, which is slightly more general than the system \eqref{w1}-\eqref{w6}:
\begin{equation}\label{wA1}
\partial_t w (t, x)  = \Sigma (x) \partial_x w (t, x)+ F(t, x) w (t, x) + \gamma(t, x) \mbox{ for } (t, x) \in \Omega,
\end{equation}
\begin{equation}\label{wA2}
w (\cdot, 1)= f  \mbox{ in } (0, T_{opt}),
\end{equation}
\begin{equation}\label{wA3}
w_+(0, \cdot) = g \mbox{ in } (0, 1),
\end{equation}
\begin{equation}\label{wA4}
w_{-, \ge k}(t, 0) = Q_k w_{<k, \ge k+m} (t, 0) \mbox{ for } t \in (T_{opt} - \tau_{k}, T_{opt} - \tau_{k-1}),
\end{equation}
\begin{equation}\label{wA5}
w_{-, \ge k-1}(t, 0) = Q_{k-1} w_{<k-1, \ge k+m-1} (t, 0) \mbox{ for } t \in  (T_{opt} - \tau_{k-1}, T_{opt} - \tau_{k-2}),
\end{equation}
\dots
\begin{equation}\label{wA6}
w_{-, \ge k-m+ 2}(t, 0) = Q_{k-m + 2} w_{<k-m+2, \ge k+2} (t, 0) \mbox{ for } t \in (T_{opt} - \tau_{k-m+2},  T_{opt} - \tau_{k - m + 1}).
\end{equation}

Given a subset  $O$ of  $\mR^2$ and a point $(t, x) \in \mR^2$,  we denote
\begin{equation*}
O_t = \Big\{ y \in \mR; (t, y) \in  O  \Big\} \quad \mbox{ and } \quad O_x = \Big\{s \in \mR ; (s, x) \in O \Big\}.
\end{equation*}

We next give the definition of the broad solutions of system \eqref{wA1}-\eqref{wA6}.

\begin{definition}\label{def-WP} Let  $F \in [L^\infty(\Omega)]^{n\times n}$, $(f, g) \in [L^2(0, T_{opt})]^n  \times [L^2(0, 1) ]^m$, and $\gamma \in [L^2(\Omega)]^n$.  A vector-valued function $w \in \cY : = \big[L^2(\Omega)\big]^n \cap C \big( [0, T_{opt}]; [L^2 (\Omega_t)]^n \big) \cap C \big( [0, 1]; [L^2 (\Omega_x)]^n \big) $ \footnote{A function $\varphi \in L^2(\Omega)$ is said to be in $C \big( [0, T_{opt}]; L^2 (\Omega_t) \big)$ if $(t_n) \subset [0, T_{opt}]$ converging to $t$ then $\lim_{n \to + \infty} \left( \|f(t_n, \cdot) - f(t, \cdot) \|_{L^2(\Omega_{t_n} \cap \Omega_t)} + \|f(t_n, \cdot) \|_{L^2(\Omega_{t_n} \setminus \Omega_t)} +  \|f(t, \cdot) \|_{L^2(\Omega_{t} \setminus \Omega_{t_n})} \right) = 0$. Similar meaning is used for $C \big( [0, 1]; L^2 (\Omega_x) \big)$.} is called a broad solution of \eqref{wA1}-\eqref{wA6} if for almost $(t_1, \xi_1) \in \Omega $, the following conditions hold
\begin{enumerate}
\item[1.]  for $ 1 \le j \le k-m+1$,
\begin{equation}\label{WP-p1}
w_j(t_1, \xi_1) = \int_{t}^{t_1} \Big( F \big(s, x_j (s, t_1, \xi_1) \big) w(s, x_j (s, t_1, \xi_1))  \Big)_j\, ds + \int_{t}^{t_1} \gamma_j \big(s, x_j (s, t_1, \xi_1) \big) \, ds  + f_j(t),
\end{equation}
where $t$ is such that $ x_j(t, t_1, \xi_1) = 1$;

\item[2.] for $ k - m + 2 \le j \le k$,
\begin{equation}\label{WP-p2}
w_j(t_1, \xi_1) = \int_{t}^{t_1} \Big( F \big(s, x_j (s, t_1, \xi_1) \big) w(s, x_j (s, t_1, \xi_1))  \Big)_j\, ds +  \int_{t}^{t_1} \gamma_j \big(s, x_j (s, t_1, \xi_1) \big) \, ds + f_j (t),
\end{equation}
if $t \in (0, T_{opt})$ where  $t$ is such that  $ x_j(t, t_1, \xi_1) = 1$,  otherwise,
\begin{multline}\label{WP-p3}
w_j(t_1, \xi_1) = \int_{\hat t}^{t_1} \Big( F \big(s, x_j (s, t_1, \xi_1) \big) w(s, x_j (s, t_1, \xi_1))  \Big)_j\, ds  \\[6pt]
+ \int_{\hat t}^{t_1} \gamma_j \big(s, x_j (s, t_1, \xi_1) \big) \, ds + \Big( Q_l w_{<l, \ge l+m}(\hat t, 0) \Big)_{j - l + 1},
\end{multline}
if $\hat t \in (T_{opt} - \tau_{l}, T_{opt} - \tau_{l-1})$ where $\hat t$ is  such that $ x_j(\hat t, t_1, \xi_1) = 0$.

\item[3.] for $k+1 \le j \le k+m$,
\begin{equation}\label{WP-p4}
w_j(t_1, \xi_1) = \int_{t}^{t_1} \Big( F \big(s, x_j (s, t_1, \xi_1) \big) w(s, x_j (s, t_1, \xi_1))  \Big)_j\, ds + \int_{t}^{t_1} \gamma_j \big(s, x_j (s, t_1, \xi_1) \big) \, ds  + f_j (t),
\end{equation}
if $t \in (0, T_{opt})$ where $t$ is such that $ x_j(t, t_1, \xi_1) = 1$,  otherwise
\begin{equation}\label{WP-p5}
w_j(t_1, \xi_1) = \int_0^{t_1} \Big( F \big(s, x_j (s, t_1, \xi_1) \big) w(s, x_j (s, t_1, \xi_1))  \Big)_j\, ds + \int_{0}^{t_1} \gamma_j \big(s, x_j (s, t_1, \xi_1) \big) \, ds + g_{j-k} (\eta),
\end{equation}
where $\eta \in (0, 1)$ is such that $ x_j(0, t_1, \xi_1) =  \eta$.
\end{enumerate}
\end{definition}

Recall that the characteristic flow $x_j$ with $1 \le j \le k+m$ is defined in \eqref{def-xi-1} and \eqref{def-xi-2}.

In this definition, the term $Q_l w_{<l, \ge l+m}(\hat t, 0)$ in \eqref{WP-p3} is required to be replaced by the corresponding expression in the RHS of \eqref{WP-p1}, or \eqref{WP-p2}, or \eqref{WP-p4}, or \eqref{WP-p5} with $(\hat t, 0)$ standing for  $(t_1, \xi_1)$.

\medskip
The well-posedness of broad solutions  of \eqref{wA1}-\eqref{wA6} is given in the following.
\begin{theorem}\label{thm-WP} Let $F \in \big[L^\infty \big(\Omega\big) \big]^{n \times n}$, $(f, g) \in  \big[L^2(0, T_{opt}) \big]^n \times \big[L^2(0, 1) \big]^m$, and $\gamma \in [L^2(\Omega)]^n$.  There exists a unique broad solution   $w \in \cY$  of \eqref{wA1}-\eqref{wA6}. Moreover,
\begin{equation}\label{thm-WP-e1}
\| w\|_{\cY} \le C \Big( \| f\|_{L^2(0, T_{opt})} + \| g\|_{L^2(0, 1)} + \|\gamma \|_{L^2(\Omega)} \Big),
\end{equation}
for some positive constant $C$ depending on an upper bound of $\| F \|_{L^\infty(\Omega)}$ and $\Sigma$.
\end{theorem}

Here we denote
$$
\| w\|_{\cY} =  \max\left\{ \sup_{x \in [0, 1]} \| w \|_{L^2(\Omega_{ x})}, \sup_{t \in [0, T_{opt}]} \| w  \|_{L^2(\Omega_{t})}; 1 \le i \le n \right\}.
$$

\begin{remark} \rm The analysis of \Cref{thm-WP} can be easily extended to cover the case where source terms in $L^2$ are added in \eqref{wA4}-\eqref{wA6}.
\end{remark}

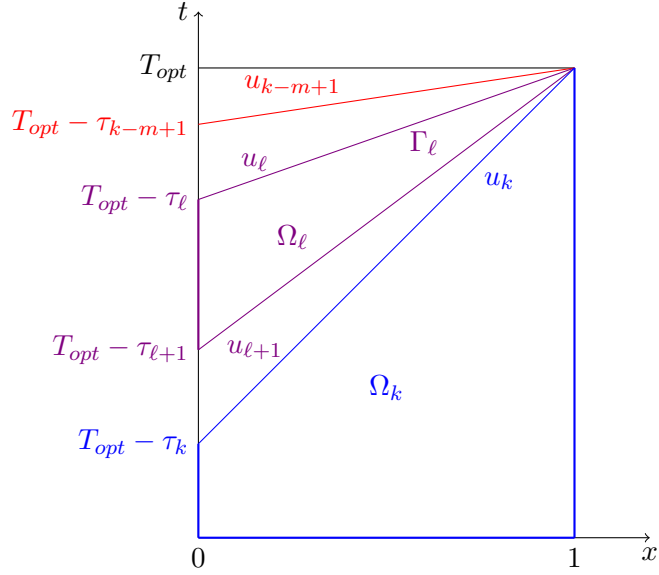
\begin{figure}
\centering
\begin{tikzpicture}[scale=2.5]

\newcommand\z{0.0}

\draw[->] (0.0 + \z,0) -- (2.4+\z,0);
\draw[->] (0.0 +\z,0) -- (0.0 +\z,2.8);
\draw[] (2+\z,0) -- (2+\z,2.5);

\draw (2+\z, 0) node[below]{$1$};

\draw (\z, 2.8) node[left]{$t$};

\draw (2.4+\z, 0) node[below]{$x$};

\draw (\z, 0) node[left, below]{$0$};

\draw[] (\z, 2.5) -- (\z + 2, 2.5);

\draw (\z, 2.5) node[left]{$T_{opt}$};

\draw[red] (\z + 2, 2.5) -- (\z , 2.2);
\draw[red]  (\z , 2.2) node[left]{$T_{opt} - \tau_{k-m+1}$};
\draw[red]  (\z + 0.5, 2.3) node[above]{$u_{k-m+1}$};

\draw[violet] (\z + 2, 2.5) -- (\z , 1.8);
\draw[violet]  (\z , 1.8) node[left]{$T_{opt} - \tau_{\ell}$};
\draw[violet]  (\z + 0.3, 1.9) node[above]{$u_{\ell}$};

\draw[violet] (\z + 2, 2.5) -- (\z , 1.);
\draw[violet]  (\z , 1.) node[left]{$T_{opt} - \tau_{\ell+1}$};
\draw[violet]  (\z + 0.3, 1.) node[]{$u_{\ell+1}$};
\draw[violet, thick] (\z, 1.) -- (\z , 1.8);
\draw[violet]  (\z + 0.5, 1.6) node[]{$\Omega_\ell$};

\draw[violet] (\z+1.2, 2.1) node[]{$\Gamma_\ell$};

\draw[blue] (\z + 2, 2.5) -- (\z , 0.5);
\draw[blue]  (\z , 0.5) node[left]{$T_{opt} - \tau_{k}$};
\draw[blue]  (\z + 1.6, 2.) node[below]{$u_{k}$};
\draw[blue]  (\z + 1, 0.8) node[]{$\Omega_k$};
\draw[blue, thick] (\z, 0) -- (\z , 0.5);
\draw[blue, thick] (\z, 0) -- (\z + 2, 0);
\draw[blue, thick] (\z + 2, 0) -- (\z + 2, 2.5);

\end{tikzpicture}
\caption{Geometry of $\Omega_\ell$ and $\Gamma_\ell$ with $k-m +1 \le \ell \le k$ for a constant $\Sigma$.}\label{fig-Omega}
\end{figure}

Before giving the proof of \Cref{thm-WP}, let us introduce some notations. For $k-m + 1 \le \ell \le k-1$, let $\Omega_\ell$ be the region of $\Omega$ between the characteristic curves of $x_\ell$ and $x_{\ell +1}$ both passing the point $(T_{opt}, 1)$ in the $xt$-plane.  We also denote $\Omega_k$ the region of $\Omega$ below the characteristic curve of $x_k$ passing the point $(T_{opt}, 1)$ in the $xt$-plane. Let $\Gamma_\ell$ with $k-m+1 \le \ell \le k$ be the boundary part of $\Omega_\ell$ formed by the characteristic curve of $x_\ell$ passing the point $(T_{opt}, 1)$. See \Cref{fig-Omega}.

\medskip
The proof of \Cref{thm-WP} is based on two lemmas below. The first one is the following.

\begin{lemma}\label{lem-WP1} Let $F \in \big[L^\infty (\Omega_k) \big]^{n \times n}$, $(f, g) \in  [L^2(0, T_{opt})]^n \times [L^2(0, 1)]^m$, and $\gamma \in [L^2(\Omega_k)]^n$.  There exists a unique board solution
$w \in \cY_k: = [L^2(\Omega_k)]^n \cap  C \big( [0, T_{opt}]; [L^2 (\Omega_{k, t}) ]^n \big) \cap C \big( [0, 1]; [L^2 (\Omega_{k, x}) ]^n \big) $ of the system
\begin{equation}\label{WP1-wA1}
\partial_t w (t, x)  = \Sigma (x) \partial_x w (t, x)+ F(t, x) w (t, x) + \gamma(t, x) \mbox{ for } (t, x) \in \Omega_k,
\end{equation}
\begin{equation}\label{WP1-wA2}
w (\cdot, 1)= f \mbox{ in } (0, T_{opt}),
\end{equation}
\begin{equation}\label{WP1-wA3}
w_+(0, \cdot) = g \mbox{ in } (0, 1).
\end{equation}
Moreover,
\begin{equation}\label{lem-WP1-stability}
\| w\|_{\cY_k} \le C \Big( \| f\|_{L^2(0, T_{opt})} + \| g\|_{L^2(0, 1)} + \| \gamma\|_{L^2(\Omega_k)} \Big),
\end{equation}
for some positive constant $C$ depending only on an upper bound of $\| F \|_{L^\infty(\Omega_k)}$ and $\Sigma$.
\end{lemma}

Here we denote
$$
\| w\|_{\cY_k} =  \max\left\{ \sup_{x \in [0, 1]} \| w \|_{L^2(\Omega_{k, x})}, \sup_{t \in [0, T_{opt}]} \| w  \|_{L^2(\Omega_{k, t})}; 1 \le i \le n \right\}.
$$

The broad solutions considered in \Cref{lem-WP1} are defined similarly as  the one of \Cref{def-WP} as follows:

\begin{definition}\label{def-WP1} Let $F \in [L^\infty (\Omega_k) ]^{n \times n}$, and  $(f, g) \in [L^2(0, T_{opt})]^n \times [L^2(0, 1)]^m$, and $\gamma \in [L^2(\Omega_k)]^n$.  A vector-valued function $w \in \cY_k$ is called a broad solution of \eqref{WP1-wA1}-\eqref{WP1-wA3} if for almost $(t_1, \xi_1) \in \Omega_k $,  the following conditions hold
\begin{enumerate}
\item[1.]  for $ 1 \le j \le k$,
\begin{equation}\label{def-WP-A1-p1}
w_j(t_1, \xi_1) = \int_{t}^{t_1} \Big( F \big(s, x_j (s, t_1, \xi_1) \big) w(s, x_j (s, t_1, \xi_1))  \Big)_j\, ds +  \int_{t}^{t_1} \gamma_j \big(s, x_j (s, t_1, \xi_1) \big)  \, ds + f_j(t),
\end{equation}
where $t$ is such that $ x_j(t, t_1, \xi_1) = 1$;

\item[2.] for $k+1 \le j \le k+m$,
\begin{equation}\label{def-WP-A1-p2}
w_j(t_1, \xi_1) = \int_{t}^{t_1} \Big( F \big(s , x_j (s, t_1, \xi_1) \big) w(s, x_j (s, t_1, \xi_1))  \Big)_j\, ds +  \int_{t}^{t_1} \gamma_j \big(s, x_j (s, t_1, \xi_1) \big)  \, ds + f_j (t),
\end{equation}
if $t \in (0, T_{opt})$ where $t$ is such that $ x_j(t, t_1, \xi_1) = 1$,  otherwise
\begin{equation}\label{def-WP-A1-p3}
w_j(t_1, \xi_1) = \int_0^{t_1} \Big( F  \big(s, x_j (s, t_1, \xi_1) \big) w(s, x_j (s, t_1, \xi_1))  \Big)_j\, ds +  \int_0^{t_1} \gamma_j \big(s, x_j (s, t_1, \xi_1) \big)  \, ds + g_{j-k} (\eta),
\end{equation}
where $\eta \in (0, 1)$ is such that $ x_j(0, t_1, \xi_1) =  \eta$.
\end{enumerate}
\end{definition}

\begin{proof}[Proof of \Cref{lem-WP1}]  For $v \in [L^2(\Omega_{k})]^n$, set
$$
 T_k(v)(t, x) = e^{L x} v (t, x)  \mbox{ for } (t, x) \in \Omega_k,
$$
where $L$ is a large positive constant determined later.

We now introduce
$$
\| v \|_{\Omega_k} : = \max\left\{ \sup_{x \in [0, 1]} \| (T_kv)_i \|_{L^2(\Omega_{k, x})}, \sup_{t \in [0, T_{opt}]} \| (Tv)_i \|_{L^2(\Omega_{k, t})}; 1 \le i \le n \right\}.
$$
One can check that $\cY_k$ equipped with the norm $\| \cdot \|_{\Omega_k}$ is a Banach space. It is also clear that $\| \cdot \|_{\Omega_k}$ is equivalent to $\| \cdot \|_{\cY_k}$.

The proof is now based on a fixed point argument. To this end, define $\cF_k$ from $\cY_k$  into itself as follows: for $v \in \cY_k$,  and for $(t_1, \xi_1) \in \Omega_k$ and $1 \le j \le k+m$:
\begin{multline}\label{WP1-def-Fk}
\big(\cF_k(v) \big)_j(t_1, \xi_1) \mbox{ is the RHS of \eqref{def-WP-A1-p1},  or \eqref{def-WP-A1-p2}, or \eqref{def-WP-A1-p3} } \\
\mbox{ under the corresponding conditions.}
\end{multline}

We claim that, for $L$ large enough, $\cF_k$ is a contraction mapping  from $\cY_k$ equipped with the norm $\| \cdot \|_{\Omega_k}$ into itself;  and the conclusion follows then.

For $v \in \cY_k$, one can check that $\cF(v) \in \cY_k$.

Let $v, w \in \cY_k$ be arbitrary.  Fix $\xi_1 \in [0, 1]$. Let $1 \le j \le k$. We have for $(t_1, \xi_1) \in \Omega_k$, by \eqref{def-WP-A1-p1},
\begin{equation}\label{lem-WP1-p0}
  \cF(v)_j (t_1, \xi_1) -  \cF(w)_j(t_1, \xi_1) =
\int_{t}^{t_1} \Big( F \big(s, x_j (s, t_1, \xi_1) \big) (v - w) (s, x_j (s, t_1, \xi_1))  \Big)_j\, ds,
\end{equation}
where $t = t(t_1, \xi_1)$ is such that $ x_j(t, t_1, \xi_1) = 1$. This implies
\begin{multline*}
\int_{\Omega_{k, \xi_1}} e^{2L \xi_1} | \cF(v)_j (t_1, \xi_1) -  \cF(w)_j(t_1, \xi_1)|^2 \, d t_1 \\[6pt]
\le  C \int_{\Omega_{k, \xi_1}} \sign(t-t_1) \int_{t_1}^t  e^{2L \xi_1} |v - w|^2 \big(s, x_j (s, t_1, \xi_1) \big) \, ds \, d t_1,
\end{multline*}
where $\sign(\theta) = 1$ if $\theta>0$ and $-1$ if $\theta<0$.  Here and in what follows in this proof, $C$ denotes a positive constant which depends only on an upper bound of $\| F \|_{L^\infty(\Omega_k)}$ and $\Sigma$,  and can change from one place to another.

Since
$$
e^{2L \xi_1} |v - w|^2 \big(s, x_j (s, t_1, \xi_1) \big)  =  e^{2 L \big(\xi_1  -x_j (s, t_1, \xi_1) \big) } e^{2 L  x_j (s, t_1, \xi_1)}  |v - w|^2 \big(s, x_j (s, t_1, \xi_1) \big) ,
$$
and, for $s$ between $t_1$ and $t$,
$$
\xi_1  - x_j (s, t_1, \xi_1) \le 0,
$$
by a change of variables $x = x_j (s, t_1, \xi_1)$  \footnote{$x_j$ is continuously differentiable with respect to  $s, t_1, \xi_1$ when  $x_j (s, t_1, \xi_1)$ is in $\bar \Omega$ since $\Sigma$ is of class $C^2$.}, one obtains, for $1 \le j \le k$,
\begin{multline}\label{lemWP1-p1}
\int_{\Omega_{k, \xi_1}} e^{2L \xi_1} | \cF(v)_j (t_1, \xi_1) -  \cF(w)_j(t_1, \xi_1)|^2 \, d t_1 \\[6pt] \le  C \int_{\Omega_k; x \ge \xi_1} e^{2L (\xi_1 - x)} e^{2 L x} |v - w|^2 (s, x) \, ds \, d x \le \frac{C}{L} \| v- w \|_{\Omega_k}^2.
\end{multline}

We next consider $k+1 \le j \le k+ m$. Using \eqref{def-WP-A1-p2} and \eqref{def-WP-A1-p3}, similar to \eqref{lemWP1-p1}  for $1 \le j \le k$, we also reach \eqref{lemWP1-p1} for $k+1 \le j \le k+m$. Combining this with \eqref{lemWP1-p1}  for $1 \le j \le k$ yields
\begin{equation}\label{lemWP1-p2}
\int_{\Omega_{k, \xi_1}} e^{2 L \xi_1} | \cF(v) (t_1, \xi_1) -  \cF(w)(t_1, \xi_1)|^2 \, d t_1  \le \frac{C}{L} \| v- w \|_{\Omega_k}^2.
\end{equation}

Fix $t_1 \in [0, T_{opt}]$.  Let $1 \le j \le k$. From \eqref{lem-WP1-p0}, we obtain, for $(t_1, \xi_1) \in \Omega_k$,
\begin{multline*}
\int_{\Omega_{k, t_1}} e^{2 L \xi_1} | \cF(v)_j (t_1, \xi_1) -  \cF(w)_j(t_1, \xi_1)|^2 \, d \xi_1 \\[6pt]
\le  C \int_{\Omega_{k, t_1}} \sign(t-t_1) \int_{t_1}^t  e^{2 L \xi_1} |v - w|^2 (s, x_j (s, t_1, \xi_1)) \, ds \, d t_1.
\end{multline*}
Similar to \eqref{lemWP1-p1}, we obtain, for $1 \le j \le k$,
\begin{multline}\label{lemWP1-p3}
\int_{\Omega_{k, t_1}} e^{2 L \xi_1} | \cF(v)_j (t_1, \xi_1) -  \cF(w)_j(t_1, \xi_1)|^2 \, d \xi_1 \\[6pt] \le  C \int_{\Omega_k; x \ge \xi_1} e^{2 L (\xi_1 - x)} e^{2 L x} |v - w|^2 (s, x) \, ds \, d t_1 \le \frac{C}{L} \| v- w \|_{\Omega_k}^2.
\end{multline}
Using \eqref{def-WP-A1-p2} and \eqref{def-WP-A1-p3}, similar to \eqref{lemWP1-p3}  for $1 \le j \le k$, we also reach \eqref{lemWP1-p3} for $k+1 \le j \le k+m$. Combining this with  \eqref{lemWP1-p3}  for $1 \le j \le k$ yields
\begin{equation}\label{lemWP1-p4}
\int_{\Omega_{t_1}} e^{- 2L \xi_1} | \cF(v) (t_1, \xi_1) -  \cF(w)(t_1, \xi_1)|^2 \, d \xi_1  \le \frac{C}{L} \| v- w \|_{\Omega_k}^2.
\end{equation}

The claim now follows from \eqref{lemWP1-p2} and \eqref{lemWP1-p4}. The proof is complete.
\end{proof}

The second lemma used in the proof of \Cref{thm-WP} is the following.

\begin{lemma}\label{lem-WP2}  Let $k-m+1 \le \ell \le k-1$,
$F \in \big[L^\infty (\Omega_\ell) \big]^{n \times n}$, $\gamma \in [L^2(\Omega_\ell)]^n$, and  $h_j \in L^2(\Gamma_{\ell+1})$ for $1 \le j   \le k+m$ and $j \neq \ell + 1$. There exists a unique board solution $w \in \cY_\ell: = [L^2(\Omega_\ell) ]^n \cap C \big( [0, T_{opt}]; [L^2 (\Omega_{\ell, t})]^n \big) \cap C \big( [0, 1]; [L^2 (\Omega_{\ell, x})]^n \big) $ of the system
\begin{equation}\label{WP2-wA1}
\partial_t w (t, x)  = \Sigma (x) \partial_x w (t, x)+ F(t, x) w (t, x) + \gamma(t, x) \mbox{ for } (t, x) \in \Omega_\ell,
\end{equation}
\begin{equation}\label{WP2-wA2}
w_j  =  h_j \mbox{ on } \Gamma_{\ell+1}, \mbox{ for $1 \le j \le k+m$ and $j \neq \ell + 1$},
\end{equation}
\begin{equation}\label{WP2-wA3}
w_{-, \ge \ell + 1}(0, \cdot) = Q_{\ell+1} w_{< \ell + 1, \ge k+\ell + 1} \mbox{ for } t \in (T_{opt} - \tau_{\ell+1}, T_{opt} - \tau_\ell).
\end{equation}
Moreover,
$$
\| w\|_{\cY_\ell} \le  C \left(  \sum_{1 \le j \le k+m; j \neq \ell + 1} \|  h_j\|_{L^2(\Gamma_\ell)} + \| \gamma\|_{L^2(\Omega_\ell)} \right),
$$
for some positive constant $C$ depending only on an upper bound of $\| F \|_{L^\infty(\Omega_\ell)}$ and $\Sigma$.
\end{lemma}

\begin{remark} \rm The analysis of \Cref{lem-WP2} can be easily extended  to cover the case where source terms in $L^2$ are added in \eqref{WP2-wA3}.
\end{remark}

The broad solutions considered in \Cref{lem-WP2}, which are in the same spirit of the ones in \Cref{lem-WP1},  are defined as follows:

\begin{definition}\label{def-WP2}   Let $k-m+1 \le \ell \le k-1$,
$F \in \big[L^\infty (\Omega_\ell) \big]^{n \times n}$, $\gamma \in [L^2(\Omega_\ell)]^n$, and  $h_j \in L^2(\Gamma_{\ell+1})$ for $1 \le j  \le  k+m$ and $j \neq \ell + 1$.  A vector-valued function $w \in \cY_\ell$ is called a broad solution of \eqref{WP2-wA1}-\eqref{WP2-wA3} if for almost $(t_1, \xi_1) \in \Omega_\ell $, the following conditions hold
\begin{enumerate}
\item[1.]  for $ 1 \le j \le \ell$ and for $k+1 \le j \le k+m$,
\begin{equation}\label{def-WP-A2-p1}
w_j(t_1, \xi_1) = \int_{t}^{t_1} \Big( F \big(s, x_j (s, t_1, \xi_1) \big) w(s, x_j (s, t_1, \xi_1))  \Big)_j\, ds + \int_{t}^{t_1}  \gamma_j \big(s, x_j (s, t_1, \xi_1) \big) \, ds  + h_j(t),
\end{equation}
where $t$ is such that $ x_j(t, t_1, \xi_1) \in \Gamma_{\ell+1}$;

\item[2.] for $\ell + 1 \le j \le k$,
\begin{multline}\label{def-WP-A2-p2}
w_j(t_1, \xi_1) = \int_{\hat t}^{t_1} \Big( F \big(s, x_j (s, t_1, \xi_1) \big) w(s, x_j (s, t_1, \xi_1))  \Big)_j\, ds \\[6pt]+ \int_{\hat t}^{t_1}  \gamma_j \big(s, x_j (s, t_1, \xi_1) \big) \, ds   + \big(Q_{\ell+1} w_{< \ell + 1, \ge \ell + m + 1} \big)_{j - \ell}(\hat t, 0)
\end{multline}
if $\hat t \in (T_{opt} - \tau_{\ell +1}, T_{opt} - \tau_\ell)$ where $\hat t $ is such that $ x_j(\hat t, t_1, \xi_1) = 0$,  otherwise,
\begin{equation}\label{def-WP-A2-p3}
w_j(t_1, \xi_1) = \int_{t}^{t_1} \Big( F \big(s, x_j (s, t_1, \xi_1) \big) w(s, x_j (s, t_1, \xi_1))  \Big)_j\, ds + \int_{t}^{t_1}  \gamma_j \big(s, x_j (s, t_1, \xi_1) \big) \, ds  + h_j(t),
\end{equation}
where $t$ is such that $ x_j(t, t_1, \xi_1) \in \Gamma_{\ell+1}$.
\end{enumerate}
\end{definition}

As in  \Cref{def-WP}, the term $Q_{\ell + 1} w_{<\ell+1, \ge \ell + m + 1}(\hat t, 0)$ in \eqref{def-WP-A2-p2} is required to be replaced by the corresponding expression in the RHS of \eqref{def-WP-A2-p2} or  \eqref{def-WP-A2-p3}  with $(\hat t, 0)$ standing for  $(t_1, \xi_1)$.

\begin{proof}[Proof of \Cref{lem-WP2}] The key part of the proof is to introduce an appropriate weighted norm, which is adapted to the geometry and the boundary conditions considered,  for which the fixed point argument works (see \Cref{rem-Weight} for comments on this point).

We begin with the case where $\Sigma$ is constant.  For $1 \le j \le k + m$, let  $\vec{v}_j$ be the unit vector parallel to the characteristic curve of $x_j$ directed to the boundary for which  the boundary condition for  $v_j$ is given ($\vec{v}_j$ is parallel to $(1, \Sigma_{jj})\tr$ in the $xt$-plane). Set
$$
G_1 = \Big\{\vec{v}_j; 1 \le j \le \ell, k +1 \le j \le k+m \Big\} \quad \mbox{ and } \quad G_2 = \Big\{\vec{v}_j; \ell + 1 \le j \le  k \Big\}.
$$
Here are some useful observations. There exist two non-zero vectors $\vec{u}_1$ and $\vec{u}_2$ such that
\begin{enumerate}
\item[a1)] $G_1 \cup G_2 \cup \{\vec{u}_1 \}$ lies strictly on one side of the line containing $\vec{u}_2$
\item[a2)] $G_1$ is a subset of  the open, solid,  cone centered at the origin and formed by  $\vec{u}_1$ and $\vec{u}_2$,  i.e.,  in the set $\big\{s_1 \vec{u}_1 + s_2 \vec{u}_2; s_1, s_2 > 0 \big\}$.

\item[a3)] $G_2$ is a subset of the open, solid,  cone centered at the origin and  formed by  $\vec{u}_1$ and $-\vec{u}_2$,  i.e.,  in the set $\big\{s_1 \vec{u}_1 - s_2 \vec{u}_2; s_1, s_2 >  0 \big\}$.
\end{enumerate}
(For example, one can choose $\vec{u}_1 = (0, -1)\tr$ and $\vec{u_2}$ is close to $\vec{v}_\ell$ but with a larger slope in the $xt$-plane, see \Cref{fig-uv}.)
\begin{figure}
\centering
\begin{tikzpicture}[scale=2.5]

\newcommand\z{0.0}

\draw[violet] (\z + 2, 2.5) -- (\z , 1.8);
\draw[violet]  (\z , 1.8) node[left]{$T_{opt} - \tau_{\ell}$};
%\draw[violet]  (\z + 0.3, 1.9) node[above]{$u_{\ell}$};

\draw[violet] (\z + 2, 2.5) -- (\z , 1.);
\draw[violet]  (\z , 1.) node[left]{$T_{opt} - \tau_{\ell+1}$};

%\draw[violet]  (\z , 1.) node[left]{$T_{opt} - \tau_{\ell+1}$};
%\draw[violet]  (\z + 0.3, 1.) node[]{$u_{\ell+1}$};
\draw[violet, thick] (\z, 1.) -- (\z , 1.8);
%\draw[violet]  (\z + 0.5, 1.6) node[]{$\Omega_\ell$};

\draw[violet] (\z+1., 2.15) node[above]{$\Gamma_\ell$};
\draw[violet] (\z+1., 1.6) node[below]{$\Gamma_{\ell+1}$};

\newcommand\zzx{0.5}
\newcommand\zzy{2.0}

\draw[orange,->] (\zzx, \zzy) -- ({\zzx + 2/sqrt(0.7*0.7 + 2*2)}, {\zzy + 0.7/sqrt(0.7*0.7 + 2*2)});
\draw[orange] ({\zzx + 1/sqrt(0.7*0.7 + 2*2)}, {\zzy + 0.35/sqrt(0.7*0.7 + 2*2)}) node[below]{$\vec{v}_\ell$};

\newcommand\zzzx{1.4}
\newcommand\zzzy{2.0}
\draw[violet,->] (\zzzx, \zzzy) -- ({\zzzx - 2/sqrt(1.5*1.5+ 2.*2.)}, {\zzzy - 1.5/sqrt(1.5*1.5+ 2.*2.)});
\draw[violet]  ({\zzzx - 1/sqrt(1.5*1.5+ 2.*2.)}, {\zzzy - 0.6/sqrt(1.5*1.5+ 2.*2.)}) node[above, left]{$\vec{v}_{\ell+1}$};

\newcommand\zx{3.0}
\newcommand\zy{2.0}

\draw[orange,->] (\zx, \zy) -- ({\zx + 2/sqrt(0.7*0.7 + 2*2)}, {\zy + 0.7/sqrt(0.7*0.7 + 2*2)});
\draw[orange] ({\zx + 1.8/sqrt(0.7*0.7 + 2*2)}, {\zy + 0.63/sqrt(0.7*0.7 + 2*2)}) node[below]{$\vec{v}_\ell$};

\draw[orange, ->] (\zx, \zy) -- ({\zx + cos(5)}, {\zy + sin(5)});
\draw[orange] ({\zx + 0.9*cos(5)}, {\zy + 0.9*sin(5)}) node[below]{$\vec{v}_{1}$};

\draw[orange, ->] (\zx, \zy) -- ({\zx + cos(-30)}, {\zy + sin(-30)});
\draw[orange] ({\zx + 0.9*cos(-30)}, {\zy + 0.9*sin(-30)}) node[below]{$\vec{v}_{k+m}$};

\draw[orange, ->] (\zx, \zy) -- ({\zx + cos(-60)}, {\zy + sin(-60)});
\draw[orange] ({\zx + 0.9*cos(-60)}, {\zy + 0.9*sin(-60)}) node[left]{$\vec{v}_{k+1}$};

\draw[blue, ->] (\zx, \zy) -- ({\zx }, {\zy -1.45});
\draw[blue] ({\zx}, {\zy -1.30}) node[left]{$\vec{u}_{1}$};

\draw[violet,->] (\zx, \zy) -- ({\zx - 2/sqrt(1.5*1.5+ 2.*2.)}, {\zy - 1.5/sqrt(1.5*1.5+ 2.*2.)});
\draw[violet]  ({\zx - 1.8/sqrt(1.5*1.5+ 2.*2.)}, {\zy - 1.35/sqrt(1.5*1.5+ 2.*2.)}) node[below,right]{$\vec{v}_{\ell+1}$};

\draw[violet, ->] (\zx, \zy) -- ({\zx + cos(-120)}, {\zy + sin(-120)});
\draw[violet] ({\zx + 0.9*cos(-120)}, {\zy + 0.9*sin(-120)}) node[right]{$\vec{v}_{k}$};

\draw[]  ({\zx - 1.6*cos(25)}, {\zy - 1.6*sin(25)}) -- ({\zx + 1.6*cos(25)}, {\zy + 1.6*sin(25)});

\draw[red, thick, ->]  (\zx, \zy) -- ({\zx + 1.4*cos(25)}, {\zy + 1.4*sin(25)});

\draw[red]  ({\zx + 1.26*cos(25)}, {\zy + 1.26*sin(25)}) node[above]{$\vec{u}_2$};

%\draw[violet] (1., 0) -- (1, {sqrt (2)});

\end{tikzpicture}
\caption{Geometry of $\vec{v}_j$ for $1 \le j \le n$,  and $\vec{u}_1$  and $\vec{u}_2$ for $\Omega_\ell$ when $\Sigma$ is constant.}\label{fig-uv}
\end{figure}
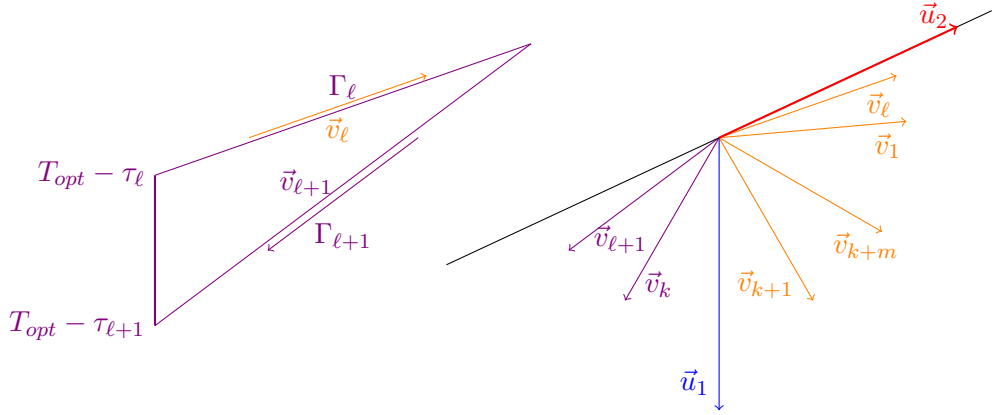

We are ready to introduce the weighted norm used. For $v \in [L^2(\Omega_{\ell})]^n$, set
\begin{equation}\label{lem-WP2-w1}
 T_\ell (v)(t, x) = e^{L  y_1(t, x)} v (t, x)  \mbox{ for } (t, x) \in \Omega_\ell,
\end{equation}
where $y_1 (t, x) $ is the first component of $(y_1, y_2)(t, x)$ which  is the coordinate of $(t, x)$ corresponding to the basis $\vec{u}_1$ and $\vec{u}_2$ (in the $xt$-plane).

We now introduce
\begin{equation}\label{lem-WP2-w}
\| v \|_{\Omega_\ell} : = \max\Big\{ \sup_{x \in [0, 1]} \| (T_\ell v)_i \|_{L^2(\Omega_{\ell, x})}, \sup_{t \in [0, T_{opt}]} \| (T_\ell v)_i \|_{L^2(\Omega_{\ell, t})}; 1 \le i \le n \Big\}.
\end{equation}
One can check that $\cY_\ell$ equipped with the norm $\| \cdot \|_{\Omega_\ell}$ is a Banach space.  It is also clear that $\| \cdot \|_{\Omega_\ell}$ is equivalent to $\| \cdot \|_{\cY_\ell}$.

The proof is now based on a fixed point argument as in the one of \Cref{lem-WP1}.  To this end, define $\cF_\ell$ from $\cY_\ell$ equipped with the norm $\| \cdot \|_{\Omega_\ell}$ into itself as follows: for $v \in \cY_\ell$ and for $(t_1, \xi_1) \in \Omega_\ell$:
\begin{multline}\label{WP1-def-Fl}
\big(\cF_\ell(v) \big)_i(t_1, \xi_1) \mbox{ is the RHS of \eqref{def-WP-A2-p1},  or \eqref{def-WP-A2-p2},  or \eqref{def-WP-A2-p3} } \\[6pt] \mbox{under the corresponding conditions.}
\end{multline}

 Fix $\xi_1 \in [0, 1]$. Let $1 \le j \le \ell$ or $k+1 \le j \le k+m$. We have, for $(t_1, \xi_1) \in \Omega_\ell$, by \eqref{def-WP-A2-p1},
\begin{equation}\label{lem-WP2-p0}
  \cF(v)_j (t_1, \xi_1) -  \cF(w)_j(t_1, \xi_1) =
\int_{t}^{t_1} \Big( F \big(s, x_j (s, t_1, \xi_1) \big) (v - w) \big(s, x_j(s, t_1, \xi_1) \big)  \Big)_j\, ds,
\end{equation}
where $t$ is such that $ x_j(t, t_1, \xi_1) \in \Gamma_{\ell + 1}$. This implies
\begin{multline}\label{lem-WP2-cc}
\int_{\Omega_{\ell, \xi_1}} e^{2 L y_1 (t_1, \xi_1)} | \cF(v)_j (t_1, \xi_1) -  \cF(w)_j(t_1, \xi_1)|^2 \, d t_1 \\[6pt] \le  C \int_{\Omega_{\ell, \xi_1}} \sign(t - t_1) \int_{t_1}^t  e^{2 L y_1(t_1, \xi_1)} |v - w|^2 \big(s, x_j (s, t_1, \xi_1) \big) \, ds \, d t_1.
\end{multline}
Here and in what follows in this proof, $C$ (resp. $c$) denotes a positive constant which depends only on an upper bound of $\| F \|_{L^\infty(\Omega_k)}$ and $\Sigma$ (resp. $\Sigma$),  and can change from one place to another.

We have
\begin{multline}\label{lem-WP2-coucou}
e^{2 L y_1(t_1, \xi_1)} |v - w|^2 \big(s, x_j (s, t_1, \xi_1) \big)  \\[6pt]
=  e^{ 2 L \big(y_1(t_1, \xi_1) - y_1 (s, x_j (s, t_1, \xi_1)) \big)}  e^{ 2L  y_1 (s, x_j (s, t_1, \xi_1) )} |v - w|^2 \big(s, x_j (s, t_1, \xi_1) \big) ,
\end{multline}
and, for $s$ between $t_1$ and $t$,
\begin{equation}\label{lem-WP2-c1}
y_1(t_1, \xi_1) - y_1 (s, x_j (s, t_1, \xi_1)) \le - c |\xi_1 - x_j (s, t_1, \xi_1)|  \mbox{ by a2) and the definition of $G_1$}.
\end{equation}
Making a change of variables $x = x_j (s, t_1, \xi_1)$, we derive from \eqref{lem-WP2-cc} that, for $1 \le j \le \ell$ or $k+1 \le j \le k + m$,
\begin{multline}\label{lemWP2-p1}
\int_{\Omega_{\ell, \xi_1}} e^{2 L y_1(t_1, \xi_1)} | \cF(v)_j (t_1, \xi_1) -  \cF(w)_j(t_1, \xi_1)|^2 \, d t_1 \\[6pt] \le  C \int_{\Omega_\ell} e^{- c L |\xi_1 - x|} e^{2 L y_1 (s, x)}  |v - w|^2 (s, x) \, ds \, d x \le \frac{C}{L} \| v- w \|_{\Omega_\ell}^2.
\end{multline}

We next deal with $\ell + 1 \le j \le k$. Set
$$
\Omega_{\ell, \xi_1, 1} = \Big\{t_1 \in [0, T_{opt}];  \mbox{ \eqref{def-WP-A2-p2} holds} \Big\} \quad \mbox{ and }
\quad \Omega_{\ell, \xi_1, 2} = \Big\{t_1 \in [0, T_{opt}];  \mbox{ \eqref{def-WP-A2-p3} holds} \Big\}.
$$
We have, by \eqref{def-WP-A2-p2},  for $t_1 \in \Omega_{\ell, \xi_1, 1}$,
\begin{multline}\label{lemWP2-p2}
\cF(v)_j(t_1, \xi_1) - \cF(w)_j(t_1, x_1) \\[6pt]
= \int_{\hat t}^{t_1} \Big( F \big(s, x_j (s, t_1, \xi_1) \big) (v-w)(s, x_j (s, t_1, \xi_1))  \Big)_j\, ds + \big(Q_{\ell+1} (v-w)_{< \ell + 1, \ge \ell + m + 1} \big)_{j - \ell} (\hat t, 0)
\end{multline}
where $\hat t  = \hat t (t_1, \xi_1)$ is such that $ x_j(\hat t, t_1, \xi_1) = 0$.

We next estimate
\begin{equation*}
\int_{\Omega_{\ell, \xi_1, 1}} \sign(\hat t - t_1) \int_{t_1}^{\hat t} e^{ 2 L y_1(t_1, \xi_1)}   |v-w|^2 \big(s, x_j (s, t_1, \xi_1) \big) \, ds  \, d t_1.
\end{equation*}
We have, for $s$ between $t_1$ and $\hat t$,
\begin{equation}\label{lem-WP2-c2}
y_1(t_1, \xi_1) - y_1 (s, x_j (s, t_1, \xi_1)) \le - c |\xi_1 - x_j (s, t_1, \xi_1)|  \mbox{ by a3) and the definition of $G_2$}.
\end{equation}
Making a change of variables $x = x_j (s, t_1, \xi_1)$, we derive from \eqref{lem-WP2-coucou} that
\begin{multline*}
\int_{\Omega_{\ell, \xi_1, 1}} \sign(\hat t - t_1)  \int_{t_1}^{\hat t} e^{2 L y_1(t_1, \xi_1)}   |v-w|^2\big(s, x_j (s, t_1, \xi_1) \big) \, ds  \, d t_1 \\[6pt] \le C \int_{\Omega_{\ell}} e^{ - c L |\xi_1 - x|} e^{- 2L  y_1 (s, x)}  |v - w|^2 (s, x) \, ds \, d x.
\end{multline*}
This implies
\begin{equation}\label{lemWP2-p3}
\int_{\Omega_{\ell, \xi_1}} \sign(\hat t - t_1)  \int_{t_1}^{\hat t} e^{2 L y_1(t_1, \xi_1)}   |v-w|^2(s, x_j (s, t_1, \xi_1)) \, ds  \, d t_1 \le \frac{C}{L} \| v- w \|_{\Omega_\ell}^2.
\end{equation}

By \eqref{lemWP2-p1}, we also have
\begin{equation}\label{lemWP2-p4}
\int_{\Omega_{\ell, 0}} e^{2L y_1(\hat t, 0)}|Q_{\ell+1} (v-w)_{< \ell + 1, \ge \ell + m  + 1}(\hat t, 0)|^2 \, d \hat t \le \frac{C}{L} \| v- w \|_{\Omega_\ell}^2.
\end{equation}
Using \eqref{lem-WP2-c2},  and making a change of variable $\hat t = \hat t(t_1, \xi_1)$, we derive that
\begin{multline}\label{lemWP2-p5}
\int_{\Omega_{\ell, \xi_1}} e^{L y_1 (t_1, \xi_1)} |Q_{\ell+1} (v-w)_{< \ell + 1, \ge \ell + m + 1}(\hat t (t_1, \xi_1), 0)|^2 \, d t_1 \\[6pt]
 \le C \int_{\Omega_{\ell, 0}} e^{2L y_1(\hat t, 0)}|Q_{\ell+1} (v-w)_{< \ell + 1, \ge \ell + m + 1}(\hat t, 0)|^2 \, d \hat t.
\end{multline}

Combining \eqref{def-WP-A2-p2}, \eqref{lemWP2-p3}, \eqref{lemWP2-p4}, and \eqref{lemWP2-p5} yields, for $\ell + 1 \le j \le k$,
\begin{equation}\label{lemWP2-p6-1}
\int_{\Omega_{\ell, \xi_1, 1}} e^{2 L y_1(t_1, \xi_1)} |\cF(v)_j(t_1, \xi_1) - \cF(w)_j(t_1, x_1)|^2 \, d t_1 \le \frac{C}{L} \| v- w \|_{\Omega_\ell}^2.
\end{equation}
Using similar arguments, we also obtain, for $\ell + 1 \le j \le k$,
\begin{equation}\label{lemWP2-p6-2}
\int_{\Omega_{\ell, \xi_1, 2}} e^{2 L y_1(t_1, \xi_1)} |\cF(v)_j(t_1, \xi_1) - \cF(w)_j(t_1, x_1)|^2 \, d t_1 \le \frac{C}{L} \| v- w \|_{\Omega_\ell}^2.
\end{equation}
We derive from \eqref{lemWP2-p6-1} and \eqref{lemWP2-p6-2} that
\begin{equation}\label{lemWP2-p6}
\int_{\Omega_{\ell, \xi_1}} e^{2 L y_1(t_1, \xi_1)} |\cF(v)_j(t_1, \xi_1) - \cF(w)_j(t_1, x_1)|^2 \, d t_1 \le \frac{C}{L} \| v- w \|_{\Omega_\ell}^2.
\end{equation}

From \eqref{lemWP2-p1} and \eqref{lemWP2-p6}, we obtain
\begin{equation}\label{lemWP2-p7}
\int_{\Omega_{\xi_1}} e^{2 L y_1(t_1, \xi_1)} |\cF(v)(t_1, \xi_1) - \cF(w)(t_1, x_1)|^2 \, d t_1 \le \frac{C}{L} \| v- w \|_{\Omega_\ell}^2.
\end{equation}

For $t_1 \in (0, T_{opt})$, by the same approach used to derive \eqref{lemWP2-p7},   we also have
\begin{equation}\label{lemWP2-p8}
\int_{\Omega_{t_1}} e^{2 L y_1(t_1, \xi_1)} |\cF(v)(t_1, \xi_1) - \cF(w)(t_1, x_1)|^2 \, d \xi_1 \le \frac{C}{L} \| v- w \|_{\Omega_\ell}^2.
\end{equation}

The conclusion in the case where $\Sigma$ is constant now follows from \eqref{lemWP2-p7} and \eqref{lemWP2-p8}.

\medskip
We next make necessary modifications to derive the conclusion in the general case. The idea is to find a replacement for $y_1(t, x)$ which is {\it increasing} when one follows the characteristic flows directed to the boundary for which the boundary conditions are imposed. To this end, for $1 \le j \le k + m$, let  $\vec{v}_j = \vec{v}_j(t, x)$ be the unit vector tangent to the characteristic curve of $x_j$ at the point $(t, x)$ directed to the boundary where the boundary condition for  $v_j$ is given. The vector $\vec{v}_j (t, x)$ is parallel to $(1, \Sigma_{jj}(x))\tr$ in the $xt$-plane and so that one can choose it independent of $t$ and in fact we will do. We will  denote it by $\vec{v}_j(x)$ from now on.
Set
$$
G_1(x) = \Big\{\vec{v}_j (x); 1 \le j \le \ell, k +1 \le j \le k+m \Big\} \quad \mbox{ and } \quad G_2 (x) = \Big\{\vec{v}_j (x); \ell + 1 \le j \le  k \Big\}.
$$
Let $\varphi(x)$ be such that
$$
\vec{v}_{\ell} (x) \mbox{ is parallel to and has the same direction with } (\varphi(x), 1)\tr.
$$
Set, in the $xt$-plane,
$$
\vec{u}_1 (x) = (0, -1)\tr,
$$
and
$$
\vec{u}_2 (x) =  (\varphi(x) - \eps, 1)\tr,
$$
where $\eps$ is a  constant which is positive and  sufficiently small, the smallness of $\eps$ is independent of $x$,  such that, $\varphi(x) > 2 \eps$, and
\begin{enumerate}
\item[a1)] $G_1 (x) \cup G_2(x) \cup \{\vec{u}_1 (x) \}$ lies on one side of the line containing $\vec{u}_2 (x)$
\item[a2)] $G_1 (x)$ is a subset of  the open solid cone centered at the origin and  formed by  $\vec{u}_1(x)$ and $\vec{u}_2 (x)$,  i.e.,  in the set $\big\{s_1 \vec{u}_1 (x) + s_2 \vec{u}_2 (x); s_1, s_2 > 0 \big\}$.

\item[a3)] $G_2(x)$ is a subset of  the open solid cone centered at the origin and  formed by  $\vec{u}_1 (x)$ and $-\vec{u}_2(x)$,  i.e.,  in the set $\big\{s_1 \vec{u}_1 (x) - s_2 \vec{u}_2 (x); s_1, s_2 > 0 \big\}$.

\end{enumerate}

Fix such a positive constant $\eps$. For a point $(x_0, t_0) \in \Omega_\ell$, let $(x(s), t(s))$ for  $ s \in [\alpha, \beta] \subset \mR$ be a (piecewise) $C^1$ regular curve \footnote{Regularity  means that $(x'(s), t'(s)) \neq (0, 0)$ for $s \in [\alpha, \beta]$ such that $(x'(s), t'(s))$ is well-defined.} in $\bar \Omega_\ell$ (in the $xt$-plane) starting from  $(0, 0)$ and arriving at $(x_0, t_0)$. We first claim that
\begin{multline}\label{claim}
\int_\alpha^\beta y_1\big(x'(s), t'(s),  x(s), t(s) \big) |(x'(s), t'(s))| \, ds \mbox{ depends  on $(t_0, x_0)$} \\
\mbox{but is independent of the curve and the parametrization.}
\end{multline}
Here
$y_1\big(t'(s), x'(s), t(s), x(s) \big)$ is the first coordinate of the vector $(t'(s), x'(s))/ |(t'(s), x'(s))|$ in the bases $\vec{u}_1 (t(s), x(s))$  and $\vec{u}_2 (t(s), x(s))$.

We now establish the claim. For notational ease, we assume that $|(t'(s), x'(s))| = 1$. We first compute $y_1\big(t'(s), x'(s), t(s), x(s) \big)$. Let $a$ and $b$ in $\mR$ be such that
$$
(x'(s), t'(s)) = a  (0, -1)  + b (\varphi(x(s)) - \eps, 1).
$$
We have
$$
a = - t'(s) +  \frac{x'(s)}{\varphi (x(s)) - \eps} \quad \mbox{ and } \quad  b = \frac{x'(s)}{\varphi (x(s)) - \eps}.
$$
Thus
$$
y_1\big(t'(s), x'(s), t(s), x(s) \big)= - t'(s) +  \frac{x'(s)}{\varphi (x(s)) - \eps}.
$$
It follows that
\begin{equation}\label{def-Phi}
\int_\alpha^\beta y_1\big(x'(s), t'(s),  x(s), t(s) \big) \, ds  = - t_0 + \Phi(x_0),
\end{equation}
where
$$
\Phi(\xi) = \int_0^\xi \frac{1}{\varphi(s) + \eps} \, ds \mbox{ for } \xi \in [0, 1].
$$
The claim is proved.

Define
\begin{equation*}
\begin{array}{cccc}
Y_1 &: \Omega_\ell & \to &  \mR, \\[6pt]
& (t, x) & \mapsto &   - t + \Phi(x).
\end{array}
\end{equation*}
The proof in the general case follows as in the constant case with $T_\ell$ now defined by
\begin{equation}\label{lem-WP2-w2}
T_\ell (v) (t, x) =  e^{L Y_1(t, x)} v(t, x).
\end{equation}
One just notes that \eqref{lem-WP2-c1} and \eqref{lem-WP2-c2} hold with $y_1$ replaced by $Y_1$. Indeed, one has
\begin{multline*}
 Y_1 \big(s, x_j (s, t_1, \xi_1) \big) - Y_1(t_1, \xi_1) = \int_{t_1}^s y_1 \big( \partial_\theta x_j (\theta, t_1, \xi_1), x_j (\theta, t_1, \xi_1) \big) |\partial_\theta x_j (\theta, t_1, \xi_1)| \, d \theta \\[6pt]
 \ge C \sign(s - t_1) \int_{t_1}^s y_1 \big( \vec{v}_j (x_j (\theta, t_1, \xi_1)), x_j (\theta, t_1, \xi_1) \big) \, d \theta \ge C |t_1 - s| \ge C |x_j (s, t_1, \xi_1) - \xi_1|.
\end{multline*}
The details are omitted.
\end{proof}

We are ready to give

\begin{proof}[Proof of \Cref{thm-WP}] We first prove the uniqueness. Assume that $f = 0$, $g = 0$, and $\gamma = 0$. Then the restriction of $w$ into $\Omega_k$ is 0 by \Cref{lem-WP1}. It follows that the restriction of $w$ into $\Omega_{k-1} = 0$ by \Cref{lem-WP2}, \dots,   the restriction of $w$ into $\Omega_{k- m + 1} = 0$ by \Cref{lem-WP2}. Therefore, $w= 0$ in $\Omega$.

To establish the existence, we proceed as follows. Let $w^{(k)}$ be the unique broad solution in $\Omega_k$ corresponding to $(f, g)$, let $w^{(k-1)}$ be the unique broad solution in $\Omega_{k-1}$  where the data on $\Gamma_k$ come from  $w^{(k)}$, \dots, let $w^{(k-m+1)}$ be the unique broad solution in $\Omega_{k-m+1}$  where the data on $\Gamma_{k-m+2}$ come from  $w^{(k-m+2)}$ \footnote{The data coming from $w^{(k)}$ on $\Gamma_k$, \dots, $w^{(k-m+2)}$ on $\Gamma_{k-m+2}$ are given by the RHS of \eqref{def-WP-A1-p1}-\eqref{def-WP-A1-p3} in \Cref{def-WP1} for $(t_1, \xi_1) \in \Gamma_k$, and \eqref{def-WP-A2-p1}-\eqref{def-WP-A2-p3} in \Cref{def-WP2} for $(t_1, \xi_1) \in \Gamma_{\ell}$ with $\ell = k-1$, \dots, $k-m+1$, respectively.}. The corresponding solution is obtained by gluing these solutions together.  The proof is complete.
\end{proof}

\begin{remark}\label{rem-Weight} \rm  The introduction of  appropriate weighted norms plays a crucial role in the proof of
the well-posedness of broad solutions considered so far in this section, in particular in the proof of \Cref{lem-WP2}.
The introduction of  weighted norms in order to be able to apply the fixed point argument used in establishing the well-posedness of hyperbolic system is not new. The standard one is $e^{-Lt}$ where $L$ is a large positive number, see e.g. \cite[(1.18), p. 78]{LY85} or \cite[(3.36), p. 50]{Bressan00}, while the weight $e^{-Lx}$ is used in \cite{2002-Xu-Sallet-COCV,CBA08} to prove exponential stability; see also \cite[$V$ defined in Section 3.2]{1999-Coron-SICON} for the Euler equations of incompressible fluids. In \cite{CoronNg19}, we used the weight $e^{-L_1 x - L_2 t}$ where $L_1$ and $L_2$ are two large positive numbers with $L_2$ being much larger than $L_1$. The introduction of $e^{-L_1 x}$ in the weight is to handle the non-local term from the boundary condition imposed on the right (at $x =1$) considered there. In these settings, $t$-direction has a  privileged role. In the settings considered in  this section, the domain is not a rectangle with respect to $t$ and $x$, and the boundary conditions are quite complicated. Therefore,   the time direction and the space direction play almost the same role here. In the setting of \Cref{lem-WP1}, the privileged direction is $x$-direction so the weighted norm is chosen of the form $e^{Lx}$. In \Cref{lem-WP2},
the {\it new} weighted norm introduced in \eqref{lem-WP2-w} with $T_\ell$ given by \eqref{lem-WP2-w1} or \eqref{lem-WP2-w2} adapts the geometry and the boundary conditions,  imposed in a nontrivial way.  It is interesting to note that $Y_1$ is a {\it non-linear} function of $t$ and $x$.  The analysis here is inspired by \cite{CoronNg19} (see also \cite{CoronNg20}).
\end{remark}

As a consequence of \Cref{thm-WP}, we can prove

\begin{proposition}\label{pro-WP} Let  $C \in \big[L^\infty \big( I \times (0, 1)\big) \big]^{n \times n}$. Define, for $\tau \in I_1$,
\begin{equation}
\begin{array}{cccc}
\cT(\tau):  & [L^2(0, T_{opt})]^n \times  [L^2(0, 1)]^m  & \to  & \cY \\[6pt]
 & (f, g) & \mapsto & w,
\end{array}
\end{equation}
where $w$ is the solution of  \eqref{w1}-\eqref{w6}. Then $\cT(\tau)$  is uniformly bounded in $I_1$. Assume in addition that $C \in \cH(I, [L^\infty(0, 1)]^{n \times n})$. Then $\cT(\tau, \cdot, \cdot)$ is analytic in $I_1$.
\end{proposition}

\begin{proof} By \Cref{thm-WP}, for each $(f, g) \in [L^2(0, 1)]^n \times  [L^2(0, 1)]^m$, there exists a unique broad solution $w \in \cY$ of \eqref{w1}-\eqref{w6}. Hence $\cT(\tau)$ is well-defined.
The uniform boundedness of $\cT$ is also a direct consequence of \Cref{thm-WP}, in particular of \eqref{thm-WP-e1}.

We next deal with the analyticity of $\cT$ and thus assume that $C \in \cH(I, [L^\infty(0, 1)]^{n \times n})$.  Fix $\tau_0$ in a sufficiently small neighborhood of $I_1$ (in the complex plane). We will prove that $\cT$ is differentiable at $\tau_0$ in the complex sense. For notational ease, we will assume that $\tau_0 = 0$.

Fix $(f, g) \in [L^2(0, T_{opt})]^n \times [L^2(0, 1)]^n$. Set $w^{(\tau)} = \cT(\tau) (f, g)$ in $\Omega$ for $\tau$ in a small neighborhood (in the complex plane) of $0$ and let $v \in \cY$ be the unique broad solution of the system
\begin{equation}\label{v1}
v_t (t, x)  = \Sigma (x) \partial_x v (t, x)+ \bC(t, x) v(t, x) +   \bC_\tau (t, x)  w^{(0)} (t, x) \mbox{ for } (t, x) \in \Omega,
\end{equation}
\begin{equation}\label{v2}
v (\cdot, 1)= 0  \mbox{ in } (0, T_{opt}),
\end{equation}
\begin{equation}\label{v3}
v_+(0, \cdot) = 0 \mbox{ in } (0, 1),
\end{equation}
\begin{equation}\label{v4}
v_{-, \ge k}(t, 0) = Q_k v_{<k, \ge k+m} (t, 0) \mbox{ for } t \in (T_{opt} - \tau_{k}, T_{opt} - \tau_{k-1}),
\end{equation}
\begin{equation}\label{v5}
v_{-, \ge k-1}(t, 0) = Q_{k-1} v_{<k-1, \ge k+m-1} (t, 0) \mbox{ for } t \in  (T_{opt} - \tau_{k-1}, T_{opt} - \tau_{k-2}),
\end{equation}
\dots
\begin{equation}\label{v6}
v_{-, \ge k-m+ 2}(t, 0) = Q_{k-m + 2} v_{<k-m+2, \ge k+2} (t, 0) \mbox{ for } t \in (T_{opt} - \tau_{k-m+2},  T_{opt} - \tau_{k - m + 1}).
\end{equation}
Here $\bC_\tau(\tau, x)$ denotes the derivative of $\bC(\tau, x)$ with respect to $\tau$ in the complex sense. The existence and uniqueness of $v$ follows from \Cref{thm-WP}.

We claim that
\begin{equation}\label{pro-WP-claim}
\mbox{ the derivative of $\cT$ at $0$ is given by $\cT_1$ where } \cT_1(f, g) = v \mbox{ in } \Omega
\end{equation}
(the derivative of $\cT$ is considered in the complex sense).  To this end, for $\tau$ in a small neighborhood (in the complex plane) of $0$ but not $0$, we consider $dw \in \cY$ defined by
$$
dw: = \frac{1}{\tau} \Big( w^{(\tau)} - w^{(0)} - \tau v \Big) \mbox{ in } \Omega.
$$
Then $dw \in \cY$ is a broad solution of the system
\begin{multline}\label{pro-WP-p1}
\partial_t dw (t, x) = \Sigma(x) \partial_x dw (t, x) + \bC(t, x) dw (t, x) \\[6pt]
+ \frac{1}{\tau} \Big( \bC(t + \tau, x) - \bC(t, x)\Big) w^{(\tau)}(t, x) - \bC_\tau (t, x) w^{(0)} (t, x) \mbox{ in } \Omega,
\end{multline}
and \eqref{v2}-\eqref{v6} with $v$ replaced by $dw$. We derive from \Cref{thm-WP} that
\begin{equation}\label{pro-WP-coucou}
\| dw\|_{\cY} \le C \Big(  \| w^{(\tau)}\|_{L^2(\Omega)} +  \| w^{(0)}\|_{L^2(\Omega)} \Big) \le C \Big( \| f\|_{L^2(0,T_{opt})} + \| g\|_{L^2(0, 1)} \Big).
\end{equation}
Using the definition of $dw$, we can write the last two terms in \eqref{pro-WP-p1} under the form
\begin{multline}\label{pro-WP-p2}
\frac{1}{\tau} \Big( \bC(t + \tau, x) - \bC(t, x)\Big) \Big(w^{(0)} + \tau dw  + \tau v \Big)  - \bC_\tau (t, x) w^{(0)} (t, x)  \\[6pt]
=  \frac{1}{\tau} \Big( \bC(t + \tau, x) - \bC(t, x) - \tau \bC_\tau(t, x) \Big) w^{(0)}(t,x) + \frac{1}{\tau} \Big( \bC(t + \tau, x) - \bC(t, x)\Big) \Big(\tau dw  + \tau v \Big).
\end{multline}
Note that the $L^2(\Omega)$-norm of the RHS of \eqref{pro-WP-p2} is bounded by
$$
C |\tau| \left( \| w^{(0)}\|_{L^2(\Omega)} + \| dw\|_{L^2(\Omega)} + \| v\|_{L^2(\Omega)} \right).
$$
Applying \Cref{thm-WP} again, we derive from \eqref{pro-WP-coucou} that
\begin{equation}\label{pro-WP-p3}
\|dw \|_{\cY} \le C |\tau| \left( \| w^{(0)}\|_{L^2(\Omega)} + \| v\|_{L^2(\Omega)} +  \| f\|_{L^2(0,T_{opt})} + \| g\|_{L^2(0, 1)} \right).
\end{equation}
By noting that
$$
 \| w^{(0)}\|_{L^2(\Omega)} + \| v\|_{L^2(\Omega)} \le C \Big( \|f \|_{L^2(0, T_{opt})} + \| g\|_{L^2(0, 1)} \Big),
$$
claim \eqref{pro-WP-claim} follows from \eqref{pro-WP-p3}. The proof is complete.
\end{proof}

\begin{remark} \label{rem-SC} \rm Let  $C \in \big[L^\infty \big( I \times (0, 1)\big) \big]^{n \times n}$. One can prove that $\cT(\tau)$ is strongly continuous, i.e., $\cT(\tau)(f, g) \to \cT(\tau_0)(f, g)$ in $\cY$ as $\tau \to \tau_0$ in $I_1$ for all $(f, g) \in [L^2(0, T_{opt})]^n \times [L^2(0, 1)]^m$. Indeed, let assume that $\tau_0 = 0$ for notational ease. Set $w^{(\tau)} = \cT(\tau) (f, g)$ in $\Omega$ for $\tau \in I_1$ and for $(f, g) \in [L^2(0, T_{opt})]^n \times [L^2(0, 1)]^m$. Denote $\delta w = w^{(\tau)} - w^{(0)}$ in $\Omega$. We have, in $\Omega$
\begin{equation*}
\partial_t \delta w (t, x) = \Sigma(x) \partial_x \delta w (t, x) + \bC(t + \tau, x) \delta w (t, x)
+ \Big( \bC(t + \tau, x) - \bC(t, x)\Big) w^{(0)}(t, x), 
\end{equation*}
and $\delta w$ satisfies the same boundary conditions as $dw$. Applying \Cref{thm-WP}, one has
\begin{equation*}
\|\delta w \|_{\cY} \le C \| g\|_{L^2(\Omega)},
\end{equation*}
where $g(t, x) = \Big( \bC(t + \tau, x) - \bC(t, x)\Big) w^{(0)}(t, x)$. Since $\| g\|_{L^2(\Omega)} \to 0$ as $\tau \to 0$, the conclusion follows.
\end{remark}

We next discuss the broad solutions used in the definition of  $\hat \cT(\tau)$ and their well-posedness. Let $F \in \big[L^\infty\big( (0, T_{opt}) \times (0, 1)  \big) \big]^{n\times n}$, $(f, g) \in [L^2(0, T_{opt})]^n \times [L^2(0, 1)]^m$, and let $q \in [L^2(0, 1)]^{k-m}$ \footnote{$q$ is irrelevant when $k=m$.}.  Consider the system
\begin{equation}\label{whA1}
\partial_t \hw (t, x)  = \Sigma (x) \partial_x \hw (t, x)+ F(t, x) \hw (t, x) + \gamma(t, x) \mbox{ for } (t, x) \in \Omega,
\end{equation}
\begin{equation}\label{whA2}
\hw (\cdot, 1)= f  \mbox{ in } (0, T_{opt}),
\end{equation}
\begin{equation}\label{whA3}
\hw_+(0, \cdot) = g \mbox{ in } (0, 1),
\end{equation}
\begin{equation}\label{whA4}
\hw_j(T_{opt}, \cdot) = q_j \mbox{ in } (0, 1), \mbox{ for } 1 \le j \le k-m,
\end{equation}
\begin{equation}\label{whA5}
\hw_{-, \ge k}(t, 0) = Q_k \hw_{<k, \ge k+m} (t, 0) \mbox{ for } t \in (T_{opt} - \tau_{k}, T_{opt} - \tau_{k-1}),
\end{equation}
\begin{equation}\label{whA6}
\hw_{-, \ge k-1}(t, 0) = Q_{k-1} \hw_{<k-1, \ge k+m-1} (t, 0) \mbox{ for } t \in  (T_{opt} - \tau_{k-1}, T_{opt} - \tau_{k-2}),
\end{equation}
\dots
\begin{equation}\label{whA7}
\hw_{-, \ge k-m+ 2}(t, 0) = Q_{k-m + 2} \hw_{<k-m+2, \ge k+2} (t, 0) \mbox{ for } t \in (T_{opt} - \tau_{k-m+2},  T_{opt} - \tau_{k - m + 1}),
\end{equation}
\begin{equation}\label{whA8}
\hw_{-, \ge k-m+ 1}(t, 0) = Q_{k-m + 1} \hw_{<k-m+1, \ge k+1} (t, 0) \mbox{ for } t \in (T_{opt} - \tau_{k-m+1},  T_{opt} ).
\end{equation}

We have the following result, which implies the well-posedness of $\hat \cT(\tau)$.

\begin{theorem}\label{thm-hWP} Let $F \in \big[L^\infty\big( (0, T_{opt}) \times (0, 1)  \big) \big]^{n\times n}$, $(f, g) \in [L^2(0, T_{opt})]^n \times [L^2(0, 1)]^m$,  $q \in [L^2(0, 1)]^{k-m}$ \footnote{$q$ is irrelevant when $k=m$.}, and $\gamma \in [L^2(\Omega)]^n$. There exists a unique broad solution
$$\hw \in \hat  \cY: =
\big[L^2\big( (0, T_{opt}) \times (0, 1) \big)\big]^n \cap C \big( [0, T_{opt}]; [L^2 (0, 1) ]^n \big) \cap C \big( [0, 1]; [L^2 (0, T_{opt}) ]^n \big)
$$
of \eqref{whA1}-\eqref{whA8}. Moreover,
$$
\| \hw\|_{\hat \cY} \le C \Big( \| f\|_{L^2(0, T_{opt})} + \| g\|_{L^2(0, 1)} + \| q\|_{L^2(0,1)} + \| \gamma\|_{L^2\big( (0, T_{opt}) \times (0, 1)\big)} \Big),
$$
for some positive constant $C$ depending only on an upper bound of $\| F \|_{L^\infty(\Omega_\ell)}$ and $\Sigma$.
\end{theorem}

Here we denote
$$
\| \hw\|_{\hat \cY} =  \max\left\{ \sup_{x \in [0, 1]} \| \hw \|_{L^2(0, T_{opt})}, \sup_{t \in [0, T_{opt}]} \| \hw  \|_{L^2(0, 1)}; 1 \le i \le n \right\}.
$$

\begin{remark} \rm The analysis of \Cref{thm-hWP} can be extended  to cover the case where source terms in $L^2$ are added in \eqref{whA5}-\eqref{whA8}.
\end{remark}

The definition of broad solutions $\hw \in \hat \cY$ of \eqref{whA1}-\eqref{whA8} is similar to the one given in \Cref{def-WP} and left to the reader.  The proof of \eqref{thm-hWP} is similar to the one of \Cref{thm-WP}. Nevertheless, in addition to \Cref{lem-WP1,lem-WP2}, we also use the following.

\begin{lemma} \label{lem-WP3} Set $\Omega_{k-m} = [0, T_{opt}] \times (0, 1) \setminus \Omega$. Let $F \in [L^\infty(\Omega_{k-m})]^{n \times n}$ $\gamma \in [L^2(\Omega_{k-m})]^n$, $h_j \in L^2(\Gamma_{k-m+1})$ for $1 \le j \le  k+m$ and $j \neq k-m+1$, and let $q_j \in L^2(\Gamma_{k-m})$ for $1 \le j \le k-m$ where $\Gamma_{k-m} = \{T_{opt}\} \times (0, 1)$.
There exists a unique broad solution
 $w \in \cY_{k-m}: = \big[L^2(\Omega_{k-m}) \big]^n \cap C \big( [0, T_{opt}]; \big[L^2 (\Omega_{k-m, t}) \big]^n \big) \cap C \big( [0, 1]; \big[L^2 (\Omega_{k-m, x}) \big]^n \big) $ of the system
\begin{equation}
\partial_t w (t, x)  = \Sigma (x) \partial_x w (t, x)+ F(t, x) w (t, x) + \gamma(t, x) \mbox{ for } (t, x) \in \Omega_{k-m},
\end{equation}
\begin{equation}
w_j  =  h_j \mbox{ on } \Gamma_{k-m+1}, \mbox{ for $1 \le j  \le k+m$ and $j \neq k-m+1$},
\end{equation}
\begin{equation}
w_j  =  q_j \mbox{ on } \Gamma_{k-m}, \mbox{ for $1 \le j  \le  k-m$},
\end{equation}
\begin{equation}\label{lem-WP3-p}
w_{-, \ge k-m + 1}(0, \cdot) = Q_{k-m+1} w_{< k - m + 1, \ge k + 1} \mbox{ for } t \in (T_{opt} - \tau_{k-m+1}, T_{opt}).
\end{equation}
Moreover,
$$
\| w\|_{\cY_\ell} \le  C \left( \sum_{1 \le j \le  k+m; j \neq k-m + 1} \|  h_j\|_{L^2(\Gamma_{k-m+1})} +  \sum_{1 \le j \le k -m} \|  q_j\|_{L^2(\Gamma_{k-m})} + \|\gamma \|_{L^2(\Omega_{k-m})} \right),
$$
for some positive constant $C$ depending only on an upper bound of $\| F \|_{L^\infty(\Omega_{k-m})}$ and $\Sigma$.
\end{lemma}

\begin{remark} \rm The analysis  of \Cref{lem-WP3} can be extended to cover the case where source terms in $L^2$ are added in \eqref{lem-WP3-p}.
\end{remark}

\begin{proof} The proof of \Cref{lem-WP3} is similar to the one of \Cref{lem-WP2}.
We just mention here how to define $G_1$, $G_2$ and determine $\vec{u_1}$ and $\vec{u}_2$ in the general case ($\Sigma$ is not required to be constant).
For $1 \le j \le k + m$, let  $\vec{v}_j = \vec{v}_j(t, x)$ be the unit vector tangent to the characteristic curve of $x_j$ at the point $(t, x)$ directed to the boundary where the boundary condition for  $v_j$ is given. The vector $\vec{v}_j (t, x)$ is parallel to $(1, \Sigma_{jj}(x))\tr$ in the $xt$-plane and so that we can choose it  independent of $t$ and in fact we will do. We denote it by $\vec{v}_j(x)$ from now on.
Set
$$
G_1(x) = \Big\{\vec{v}_j (x); 1 \le j \le k-m, k +1 \le j \le k+m \Big\} \quad \mbox{ and } \quad G_2 (x) = \Big\{\vec{v}_j (x); k-m + 1 \le j \le  k \Big\}.
$$
Let $\varphi(x)$ be such that
$$
\vec{v}_{1} (x) \mbox{ is parallel to and has the same direction with } (\varphi(x), 1)\tr.
$$
Set, in the $xt$-plane,
$$
\vec{u}_1 (x) = (0, -1)\tr,
$$
and
$$
\vec{u}_2 (x) =  (\varphi(x) - \eps, 1)\tr \mbox{ if } k > m,  \mbox{ otherwise } \vec{u_2} = (1, 0)\tr,
$$
where $\eps$ is positive and  sufficiently small, the smallness of $\eps$ is independent of $x$,  such that, $\varphi(x) > 2 \eps$ (the choice of $\eps$ is irrelevant when $k=m$), and
\begin{enumerate}
\item[a1)] $G_1 (x) \cup G_2(x) \cup \{\vec{u}_1 (x) \}$ lies on one side of the line containing $\vec{u}_2 (x)$.

\item[a2)] $G_1 (x)$ is a subset of  the open solid cone centered at the origin and  formed by  $\vec{u}_1(x)$ and $\vec{u}_2 (x)$.

\item[a3)] $G_2(x)$ is a subset of  the open solid cone centered at the origin and  formed by  $\vec{u}_1 (x)$ and $-\vec{u}_2(x)$.
\end{enumerate}
The rest of the proof is then almost unchanged and left to the reader.
\end{proof}

\medskip
\noindent \textbf{Acknowledgments.} The authors were partially supported by  ANR Finite4SoS ANR-15-CE23-0007. H.-M. Nguyen thanks Fondation des Sciences Math\'ematiques de Paris (FSMP) for the Chaire d'excellence which allowed him to visit  Laboratoire Jacques Louis Lions  and Mines ParisTech. Part of this work has been done during this visit.

\providecommand{\bysame}{\leavevmode\hbox to3em{\hrulefill}\thinspace}
\providecommand{\MR}{\relax\ifhmode\unskip\space\fi MR }
\providecommand{\MRhref}[2]{%
  \href{http://www.ams.org/mathscinet-getitem?mr=#1}{#2}
}
\providecommand{\href}[2]{#2}


\begin{thebibliography}{10}

\bibitem{Ahlfors}
Lars~V. Ahlfors, \emph{Complex analysis}, third ed., McGraw-Hill Book Co., New
  York, 1978, An introduction to the theory of analytic functions of one
  complex variable, International Series in Pure and Applied Mathematics.
  \MR{510197}

\bibitem{SFA08}
Saurabh Amin, Falk~M. Hante, and Alexandre~M. Bayen, \emph{On stability of
  switched linear hyperbolic conservation laws with reflecting boundaries},
  Hybrid systems: computation and control, Lecture Notes in Comput. Sci., vol.
  4981, Springer, Berlin, 2008, pp.~602--605. \MR{2728909}

\bibitem{ALM16}
Nalini Anantharaman, Matthieu L\'{e}autaud, and Fabricio Maci\`a, \emph{Wigner
  measures and observability for the {S}chr\"{o}dinger equation on the disk},
  Invent. Math. \textbf{206} (2016), no.~2, 485--599. \MR{3570298}

\bibitem{AM16}
Jean Auriol and Florent Di~Meglio, \emph{Minimum time control of
  heterodirectional linear coupled hyperbolic {PDE}s}, Automatica J. IFAC
  \textbf{71} (2016), 300--307. \MR{3521981}

\bibitem{BLR92}
Claude Bardos, Gilles Lebeau, and Jeffrey Rauch, \emph{Sharp sufficient
  conditions for the observation, control, and stabilization of waves from the
  boundary}, SIAM J. Control Optim. \textbf{30} (1992), no.~5, 1024--1065.
  \MR{1178650}

\bibitem{BC16}
Georges Bastin and Jean-Michel Coron, \emph{Stability and boundary
  stabilization of 1-{D} hyperbolic systems}, Progress in Nonlinear
  Differential Equations and their Applications, vol.~88,
  Birkh\"auser/Springer, [Cham], 2016, Subseries in Control. \MR{3561145}

\bibitem{Bressan00}
Alberto Bressan, \emph{Hyperbolic systems of conservation laws}, Oxford Lecture
  Series in Mathematics and its Applications, vol.~20, Oxford University Press,
  Oxford, 2000, The one-dimensional Cauchy problem. \MR{1816648}

\bibitem{CoronC13}
Eduardo Cerpa and Jean-Michel Coron, \emph{Rapid stabilization for a
  {K}orteweg-de {V}ries equation from the left {D}irichlet boundary condition},
  IEEE Trans. Automat. Control \textbf{58} (2013), no.~7, 1688--1695.
  \MR{3072853}

\bibitem{1999-Coron-SICON}
Jean-Michel Coron, \emph{On the null asymptotic stabilization of the
  two-dimensional incompressible {E}uler equations in a simply connected
  domain}, SIAM J. Control Optim. \textbf{37} (1999), no.~6, 1874--1896.
  \MR{1720143}

\bibitem{Coron07}
Jean-Michel Coron, \emph{Control and nonlinearity}, Mathematical Surveys and Monographs,
  vol. 136, American Mathematical Society, Providence, RI, 2007. \MR{2302744}

\bibitem{CBA08}
Jean-Michel Coron, Georges Bastin, and Brigitte d'Andr\'{e}a Novel,
  \emph{Dissipative boundary conditions for one-dimensional nonlinear
  hyperbolic systems}, SIAM J. Control Optim. \textbf{47} (2008), no.~3,
  1460--1498. \MR{2407024}

\bibitem{CoronGM18}
Jean-Michel Coron, Ludovick Gagnon, and Morgan Morancey, \emph{Rapid
  stabilization of a linearized bilinear 1-{D} {S}chr\"{o}dinger equation}, J.
  Math. Pures Appl. (9) \textbf{115} (2018), 24--73. \MR{3808341}

\bibitem{CHO17}
Jean-Michel Coron, Long Hu, and Guillaume Olive, \emph{Finite-time boundary
  stabilization of general linear hyperbolic balance laws via {F}redholm
  backstepping transformation}, Automatica J. IFAC \textbf{84} (2017), 95--100.
  \MR{3689872}

\bibitem{CHOS20}
Jean-Michel Coron, Long Hu, Guillaume Olive, and Peipei Shang, \emph{Boundary
  stabilization in finite time of one-dimensional linear hyperbolic balance
  laws with coefficients depending on time and space}, J. Differential
  Equations \textbf{271} (2021), 1109--1170. \MR{4160602}

\bibitem{Coron15}
Jean-Michel Coron and Qi~L\"{u}, \emph{Local rapid stabilization for a
  {K}orteweg-de {V}ries equation with a {N}eumann boundary control on the
  right}, J. Math. Pures Appl. (9) \textbf{102} (2014), no.~6, 1080--1120.
  \MR{3277436}

\bibitem{CoronNg17}
Jean-Michel Coron and Hoai-Minh Nguyen, \emph{{Null controllability and finite
  time stabilization for the heat equations with variable coefficients in space
  in one dimension via backstepping approach}}, Arch. Rational Mech. Anal.
  \textbf{225} (2017), 993--1023.

\bibitem{CoronNg19}
Jean-Michel Coron and Hoai-Minh Nguyen,  \emph{Optimal time for the controllability of linear hyperbolic
  systems in one-dimensional space}, SIAM J. Control Optim. \textbf{57} (2019),
  no.~2, 1127--1156. \MR{3932617}

\bibitem{CoronNg20}
Jean-Michel Coron and Hoai-Minh Nguyen,  \emph{Finite-time stabilization in optimal time of homogeneous
  quasilinear hyperbolic systems in one dimensional space}, ESAIM Control
  Optim. Calc. Var. \textbf{26} (2020), Paper No. 119, 24. \MR{4188825}

\bibitem{CoronNg20-L}
Jean-Michel Coron and Hoai-Minh Nguyen,   \emph{Lyapunov functions and finite time stabilization in optimal time
  for homogeneous linear and quasilinear hyperbolic systems},  (2020),
  https://arxiv.org/abs/2007.04104.

\bibitem{CoronNg19-2}
Jean-Michel Coron and Hoai-Minh Nguyen,  \emph{Null-controllability of linear hyperbolic systems in one
  dimensional space}, Systems Control Lett. \textbf{148} (2021), 104851.

\bibitem{CVKB13}
Jean-Michel Coron, Rafael Vazquez, Miroslav Krstic, and Georges Bastin,
  \emph{Local exponential {$H^2$} stabilization of a {$2\times 2$} quasilinear
  hyperbolic system using backstepping}, SIAM J. Control Optim. \textbf{51}
  (2013), no.~3, 2005--2035. \MR{3049647}

\bibitem{dHPCAB03}
Jonathan de~Halleux, Christophe Prieur, Jean-Michel Coron, Brigitte
  d'Andr\'{e}a Novel, and Georges Bastin, \emph{Boundary feedback control in
  networks of open channels}, Automatica J. IFAC \textbf{39} (2003), no.~8,
  1365--1376. \MR{2141681}

\bibitem{MVK13}
Florent Di~Meglio, Rafael Vazquez, and Miroslav Krstic, \emph{Stabilization of
  a system of {$n+1$} coupled first-order hyperbolic linear {PDE}s with a
  single boundary input}, IEEE Trans. Automat. Control \textbf{58} (2013),
  no.~12, 3097--3111. \MR{3152271}

\bibitem{DP08}
Val{\'e}rie Dos~Santos and Christophe Prieur, \emph{Boundary control of open
  channels with numerical and experimental validations}, IEEE transactions on
  Control systems technology \textbf{16} (2008), no.~6, 1252--1264.

\bibitem{Goatin06}
Paola Goatin, \emph{The {A}w-{R}ascle vehicular traffic flow model with phase
  transitions}, Math. Comput. Modelling \textbf{44} (2006), no.~3-4, 287--303.
  \MR{2239057}

\bibitem{GLR-Matrix}
Israel Gohberg, Peter Lancaster, and Leiba Rodman, \emph{Matrix polynomials},
  Classics in Applied Mathematics, vol.~58, Society for Industrial and Applied
  Mathematics (SIAM), Philadelphia, PA, 2009, Reprint of the 1982 original [
  MR0662418]. \MR{3396732}

\bibitem{GL84}
James~M. Greenberg and Ta~Tsien Li, \emph{The effect of boundary damping for
  the quasilinear wave equation}, J. Differential Equations \textbf{52} (1984),
  no.~1, 66--75. \MR{737964}

\bibitem{GL03}
Martin Gugat and G\"{u}nter Leugering, \emph{Global boundary controllability of
  the de {S}t. {V}enant equations between steady states}, Ann. Inst. H.
  Poincar\'{e} Anal. Non Lin\'{e}aire \textbf{20} (2003), no.~1, 1--11.
  \MR{1958159}

\bibitem{GLS04}
Martin Gugat, G\"{u}nter Leugering, and E.~J.~P. Georg~Schmidt, \emph{Global
  controllability between steady supercritical flows in channel networks},
  Math. Methods Appl. Sci. \textbf{27} (2004), no.~7, 781--802. \MR{2055319}

\bibitem{Hor97}
Lars H\"{o}rmander, \emph{On the uniqueness of the {C}auchy problem under
  partial analyticity assumptions}, Geometrical optics and related topics
  ({C}ortona, 1996), Progr. Nonlinear Differential Equations Appl., vol.~32,
  Birkh\"{a}user Boston, Boston, MA, 1997, pp.~179--219. \MR{2033496}

\bibitem{HMVK16}
Long Hu, Florent Di~Meglio, Rafael Vazquez, and Miroslav Krstic, \emph{Control
  of homodirectional and general heterodirectional linear coupled hyperbolic
  {PDE}s}, IEEE Trans. Automat. Control \textbf{61} (2016), no.~11, 3301--3314.
  \MR{3571452}

\bibitem{HO19}
Long Hu and Guillaume Olive, \emph{{Minimal time for the exact controllability
  of one-dimensional first-order linear hyperbolic systems by one-sided
  boundary controls}}, Preprint, hal-01982662.

\bibitem{Kato}
Tosio Kato, \emph{Perturbation theory for linear operators}, Classics in
  Mathematics, Springer-Verlag, Berlin, 1995, Reprint of the 1980 edition.
  \MR{1335452}

\bibitem{KGBS08}
Miroslav Krstic, Bao-Zhu Guo, Andras Balogh, and Andrey Smyshlyaev,
  \emph{Output-feedback stabilization of an unstable wave equation}, Automatica
  J. IFAC \textbf{44} (2008), no.~1, 63--74. \MR{2530469}

\bibitem{Krstic08}
Miroslav Krstic and Andrey Smyshlyaev, \emph{Boundary control of {PDE}s},
  Advances in Design and Control, vol.~16, Society for Industrial and Applied
  Mathematics (SIAM), Philadelphia, PA, 2008, A course on backstepping designs.
  \MR{2412038}

\bibitem{LL19}
Camille Laurent and Matthieu L\'{e}autaud, \emph{Quantitative unique
  continuation for operators with partially analytic coefficients.
  {A}pplication to approximate control for waves}, J. Eur. Math. Soc. (JEMS)
  \textbf{21} (2019), no.~4, 957--1069. \MR{3941459}

\bibitem{1981-Li-Wen-Chi-Shen-CAM}
Da~Qian Li, Wen~Chi Yu, and Wei~Xi Shen, \emph{Second initial-boundary value
  problems for quasilinear hyperbolic-parabolic coupled systems}, Chinese Ann.
  Math. \textbf{2} (1981), no.~1, 65--90. \MR{619464}

\bibitem{LY85}
Ta~Tsien Li and Wen~Ci Yu, \emph{Boundary value problems for quasilinear
  hyperbolic systems}, Duke University Mathematics Series, V, Duke University,
  Mathematics Department, Durham, NC, 1985. \MR{823237}

\bibitem{Lions88-VolI}
Jacques-Louis Lions, \emph{Contr\^{o}labilit\'{e} exacte, perturbations et
  stabilisation de syst\`emes distribu\'{e}s. {T}ome 1}, Recherches en
  Math\'{e}matiques Appliqu\'{e}es [Research in Applied Mathematics], vol.~8,
  Masson, Paris, 1988, Contr\^{o}labilit\'{e} exacte. [Exact controllability],
  With appendices by E. Zuazua, C. Bardos, G. Lebeau and J. Rauch. \MR{953547}

\bibitem{RT74}
Jeffrey Rauch and Michael Taylor, \emph{Exponential decay of solutions to
  hyperbolic equations in bounded domains}, Indiana Univ. Math. J. \textbf{24}
  (1974), 79--86. \MR{361461}

\bibitem{RZ98}
Luc Robbiano and Claude Zuily, \emph{Uniqueness in the {C}auchy problem for
  operators with partially holomorphic coefficients}, Invent. Math.
  \textbf{131} (1998), no.~3, 493--539. \MR{1614547}

\bibitem{Russell78}
David~L. Russell, \emph{Controllability and stabilizability theory for linear
  partial differential equations: recent progress and open questions}, SIAM
  Rev. \textbf{20} (1978), no.~4, 639--739. \MR{508380}

\bibitem{Slemrod83}
Marshall Slemrod, \emph{Boundary feedback stabilization for a quasilinear wave
  equation}, Control theory for distributed parameter systems and applications
  ({V}orau, 1982), Lect. Notes Control Inf. Sci., vol.~54, Springer, Berlin,
  1983, pp.~221--237. \MR{793047}

\bibitem{SCK10}
Andrey Smyshlyaev, Eduardo Cerpa, and Miroslav Krstic, \emph{Boundary
  stabilization of a 1-{D} wave equation with in-domain antidamping}, SIAM J.
  Control Optim. \textbf{48} (2010), no.~6, 4014--4031. \MR{2645471}

\bibitem{SK04}
Andrey Smyshlyaev and Miroslav Krstic, \emph{Closed-form boundary state
  feedbacks for a class of 1-{D} partial integro-differential equations}, IEEE
  Trans. Automat. Control \textbf{49} (2004), no.~12, 2185--2202. \MR{2106749}

\bibitem{SK05}
Andrey Smyshlyaev and Miroslav Krstic,  \emph{On control design for {PDE}s with space-dependent diffusivity or
  time-dependent reactivity}, Automatica J. IFAC \textbf{41} (2005), no.~9,
  1601--1608. \MR{2161123}

\bibitem{SK09}
Andrey Smyshlyaev and Miroslav Krstic, \emph{Boundary control of an anti-stable wave equation with
  anti-damping on the uncontrolled boundary}, Systems Control Lett. \textbf{58}
  (2009), no.~8, 617--623. \MR{2542119}

\bibitem{Sontag98}
Eduardo~D. Sontag, \emph{Mathematical control theory}, second ed., Texts in
  Applied Mathematics, vol.~6, Springer-Verlag, New York, 1998, Deterministic
  finite-dimensional systems. \MR{1640001}

\bibitem{Tataru95}
Daniel Tataru, \emph{Unique continuation for solutions to {PDE}'s; between
  {H}\"{o}rmander's theorem and {H}olmgren's theorem}, Comm. Partial
  Differential Equations \textbf{20} (1995), no.~5-6, 855--884. \MR{1326909}

\bibitem{Vazquez08}
Rafael Vazquez and Miroslav Krstic, \emph{Control of 1-{D} parabolic {PDE}s
  with {V}olterra nonlinearities. {I}. {D}esign}, Automatica J. IFAC
  \textbf{44} (2008), no.~11, 2778--2790. \MR{2527199}

\bibitem{XS02}
Cheng-Zhong Xu and Gauthier Sallet, \emph{Exponential stability and transfer
  functions of processes governed by symmetric hyperbolic systems}, ESAIM
  Control Optim. Calc. Var. \textbf{7} (2002), 421--442. \MR{1925036}

\bibitem{2002-Xu-Sallet-COCV}
Cheng-Zhong Xu and Gauthier Sallet, \emph{Exponential stability and transfer functions of processes
  governed by symmetric hyperbolic systems}, ESAIM Control Optim. Calc. Var.
  \textbf{7} (2002), 421--442. \MR{1925036}

\bibitem{2019-Zhang-preprint}
Christophe Zhang, \emph{Finite-time internal stabilization of a linear 1-{D}
  transport equation}, Systems Control Lett. \textbf{133} (2019), 104529, 8.
  \MR{4001127}

\end{thebibliography}
\end{document}